\documentclass[letterpaper,11pt,twoside]{article}
\usepackage{amsmath,amsfonts,amssymb,amsthm,mathabx,secdot}
\usepackage[utf8]{inputenc}
\usepackage[margin=2cm]{geometry}
\usepackage[scr=boondoxo]{mathalpha}

\usepackage{xcolor}

\newcommand{\mytitle}{ Persistence of hyperbolic solutions of
  ODE's under functional perturbations: Applications
  to  the  motion of relativistic charged particles 

}

\newcommand{\myshortitle}{%
New Solutions from Uniformly Hyperbolic Trajectories under Functional Perturbations%
}

\usepackage{graphicx}
\graphicspath{{images/}}
\usepackage[hidelinks]{hyperref}
\hypersetup{
bookmarksnumbered=true,
bookmarksopen=true,
pdftitle={\myshortitle{}},
pdfauthor={J. Gimeno, J. Yang, and R. de la Llave},
}
\usepackage{orcidlink}

\usepackage{titlecaps}
\Addlcwords{is with of for on in are and its to the by}
\usepackage{comment} 
\newtheorem{thm}{Theorem}[section]
\newtheorem{lem}[thm]{Lemma}
\newtheorem{cor}[thm]{Corollary}
\newtheorem{pro}[thm]{Proposition}

\theoremstyle{remark}
\newtheorem{rmk}[thm]{Remark}
\theoremstyle{definition}
\newtheorem{dfn}[thm]{Definition}

\newcommand{\R}{\mathbb{R}}
\newcommand{\N}{\mathbb{N}}

\newcommand{\sC}{\mathscr{R}}

\newcommand{\X}{\mathcal{X}}

\newcommand{\dop}{\mathrm{D}} 

\newcommand{\Id}{Id}

\newcommand{\prm}{\mu}

\newcommand{\Et}{E}
\newcommand{\Pop}{\mathscr{P}}
\newcommand{\ep}{\varepsilon}

\newcommand{\ball}{\mathscr{B}}

\newcommand{\vf}{f} 
\newcommand{\pp}{P} 
\newcommand{\rr}{r} 

\newcommand{\op}{\Gamma ^\ep}  

\newcommand{\ee}{\mathcal E}

\newcommand{\xx}{x}
\newcommand{\xo}{\xx _0}
\newcommand{\hx}{\widehat\xx}   

\newcommand{\yy}{y}

\newcommand{\hy}{\widehat\yy}

\newcommand{\hatt}{\phi}
\newcommand{\vft}{X}
\newcommand{\hvft}{\vft}

\newcommand{\hatty}{\psi}
\newcommand{\vfty}{Y}

\newcommand{\B}{\mathcal{B}}

\newcommand{\bt}{\mathfrak{t}}
\newcommand{\bs}{\mathfrak{s}}
\newcommand{\bu}{\mathfrak{u}}
\newcommand{\ct}[1]{\mathfrak{#1}}

\DeclareMathOperator{\Span}{Span}
\DeclareMathOperator{\Lip}{Lip}

\newcommand{\bydef}{\,\stackrel{\mbox{\tiny\textnormal{\raisebox{0ex}[0ex][0ex]{def}}}}{=}\,}

\title{\mytitle{}}

\author{
  Joan Gimeno\,\orcidlink{0000-0002-8707-6379}\textsuperscript{(1)} \and 
  Rafael de la Llave\,\orcidlink{0000-0002-0286-6233}\textsuperscript{(2)} \and 
  Jiaqi Yang\,\orcidlink{0000-0002-5282-2722}\textsuperscript{(3,\textasteriskcentered)}
  }

\begin{document}

\maketitle
{\small
\begin{itemize}
\renewcommand{\itemsep}{0pt}
\item [(1)] Departament de Matem\`atiques i Inform\`atica, %
Universitat de Barcelona, %
Gran Via de les Corts Catalanes, 585, 08007 Barcelona, %
Spain, \verb+joan@maia.ub.es+
\item [(2)] School of Mathematics (Emeritus) %
Georgia Institute of Technology, %
686 Cherry St., Atlanta GA. 30332-0160, %
USA, \verb+rafael.delallave@math.gatech.edu+
\item [(3)] School of Mathematical Sciences,
Key Laboratory of Intelligent Computing and Applications (Ministry of Education), 
Tongji University, Shanghai 200092, China, \verb+jqyang@tongji.edu.cn+
\item[\textsuperscript\textasteriskcentered] Corresponding Author
\end{itemize}
}

\begin{abstract}

We rigorously construct a variety of orbits for certain delay differential equations, including the electrodynamic equations formulated by Wheeler and Feynman in 1949. These equations involve delays and advances that depend on the trajectory itself, making it unclear how to formulate them as evolution equations in a conventional phase space. Despite their fundamental significance in physics, their mathematical treatment remains limited.

Our method applies broadly to various functional differential equations that have appeared in the literature, including advanced/delayed equations, neutral or state-dependent delay equations, and nested delay equations, under appropriate regularity assumptions.

Rather than addressing the notoriously difficult problem of proving the existence of solutions for all the initial conditions in a set, we focus on the direct construction of a diverse collection of solutions. This approach is often sufficient to describe physical phenomena. For instance, in certain models, we establish the existence of families of solutions exhibiting symbolic dynamics.

Our method is based on the assumption that the system is, in a weak sense, close to an ordinary differential equation (ODE) with ``hyperbolic'' solutions as defined in dynamical systems. We then derive functional equations to obtain space-time corrections.

As a byproduct  of the method,
we obtain that
the solutions constructed depend very smoothly on parameters of the
model. Also, we show that  many formal approximations
currently used in physics are  valid with explicit error terms.
Several of the relations between different orbits of
the ODE persist qualitatively  in the full problem.

\end{abstract}

{\footnotesize
\smallskip

{\em 2010 Mathematics Subject Classification.} %
  34K05
, 34K19
, 34K13
, 34D15
, 34K26

{\bf Keywords:} Functional Differential Equations; %
Uniformly Hyperbolic Orbit; %
Perturbation Theory; Electrodynamics; %
Neutral Equations; %
Electrodynamics %
}

\pagestyle{myheadings}
\markboth{\myshortitle{}}{J. Gimeno, R. de la Llave, and J. Yang}
\clearpage
\phantomsection\pdfbookmark[1]{\contentsname}{toc}
\tableofcontents
\newpage

\section{Introduction}
We consider 
the problem of finding trajectories $x:\R \rightarrow \R^n$
solving an equation of the form:
\begin{equation} \label{model}
  \dot x(t) = \vf \circ x(t) + \ep \pp(t, x_t, \ep, \prm),
\end{equation}
where $\vf \colon \R ^n \to \R ^n$ is a smooth function (possibly only defined
 on an open subset in $\R^n$), $\ep$ a perturbative parameter, $\prm$ an additional 
parameter, $\pp \colon \R \times \sC[-h,h] \times (0,1)^2 \to \R^n$  is a smooth map. Here,  $\sC$ is a space of differentiable functions which
will be specified later.
The notation  $x _t$ is a ``history segment" of size $h \geq 0$ of the solution which is defined as:
\begin{equation} \label{eq.xt}
 x_t(s) \bydef x (t + s), \qquad s \in [-h,h].
\end{equation}

The term $\pp$ changes the nature of \eqref{model} for $\ep \ne 0$.
Hence, in spite of having the small parameter $\ep$,
the perturbation is not small.
Even for common delays
\[
\dot x(t) = \vf \circ x(t) + \ep x(t - 1),
\]
the natural phase space for $\ep \ne 0$ is infinite
dimensional, but for $\ep = 0$ it is just $\R^n$. Our result can also deal
with equations like $\dot x(t) = \vf \circ x(t) + \ep {\dot x}(t - 1)$ which makes the unboundedness even more apparent.

\smallskip

We will take solutions of \eqref{model}
for $\ep = 0$ which are \emph{``hyperbolic''}
in the sense of dynamical systems
(see Definition~\ref{dfn.unihypersol}) 
and show that one can find solutions that resemble
them after corrections (both in the positions occupied
and the speed of travel) when $\ep\neq 0$.
This allows to show that the equations~\eqref{model}
contain sets of solutions that support symbolic dynamics.

In Section~\ref{sec.main},
we formulate the precise regularity  assumptions on $\pp$.
We anticipate that, roughly, the main assumptions are
that applying $\pp$ to functions with $\ell+1$ derivatives
produce functions with $\ell$ derivatives and that there
are some Lipschitz bounds in $C^0$  when the arguments 
lie in  spaces of smooth
functions. 

A case that served as a motivation for us is the electrodynamics
of point charges. See Section~\ref{sec:electrodynamics}.
In the model of \cite{WheelerF49},
the particles move in the Li\'enard-Wiechert
potentials (relativistic analogues of Coulomb-Ampere formulas)
generated by the other particles.
This leads to advanced/delayed equations with several delays
which are obtained solving implicit equations that
involve the trajectories. Given the physical importance of
these equations, there have been several results
establishing existence of solutions in the literature, mainly in the
one dimensional case \cite{HoagD90, BauerDDH17, DeckertH16}.  Other models
of electrodynamics (notably several versions of
the Post-Newtonian formalism) can be accommodated. The results
here provide existence of various solutions for all of
them and allow to discuss how approximate are the solutions.
We note that the effect of the delays/advances are formally of
size inversely to the speed of light, $1/c$, which is much larger than radiation effects 
that are of size $(1/c)^3$.

\smallskip

Our results also cover other cases in the literature
in which the
perturbation is clearly singular. Including for an arbitrary $\vartheta\in\sC[-h,h]$, 
\begin{itemize}
 \item $\pp(t, \vartheta, \ep, \prm) = g \circ \vartheta(0) $, an
   ODE perturbation as $\pp\bigl(t, \xx _t, \ep, \prm\bigr) =
   g\circ \xx(t)$;
 \item $\pp(t, \vartheta, \ep, \prm) = \vartheta(-1)$, a perturbation with a
   constant delay as $\pp\bigl(t, \xx _t, \ep, \prm \bigr)
   = \xx(t-1)$;
 \item $\pp(t, \vartheta, \ep, \prm) = \frac{1}{\ep} \vf \circ
   \vartheta(-\ep) - \frac{1}{\ep} \vf \circ\vartheta(0)$, the
   small delay system  $\dot \xx(t) = \vf \circ \xx (t - \ep)$;
 \item $\pp(t, \vartheta, \ep, \prm) = \vartheta \circ \rr \circ
   \vartheta(0)$, a state-dependent delay perturbation as $\pp\bigl(t,
   \xx _t, \ep, \prm\bigr) = x(t + \rr \circ x(t))$;
 \item $\pp(t, \vartheta, \ep, \prm) = \vartheta \circ \rr \circ
   \vartheta \circ \rr _1 \circ \vartheta(0)$ containing nested delays
   as $\pp\bigl(t, \xx _t, \ep, \prm\bigr) = x(t + r \circ x(t +
   \rr _1 \circ x(t)))$;
 \item $\pp(t, \vartheta, \ep, \prm) = (\frac{d}{ds}\vartheta)(0)$, an implicitly defined ODE as
   $\pp\bigl(t, \xx _t, \ep, \prm   \bigr) = \dot \xx(t)$;
 \item $\pp(t, \vartheta, \ep, \prm) = (\frac{d}{ds}\vartheta)(-1)$, a neutral equation with a constant delay as $\pp\bigl(t,
   \xx _t, \ep, \prm \bigr) = \dot \xx(t-1)$;
 \item $\pp(t, \vartheta, \ep, \prm) = \vartheta \circ \rr \circ
   (\frac{d}{ds}\vartheta)(0)$, a first
   order neutral equation as $\pp\bigl(t, \xx _t, \ep, \prm\bigr) = \xx(t + \rr
   \circ \dot x(t))$;
 \item $\pp(t, \vartheta, \ep, \prm) = \vartheta\circ \tau(t)$
   containing an explicit time-dependent delay as $\pp\bigl(t,
   \xx _t, \ep, \prm \bigr) = \xx(t +\tau (t))$;
   
 \item $\pp(t, \vartheta, \ep, \prm) $ obtained by solving an
   implicit equation of $\vartheta$, e.g. \eqref{delaydefined} in electrodynamics;
   
 \item a bounded time-dependent map $\pp$, so that we have a non-autonomous perturbation;
 
 \item $\pp(t, \vartheta, \ep, \prm) = \vartheta(-1)+\vartheta(-2) + \vartheta(+1) $, a perturbation containing several delays or advances. 
\end{itemize}

\smallskip

In the theory of delay differential equations (DDE's) with constant delay, it is customary
to take $\sC$ to be the space of continuous functions, whereas here we
 find it useful to consider spaces of more differentiable
functions so that the functional $\pp$ can involve first
derivatives (neutral equations) or more complicated forms.  Note that
we are considering ``history segments" that include both the past and
the future of the trajectory so that our theory works just as well for
delayed, advanced, or mixed expressions. The value $h$ will by default
belong to $[0,+\infty)$ even if some existence results may
  hold  for $h=\infty$. 
The assumption $h < \infty$ seems to be essential for the
  local uniqueness and for the a-posteriori results.

\subsection{Informal main results} \label{sec.informal}
Our results show that, if  equation \eqref{model} admits a set
consisting of uniformly hyperbolic solutions (see
Definition~\ref{dfn.unihypersol}) when $\ep=0$, under appropriate regularity
assumptions on the functional $\pp$, equation \eqref{model} admits a
set of solutions which are close to the hyperbolic solutions of the ODE as
long as $\ep$ is small enough. Our result is similar to
structural stability interpreted in the functional analysis
formulation of the perturbed problem \eqref{model}, see
Theorem~\ref{main}. We stress that our approach does not need to discuss the
phase space of solutions to solve any possible initial value problem of
\eqref{model}, what we do is to search for solutions with a specific
structure. In this work we look for solutions of a form based on the uniformly hyperbolic solution of the ODE. We consider a space of
functions $x(t) $ of
a specific form and formulate equations which imply that such
$x$ satisfies \eqref{model}.

The above  strategy bypasses the study of general existence, uniqueness
and dependence on initial conditions of the solutions.
The solutions of \eqref{model} we construct could fail
to be surrounded by other solutions. This strategy was used
already in \cite{HR, HR2, HeY20, YangGL21, YangGL22}.
We also note that the equations considered are numerically
well conditioned and can be implemented to produce approximate
solutions. 

\smallskip

The existence of hyperbolic sets (collections of hyperbolic orbits)
in differential
equations is a rather common situation.  Notably, existence of a 
transverse homoclinic intersections implies the existence of a  uniformly hyperbolic set (the horseshoe) which has 
a very rich dynamics including an uncountable set 
of hyperbolic orbits described by symbolic dynamics.
Other famous attractors (Lorenz, R\"ossler,
Chua, \ldots) have also been documented.
These attractors include uniformly hyperbolic sets. 
Our  results imply that all these uniformly hyperbolic sets 
persist   when we add a perturbation with a sufficiently small parameter. The perturbations allowed are very general
and include perturbations that are singular from a conventional point of view. See an informal presentation in Theorem~\ref{thm.informal}. A precise formulation 
is in Theorem~\ref{main}.

We call attention to 
\cite{Wayda95} which uses Poincare returns to 
establish persistence of hyperbolic sets 
in $C^1$ perturbations of Functional Differential Equations  that generate  a $C^1$ evolution.  
In \cite{Lani-WaydaW16},  existence of chaotic motion was established for an SDDE analyzing the evolution and finding an
analogue of Shilnikov phenomenon. In \cite{WojcikZ05}, the authors focused on small constant delay perturbations of an ODE and obtained persistence of topological horseshoes. In contrast, 
the present method is not based on analysis of evolution and
applies to advanced/delayed equations that do not define 
any evolution. Indeed, we do not need to study regularity properties 
of the evolution and not even the space where evolution is 
defined (which may involve the study of solution manifolds). 

We also call attention to the papers \cite{HumpheisBCH16} and \cite{CallejaHK17}. These papers
use numerics and bifurcation analysis to study singularly perturbed state-dependent delay equations,
where they also find that state-dependence of the delays can generate very complex dynamics.
We think it would be interesting to  reformulate our fixed point problems so that they could validate
the new solutions found  and, specially the bifurcation point (where hyperbolicity is lost).

We are not assuming many properties of the hyperbolic sets beyond
requiring that the hyperbolicity constants are uniform (sometimes
called Pesin sets \cite{BarreiraP} in non-uniform hyperbolic theory).
The hyperbolic sets we considered
may fail to be closed or \emph{maximal}.

The proofs are rater explicit and the sizes of the perturbations allowed are computable in concrete examples. (Some related calculations were achieved in \cite{GimenoLMY23}.)
\smallskip

One downside of the program presented here is that, by design, we
cannot discuss properties of the evolution for all initial data.
Nevertheless, one could remark that, even in the qualitative theory of
ODE's one often relies only on developing landmarks that organize the
behavior of all the solutions.  The analogue of the qualitative theory
in our program would be to develop a theory indicating that solutions
of a certain kind imply the existence of others. In that respect, we
hope to come back to the study of stable manifolds and the study of existence of 
symbolic dynamics from some finite calculations.

Another downside of the present treatment is that we are constrained
by the regime of solutions close to the solutions of the ODE's.  It is
well known that many of the equations we study will have many
solutions which do not resemble the solutions of the ODE.
Nevertheless, in physical applications (e.g. motions of charged
particles) the delays are small so that the effects of the delay are
small and hard to observe. (This is why relativity was only discovered
in the XX century).

Our results will apply to the space of functions that are finitely
differentiable with finite norm in a segment, namely $[-h,h]$, where
$h > 0$ is the domain of the ``history segment''.  Notice that if we
have a finite differentiable function $x \colon [-\tilde h,\tilde h]
\to \R ^n$ in a slightly bigger history segment, $\tilde h > h$, then
the function $x _t$ in \eqref{eq.xt} is defined for all
$t \in[-(\tilde h - h), \tilde h - h]$. Therefore, we can apply functionals
to all $x _t$ in an open interval of $t$.

\medskip

To give in a glimpse of the precise main result on this paper (see
Theorem~\ref{main}), let us first provide an informal result omitting
many precise formalism: The result shows that under a mild set of
hypotheses on the perturbation, the system \eqref{model} has
solutions nearby the hyperbolic orbit of the unperturbed ODE and such solution will be unique
in a suitable neighborhood.

\begin{thm}[informal result] \label{thm.informal}
 We consider perturbation of an ODE as in \eqref{model} and let $\ell
 \geq 0$ be an integer. Assume that:
 \begin{enumerate}
  \item The unperturbed ODE admits a solution $\{\xo(t)\}_{t\in \R}$ which is
    uniformly hyperbolic, see Definition~\ref{dfn.unihypersol}.
  \item The function $\vf$ is uniformly bounded as well as its
    derivatives up to order $\ell+3$ in a $\delta$-neighborhood of
    the orbit $\{\xo(t)\}$.
  \item The perturbation functional $\pp$ in \eqref{model} satisfies 
    ``propagated bounds''
    (i.e. when $\xx_t$ ranges in a ball in a space of
    $C^{\ell + 2 }$ functions, $\pp$ lies in a ball of $C^{\ell+1}$
    functions).
  \item The functional $\pp$ is Lipschitz in a low regularity,
    i.e. for all $u $ and $v$ in
    a $C^{\ell + 2 }$ ball and all $t$ and $s$, there are constants $B_1, B_2$ such that
    \[
    | \pp(t, u,  \mu) - \pp(s, v , \mu) |
    \le B_1 |t - s| + B_2 \| u - v \|_{C^1}.
    \]
 \end{enumerate}
 Then there is an $\ep _0 > 0$ such that for all $|\ep| \leq \ep _0$,
 there exist differentiable maps $\hx ^\ep \colon \R \to \R ^n$ and
 $\hatt ^\ep \colon \R \to \R$ such that $\hx ^\ep$ is in $C
 ^{\ell+1}$, $\dop \hatt^\ep$ is in $ C ^{\ell}$, and
 \begin{equation} \label{solutionformep} 
  \xx(t) = (\xo + \hx ^\ep) \circ \hatt ^\ep(t)
 \end{equation}
 is a $C ^{\ell+1}$ solution of \eqref{model}. Moreover, if
 \eqref{model} depends smoothly on parameters in an appropriate sense,
 so does the new solution $\xx$.
\end{thm}

The informal Theorem~\ref{thm.informal} is slightly different from the
formal (see Theorem~\ref{main}). We simplified
the statement to avoid introducing $C^{\ell + \Lip}$ spaces
and we did not include some technical
considerations.  All of these are discussed
in Theorem~\ref{main} along with the a-posteriori
formulation.  We also omitted several properties
that are used to prove  smooth dependence on parameters, see Section \ref{sec.smoothprm}.

\begin{rmk}[on perturbative regularity loose]
Note that we are allowing that the functional $\pp$ appearing as a
perturbation \emph{looses one derivative}, which means $\pp$ may have one derivative less than its second argument.
This will be key to apply the result to neutral differential equations
and to equations with small delay. 
\end{rmk}

\begin{rmk}[on unknown corrections]
We have two unknowns  in our existence problem: $\hx$ acting
as an additive space correction term and $\hatt$ as a time reparametrization
correction term. Both of them depend on the perturbative
parameter $\ep$ and their regularity properties are derived from a
fixed point scheme. Indeed, we will write an equation where
$x(t)$ in Theorem~\ref{thm.informal} is a solution and arrive at a fixed point problem after
manipulation.
\end{rmk}

\begin{rmk}[on further conclusions]
As a consequence of the fixed point method, we are able to
formulate and prove that the solutions we build in the perturbative
system depend smoothly on parameters. The parameters can be in the
unperturbed ODE or in the perturbative map. Moreover, after a mild change of hypotheses, we will also see how the solution can admit an
exponential derivative growth, see Section~\ref{sec.exposol}.
\end{rmk}

\begin{rmk}[on interesting applications]\label{rmk:anosov}
Among the applications of our theorem, a particular case  is 
when the perturbed equation is another ODE. 
The results on persistence
of such solutions were studied originally in 
\cite{Anosov69} and later in \cite{Moser69}.
We, however, obtain smooth dependence on parameters,
which is not true for the formulation in the above references.
The formulation we use is slightly different and
follows more closely the formulation in  \cite[Appendix A]{LlaveMM86},
which also obtained smooth dependence on parameters
for the objects considered (the objects considered
in \cite{LlaveMM86} are roughly, inverses of
the objects considered in \cite{Anosov69, Moser69}). 
\end{rmk}

\begin{rmk}[on details from previous works]
In previous papers \cite{YangGL21,YangGL22}, we formulated the results
in an a-posteriori format, meaning that if we start with an initial
guess of the correction whose error is small enough, then
the theorems conclude that there is a true solution nearby. 
The a-posteriori formulation is suitable for performing computer-assisted proofs
(the approximate solution is produced by a numerical calculation
and the needed estimates are verified using a computer
by taking care of truncation and round-off error) \cite{GimenoLMY23}. Here the a-posteriori formulation
will be more delicate. Indeed, we will use a
different norm  to obtain contractions
and the a-posteriori argument will be valid
on segments of times $t$, see Section~\ref{sec.aposteriori}.

  This paper has similarities with \cite{YangGL22} in that we
  seek  both an embedding and an inner dynamics.
  However, this paper is significantly more difficult.

  The main reason is that, given any vector field without zeros in the circle,
  there is a change of variables  that reduce it to a constant, so
  that,  in \cite{YangGL22} the inner dynamics was just a number.
  In the present case, vector fields in the line cannot,
  in general, be reduced to constants (or even be approximated
  well by periodic; take for example, vector fields who oscillate
  between two values over longer and longer intervals).

  Hence in the present case, rather than dealing with just a number,
  we have to deal with an infinite dimensional unknown that, furthermore
  appears in the functional equations as a composition on the right.

  If we apply the formalism in this paper to the case of periodic
  solutions, we will obtain a periodic vector field $\vft$, not
  necessarily constant as it happens when one applies the formalism of
  \cite{YangGL22}. Furthermore, the formalism in \cite{YangGL22} does
  not satisfy the normalizations \eqref{normalization2}. 
\end{rmk}

\subsection{Organization of the Paper}
\label{sec:organization} 

In Section~\ref{sec.unperassum},  we describe
precisely the assumptions on the unperturbed system. The main
part is the (rather standard)  definition of hyperbolic orbits,
which we use to set the notation. We also present
a characterization of the  invariant bundles.

Section~\ref{sec.constsol} presents the formalism we use to describe
the solutions 
of  perturbed system. We present class of functions we will
consider, and perform manipulations to derive a functional
equation (called invariance equation, see \eqref{Gammac}, \eqref{Gammasu}
whose solutions give solutions of  \eqref{model} when substituted
in \eqref{solutionformep}). 
As it turns out,
this invariance  equation has symmetries under changes of variables and
we  also present normalization conditions that lead to  unique solution.

Section~\ref{sec.formulation} is devoted 
to completing the formulation of the fixed point problem, specifying
the fixed point operator, its domain and range,
and the rigorous formulation of the main result, Theorem~\ref{main}.

In Section~\ref{sec.proofstrategy}, we present  the proof of Theorem~\ref{main}
starting with an overview of the strategy. 
Section~\ref{sec.further} contains results not explicitly covered in
the main formulation which require
cumbersome notations.

Finally, Section~\ref{sec.examples} provides
examples of physical interest where our results apply. In
particular, it includes the case in which the delay is small
and the motions of charged
particles with electromagnetic interactions.

\section{The unperturbed ODE} \label{sec.unperassum}
When $\ep =0$, the unperturbed system \eqref{model} is an autonomous
ODE. We assume that such a system has a uniformly hyperbolic solution
$\{\xo(t)\}_{t\in \R}$, see Definition~\ref{dfn.unihypersol}.

\medskip

First, we recall that the variational equation (or equation of
variation) of the unperturbed ODE \eqref{model}, around a solution
$\{\xo(t)\}_{t\in \R}$, is the time-dependent linear equation,
\begin{equation} \label{eq.variation}
 \dot \xi(t) = \dop \vf \circ \xo(t) \xi(t),
\end{equation}
where $\xi(t) \in \R ^n$ has the heuristic meaning of small deviations
from the baseline trajectory.  The linear equation
\eqref{eq.variation} has a family of fundamental matrix solutions
$\{U(v;t)\} _{v,t\in \R}$ such that:
\begin{equation} \label{eq.fundsol}
 \frac{d}{dv} U(v;t) = \dop \vf \circ \xo(v) U(v;t), \qquad U(t;t) = \Id _n,
\end{equation}
where $\Id _n$ denotes the $n\times n$ identity matrix.
Note that due to the existence and uniqueness of
the variational equations,
we have
\begin{equation}\label{cocycle}
  U(v;t) = U(v, s) U(s;t).
\end{equation} 

\subsection{Uniformly hyperbolic solutions  of an ODE
  and their  quality measures}

In this section, we present Definition~\ref{dfn.unihypersol}.
The starting point of Theorem~\ref{main} is
precisely that we have a solution of the unperturbed
problem satisfying Definition~\ref{dfn.unihypersol}. 
We note that this definition has qualitative aspects called \emph{``quality measures of the
hyperbolicity''}. The ranges of perturbation parameters
that are allowed depend on the values of these numbers. 

It is well known that these uniformly hyperbolic orbits often appear
together in hyperbolic sets (e.g. horseshoes, Lorenz attractor,
etc.) but the quality measures may deteriorate
as we consider orbits in the attractors. This is also common
in the theory of non-uniformly hyperbolic sets. 

\begin{dfn}[Uniformly hyperbolic solution of an ODE] \label{dfn.unihypersol}
 Let $\{\xo(t)\}_{t \in \R}$ be a solution of the unperturbed system
 \eqref{model}.  We say that $\{\xo(t)\}_{t \in \R}$ is
 \emph{uniformly hyperbolic} if, and only if, it satisfies:
 \begin{enumerate}
  \renewcommand{\theenumi}{\roman{enumi}}
  \renewcommand{\labelenumi}{\theenumi.)}
  \item For each $t \in \R$, there exists a decomposition of the
    tangent space at $\xo(t)$,
  \begin{equation} \label{eq.splitting}
   \R ^n \cong T _{\xo(t)} \R ^n = \Et _t^c \oplus \Et _t^s \oplus \Et _t^u,
  \end{equation}
  such that 
  \begin{enumerate}
  \renewcommand{\theenumii}{\alph{enumii}}
  \renewcommand{\labelenumii}{\theenumi.\theenumii)}
   \item $\Et _t ^c = \Span\{\vf \circ \xo(t)\}$ has dimension $n _c =
     1$.
   \item $\Et _t ^s$ and $\Et _t^u$ have dimensions $n _s$ and $n _u$
     respectively. Thus, $n = 1 + n _s + n _u$.
    \item $\Et _t ^\sigma$ depends on $t$ continuously for $\sigma \in \{c,s,u\}$.
  \end{enumerate}
  Moreover, the forward (resp. backward) semiflow of the variational
  equation is contractive on $\Et _t ^s$ (resp. $\Et _t ^u$).  More
  precisely, the fundamental matrices $\{U(v;t)\}_{v,t\in \R}$
  in \eqref{eq.fundsol} admit center, stable, and unstable
  families of linear operators
  \begin{equation*}
   \begin{split}
    \{U ^c(v;t)\}_{v,t\in \R}, \quad U ^c(v;t) \colon \Et _t ^c& \to \Et _v ^c, \\
    \{U ^s(v;t)\}_{v,t\in \R}, \quad U ^s(v;t) \colon \Et _t ^s& \to \Et _v ^s, \\
    \{U ^u(v;t)\}_{v,t\in \R}, \quad U ^u(v;t) \colon \Et _t ^u& \to \Et _v ^u ,
   \end{split}
  \end{equation*}
  satisfying for $\sigma \in \{c,s,u\}$,
  \begin{equation} \label{eq.variEt}
    \frac{d}{dv} U^\sigma(v;t) = \dop \vf \circ \xo(v) U^\sigma(v;t),
    \qquad U^\sigma(t;t) = \Id|_{\Et _t ^\sigma},
  \end{equation}
  where $\Id|_{\Et _t ^\sigma}$ denotes the identity operator
  restricted to the linear subspace $\Et _t ^\sigma$.
  \item There exist $C _U, \lambda _s, \lambda _u > 0$ such that 
  \begin{equation} \label{eq.expoU}
   \begin{aligned}
    | U ^s (v;t) | &\leq e ^{-\lambda _s(v - t)} C _U & v &\geq t, \\
    | U ^u (v;t) | &\leq e ^{ \lambda _u(v - t)} C _U & v &\leq t,
   \end{aligned}
  \end{equation}
  where $|\cdot |$ is the operator norm.
  \item There exists $C _\Pi > 0$ such that 
  \begin{equation}
   \sup _{t \in \R} \| \Pi _t^\sigma \| \leq C _\Pi, \qquad \sigma \in \{c,s,u\},
  \end{equation}
  where $\Pi _t ^\sigma \colon \R ^n \to \Et _t ^\sigma$ denotes the
  projection corresponding to the splitting \eqref{eq.splitting}.
 \end{enumerate}
\end{dfn}

\medskip

 In particular, for the center direction projection, for any given vector $V
 \in T_{\xo(t)}\R ^n$, there exists $A_V \in \R$
 \[ \Pi ^c _t V = A_V \vf \circ \xo(t). \]

\begin{dfn}[Quality measures of Uniformly Hyperbolic orbit]\label{dfn.qual}
 The quantities $C _U$, $C _\Pi$, $\lambda _s$, and $\lambda _u$
 appearing in Definition~\ref{dfn.unihypersol} are referred to as {\em
   quality measures} of the hyperbolic solution $\{\xo(t)\}_{t\in\R}$.
\end{dfn}

\begin{rmk}[on the variational on the $\Et_t^\sigma$]
 The ODE's in \eqref{eq.variEt} are understood as equations for
 operators on $\R ^n$. In particular, the restriction of $\Id$ on $E
 ^\sigma _t$ must be interpreted using the injection from $E ^\sigma
 _t$ to $\R^n$ for $\sigma \in \{c,s,u\}$.
\end{rmk}

\begin{rmk}[on the projections] \label{rmk.projections}
 The projections $\Pi ^\sigma_t$ in Definition~\ref{dfn.unihypersol}
 may not be orthogonal projections on the space and they depend on the
 decomposition, e.g., $\Pi ^s _t$ could change if $\Et _t^u$ changes
 even if
 $\Et _t^s$  remains fixed.

 The constant $C _\Pi$ can be interpreted as the inverse of a measure
 of the angles between the spaces in the decomposition
 \eqref{eq.splitting}.
\end{rmk}

\begin{rmk}[on hyperbolic regularity]
  We have formulated Definition~\ref{dfn.unihypersol} including
  continuous dependence of the bundles on the base point along the
  orbit to keep compatibility with the standard definitions of
  normally hyperbolic manifolds in \cite{Fenichel71, HirschPS77}.

  These references show that, when the vector field is $C ^{r}$ and
  bounded and the hyperbolicity is uniform, then the continuous
  splittings are actually $C^{r-1}$.  We will provide the details in
  the formal result section. 
\end{rmk}

\begin{rmk}[on the quality measures]
  Note that the quality measures $C _\Pi$ and $C _U$ depend on the
  metric used.  In the theoretical literature on hyperbolic systems,
  it is standard to define a metric (and modify slightly the exponents
  of contraction) called \emph{adapted metric} so that $C _\Pi = C _U
  = 1$ and the center, stable, and unstable directions are
  orthogonal. This adapted metric is equivalent to the original one.
  Some rigorous proofs get simplified by using the adapted metric.
  
  Nevertheless, we do not use an adapted metric in this work for
  several reasons. The use of adapted metric would be confusing for us
  since the formulas for state-dependent delays are affected by the
  metric as well.  Moreover, the strength of the perturbations allowed
  in this paper depends on the values of the quality measures and the
  sizes of the derivatives of the perturbations. Changing the metric
  would require measuring the properties of the perturbing function in
  the adapted metric.

  Besides, the use of an adapted metric also obscures the study of
  phenomena that happen in the boundary of hyperbolicity. Notably
  \cite{HaroL06, HaroL07} identified numerically a boundary of
  hyperbolicity characterized by $C_\Pi $ blowing up (the angle
  between the splittings going to zero) while the exponents of
  contraction remain uniformly bounded away from zero.
\end{rmk}

\subsection{Infinitesimal characterization of the invariant bundles of a hyperbolic orbit} 

It will be useful for us  to characterize trajectories in $\Et _t ^\sigma$ for
$\sigma \in \{c,s,u\}$ as solutions of ODE's.  In
Lemma~\ref{lem:caracinv},   we use $U^\sigma(0;t)$ to connect $\xi(t) \in \Et
_t^\sigma$ with $\Et _0 ^\sigma$. This gives an infinitesimal
characterization of the invariant bundles.

\begin{lem} [Bundle Characterization]\label{lem:caracinv}
 Let $\{\xo(t)\}_{t\in\R}$ be a uniformly hyperbolic orbit and let
 $\sigma \in \{s,c,u\}$. If $\vf$ is differentiable enough (at least
 $C ^1$), then $\xi(t) \in \Et ^\sigma _t$ if, and only if,
 \begin{enumerate}
  \item $\xi (0) \in \Et ^\sigma _0$; and 
  \item $\dot \xi(t) = \dop \vf \circ \xo(t) \xi(t) + a(t)$ with $a(t)
    \in \Et ^\sigma _t$.
 \end{enumerate}
\end{lem}
 \begin{proof}
 \begin{enumerate}
  \item [$\Rightarrow)$] Let $\alpha(t) = U ^\sigma (0; t) \xi(t)$ or,
    equivalently, $\xi(t) = U ^\sigma(t;0) \alpha (t)$. Note that
    $\alpha (t) \in \Et ^\sigma _0$ for all $t$ and consequently $\dot
    \alpha(t)$ belongs to $\Et^\sigma _0$ as well. By taking
    derivatives we obtain
  \[
  \dot \xi(t) = \dop \vf \circ \xo (t) \xi(t) + U ^\sigma (t; 0) \dot \alpha(t).
  \]
  Thus, $a(t)\bydef U^\sigma(t;0) \dot \alpha(t) \in \Et ^\sigma _t$.
\item [$\Leftarrow)$] By the variation of parameters formula, 
  \[
  \xi(t) = U ^\sigma (t; 0)  \xi (0) + \int _0 ^t U ^\sigma (t;s) a(s)\, ds.
  \]
  And by the invariance of the $\sigma$-space; i.e. $U ^\sigma (t;s) \Et
  ^\sigma _s = \Et ^\sigma _t$, then $\xi(t) \in \Et ^\sigma_t$.
  \qedhere
 \end{enumerate}
\end{proof}
 
\subsection{Uniformly Hyperbolic Set}

\begin{dfn}[Uniformly Hyperbolic Set]
  \label{dfn:hyperbolicset}
  We say that a set $\Sigma \subset \R \times \R^n$ is a {\em
    uniformly hyperbolic set} when there exist constants $C _\Pi$, $C
  _U$, $\lambda_s$, and $ \lambda_u$ so that all the orbits in the set
  $\Sigma$ are uniformly hyperbolic with the above constants as
  quality measures.
\end{dfn}

If the hyperbolic sets considered lie in a subset of $\R^n$, we just
need to assume that the derivatives of the vector field $\vf$ is
uniformly bounded in a suitable open set containing the hyperbolic set
(it needs to contain all balls of a certain radius centered within the
hyperbolic set).

One interesting example is the Lorenz attractor.  The Lorenz equations
are not bounded in the whole space, but they are bounded in a
neighborhood of the Lorenz attractor.  The Lorenz attractor is not
uniformly hyperbolic but it contains many uniformly hyperbolic sets to
which our theory applies.

We do not assume that the set $\Sigma$ has any particular
structure. In particular, we do not need that the set is closed nor
that it is \emph{locally maximal}, assumptions that are very common in
the theory of hyperbolic systems.

\begin{rmk}[Translation Invariance]
  \label{translation}
  By the uniqueness of solutions of ODE, we can identify an orbit
  with its initial condition.
If an orbit $\{\xx(t)\}$ is hyperbolic, by definition, so are all the
translates $\{\xx(t + \tau)\}$ and they have the same quality
measures.
Hence, when considering a hyperbolic set for an ODE, we can
identify the set of hyperbolic trajectories
with a set in $\R^n$ invariant under the flow. 
\end{rmk}

\begin{rmk}[Loss of invariance]
  \label{noinitial}  
  When the perturbation $\pp$ is time-dependent, then the translation
  invariance of hyperbolic solutions may be lost.  For the
  applications to state-dependent delay or advance equations, where
  the phase space is not clear, there is no easy way to identify the
  space of solutions with the space of initial conditions. Hence, for
  our goal in this paper, it is better to think of a hyperbolic set as
  a collection of trajectories rather than as a set of initial
  conditions.
\end{rmk}

\begin{rmk} 
  In the standard theory of uniformly hyperbolic sets, it is natural
  to consider the splittings along a trajectory $\{\xx(t)\}$ not as
  functions of the time, but as functions of the base point.  It is a
  standard result in hyperbolic systems \cite{Anosov69,
    KatokH95,FisherH19} that the stable and unstable bundles depend on
  the base point in a H\"older way.

 For reasons indicated in Remark~\ref{noinitial}, we, instead, choose
 to study the splittings as functions of time.  We will be able to
 prove some regularity from the space of trajectories of the
 unperturbed system to the space of trajectories of the perturbed
 system. The regularity is somewhat technical since it involves
 weighted spaces.
\end{rmk} 

\section{Construction of perturbative solutions} \label{sec.constsol}
In this section, we introduce the main idea of our result.  We will
describe the geometric motivations and the manipulations needed to
transform the problem considered into a fixed point problem. We
postpone a precise discussion of the regularity assumptions and other
sophistication.  Indeed, those precise assumptions are motivated to
make the arguments in this section work.

Of course, readers interested only in precise formulations can move
directly to Section~\ref{sec.formulation} and use the present section
as a reference for the notations we introduce.

\medskip

The formalism we present resembles the proof of structural stability
for Anosov Flows, which involves a reparametrization of time and a
geometric change of the trajectories. These are the two main ideas we
apply. Nevertheless, in contrast with many proofs of the structural
stability, the reparameterization and the corrections are done
differently for each trajectory and we formulate a different
functional equation for each trajectory. As mentioned before,
given the fact that the functional equations we need involve
the composition operator, several formally equivalent equations 
may have different analytic properties.  We have
carefully chosen a formulation that leads to smooth dependence on
parameters.

\subsection{Form of the correction}
Let $\{\xo(t)\}$ be a uniformly hyperbolic orbit of the unperturbed
system in \eqref{model}, then for $\ep \ne 0$ we consider a solution
of \eqref{model} close to $\{\xo(t)\}$ of the form
\begin{equation} \label{newform}
  \xx(t) =  (\xo + \hx) \circ \hatt(t), 
\end{equation}
where $\hx = O(\ep)$ and $\dop\hatt(t) = 1 + O(\ep)$ are the
correcting unknowns. The term $\hatt$ encodes the internal dynamics of
the new solution while $\hx$ is the displacement from $\xo$. Note that
the form \eqref{newform} is reminiscent of the Anosov Shadowing
Theorem using the functional analysis approach, see
Remark~\ref{rmk:anosov}.

We formulate functional equations (invariance equations) for the
unknown pairs $\hx \colon \R \to \R ^n$ and $\hatt\colon \R \to
\R$. These equations require that $\xx(t)$ in \eqref{newform} is a
solution of the equation \eqref{model}.  We then solve the invariance
equations by fixed point methods using geometric assumptions on
$\{\xo(t)\}$.

We find it more convenient to use the
unknown vector field $\vft$ 
\begin{equation}\label{newunknown}
  \dot \hatt(t) = \vft \circ \hatt(t),
\end{equation}
associated to $\hatt$ instead of $\hatt$ itself, since the invariance
equations are simpler in terms of $\vft$.  It is clear that given
$\hatt$ we can obtain $\vft$ by taking derivatives.  Conversely, given
$\vft$, we can recover $\hatt$ using the differential equation
\eqref{newunknown}. After a normalization condition that sets
$\hatt(0)=0$, see Section~\ref{nonunique}, we see that if $\vft$ is
$C^1$, which implies that we can solve the ODE uniquely for all time, a flow $\hatt \colon \R \to
\R$ is determined. Therefore, we can consider $\vft$ and $\hatt$ as equivalent
unknowns. Although, of course, the natural function spaces for them
are different. Going from $\hatt$ to $\vft$ involves a loss of
derivative, but going from $\vft$ to $\hatt$ gains a derivative.  We
have collected some of these subtleties in Lemma~\ref{lem.vfregular}.

\smallskip

From here, we adopt the convention that $\vft$ and $\hatt$ are related
as indicated above in \eqref{newunknown}. When there is a need to
discuss the dependence of $\hatt$ on $\vft$, we will write
\begin{equation}\label{solution}
  \hatt = {\cal S}[\vft],
 \end{equation}
to indicate that $\hatt$ is the solution of \eqref{newunknown} with a
fixed initial condition.  We refer to ${\cal S}$ as the solution
operator. The operator is defined for $C^1$ vector fields.  If $\vft$
is bounded away from zero, ${\cal S}[\vft]$ will be a diffeomorphism
on $\R$.

\smallskip

To find a locally unique pair $(\vft,\hx)$ from the invariance
equations, we require appropriate normalizations, see
Section~\ref{nonunique} later on. As a consequence of that uniqueness,
we will be able to discuss smooth dependence on parameters around the
initial orbit $\{\xo(t)\}$. However, since we are dealing with
functions defined on the whole line, this will involve some
subtleties.

The strategy of invariance equations treated by functional analysis is
very different from the strategy based on defining an evolution in a
space of functions associated to \eqref{model} and finding hyperbolic
solutions.  Notably, we start by fixing the form, \eqref{newform}, and
finding functions of this form that satisfy \eqref{model}.  There are
cases where invariant objects of systems without globally defined
solutions have been studied \cite{Llave09, ChengL20b} and we will use
some of the techniques developed there.

We also note that the invariance equations in this strategy can be
studied numerically or using formal expansions.  Numerical treatments
of the equations for periodic orbits and their stable manifolds in
simple models were done in \cite{GimenoYL21}. An interesting problem
is to extend the above numerical methods for periodic solutions to the
solutions with arbitrary time dependence considered here.

\subsection{Non-Uniqueness of the parametrization and normalization conditions} \label{nonunique} 
The expressions in \eqref{newform} is underdetermined.  There are many
representations of the same function $\xx(t)$ using different unknown
pairs. Indeed, given a solution $(\hatt, \hx)$ for the invariance
equations and any diffeomorphism $w$ of $\R$,
\begin{equation} \label{underdeterminacy}
\begin{split}
  \xx(t) &= ((\xo + \hx)\circ w) \circ (w^{-1} \circ \hatt) (t)  \\
         &= \bigl(\xo + (\xo \circ w - \xo + \hx \circ w)\bigr) \circ (w^{-1} \circ \hatt) (t),
\end{split}
\end{equation}
provides another choice $(\hatty,\hy)\bydef (w ^{-1} \circ \hatt,~ \xo
\circ w - \xo + \hx \circ w)$ for the solution.

\smallskip

The underdeterminacy \eqref{underdeterminacy} can be avoided by
imposing normalizations that simplify the treatment and lead to 
local uniqueness of the solution. We choose two normalizing conditions:
\begin{align}
 \hx \circ\hatt(t) &\in \Et^s_{\hatt(t)} \oplus
 \Et^u_{\hatt(t)}, \label{normalization1} \\ \hatt(0) &=
 0. \label{normalization2}
\end{align}
Of course, the fact that
\eqref{normalization1} and \eqref{normalization2} are good
normalizations will become apparent when we show that we can find
locally unique solution of the invariance equations satisfying them.

The heuristic reason for \eqref{normalization1} is that adding a
component of $\hx$ in the direction of the flow is roughly equivalent
to adjusting $\hatt$, which can be seen from \eqref{underdeterminacy}
and $\xo \circ w - \xo \approx \xo' (w-\Id)$. Meanwhile,
\eqref{normalization2} can be justified by choosing the origin of $t$
in the reference line.

\subsection{Formulation of the functional equations characterizing a solution of \eqref{model} } \label{sec.inveqs}

To derive functional equations for the unknowns $(\vft, \hx)$, we
substitute \eqref{newform} into \eqref{model}, yielding
\begin{equation}\label{goal0} 
   \vft(\hatt(t))  (\xo + \hx)'  (\hatt(t)) =
  \vf \circ(\xo  + \hx)\circ \hatt(t) +
  \ep \Pop[ (\xo + \hx)\circ \hatt, \ep,\prm](t),
\end{equation} 
where $\hatt={\cal S}[\vft]$ as in \eqref{solution}, and ${}'$ denotes the derivative, and
$\Pop$ is a functional operator. Explicitly,
\begin{equation} \label{Pop}
  \Pop[u,\ep, \prm](t) \bydef \pp (t, u _t, \ep, \prm),
\end{equation}
where $\pp$ is the perturbative map in \eqref{model}.

Note that we are using that $\vft(\hatt(t))$ is a number, so
that we can put the product by it either as a prefactor
or as a postfactor as would come from the chain rule. 
\smallskip

Now, we start to rewrite the equation \eqref{goal0} separating the
small terms.  We first consider the linear approximation of the vector
field $\vf$ along the hyperbolic orbit $\xo$ of the ODE
\begin{align}
     \vf \circ(\xo + \hx)\circ \hatt(t) &= \vf \circ \xo \circ
     \hatt(t) + \dop \vf \circ \xo\circ \hatt(t) \hx \circ \hatt(t) +
     T[\xo, \hx](\hatt(t)), \nonumber \intertext{where} T[\xo,
       \hx](\hatt(t))&\bydef \vf \circ(\xo + \hx)\circ \hatt(t) - \vf
     \circ \xo \circ \hatt(t) - \dop \vf \circ \xo\circ \hatt(t) \hx
     \circ \hatt(t) \label{eq.T}
\end{align}
is the remainder of the first order Taylor expansion. 

\medskip

Using \eqref{eq.T}, 
equation \eqref{goal0} is rewritten as
\begin{equation} \label{goal1}
\begin{split}
    \vft(\hatt(t))  \hx'(\hatt(t)) 
     = 
       (1 - \vft \circ \hatt(t)) & \vf \circ \xo\circ \hatt(t) + \dop\vf\circ \xo\circ \hatt(t)  \hx\circ \hatt(t) \\
       &+T[\xo, \hx](\hatt(t))+ \ep \Pop[(\xo + \hx) \circ \hatt,\ep,\prm](t) .
\end{split}
\end{equation}
We apply the time change $\rho = \hatt(t)$, add and subtract
$\vft(\rho) \dop\vf\circ \xo(\rho) \hx(\rho)$ in \eqref{goal1} to
obtain
\begin{equation} \label{goal2}
 \vft (\rho) \hx'(\rho) = \vft(\rho)  \dop\vf\circ \xo(\rho) \hx(\rho) + (1 - \vft(\rho))  \vf\circ \xo(\rho) +
   \B[\vft, \hx](\rho) + \ep \varphi[\vft,\hx](\rho)
\end{equation}
where for typographical reasons, we introduce $\B$ to capture the 
``quadratically'' small terms, i.e.
\begin{equation}\label{Bdefined}
  \B[\vft, \hx](\rho)\bydef (1 - \vft(\rho))  \dop\vf\circ \xo(\rho) \hx (\rho) + T[\xo, \hx](\rho),
\end{equation}
and $\varphi$ to represent the term from $\Pop$,
\begin{equation} \label{eq.Px}
 \varphi[\vft, \hx](\rho) \bydef \Pop[(\xo + \hx)\circ \hatt, \ep , \prm](\hatt^{-1}(\rho)) = \pp \bigl(\hatt^{-1}(\rho), ((\xo + \hx)\circ \hatt)_{\hatt^{-1}(\rho)}, \ep , \prm\bigr).
\end{equation}

\begin{rmk}
Note that $\varphi$ depends on $\hx$, $\hatt$, $\xo$, the perturbation
$\pp$, the perturbative parameter $\ep$, and the parameter $\prm$ in
equation \eqref{model}. By using $\hatt={\cal S}[\vft]$ in
\eqref{solution}, we consider $\varphi$ as a functional which produces
a function from $\R$ to $\R^n$ given the vector field $\vft$ and the
deformation $\hx$. To simplify the notation, we denote
$\varphi[\vft,\hx]$ without writing explicitly other dependencies.
\end{rmk}

\smallskip

We first consider the center direction of equation \eqref{goal2} to
obtain equation \eqref{c-proj}. Then we use the uniform hyperbolicity
of $\xo$ and the normalization $\hx = \hx ^s + \hx ^u$ in
\eqref{normalization1} to derive \eqref{s-proj} and \eqref{u-proj}.
Thus, by Lemma~\ref{lem:caracinv}, \eqref{goal2} is equivalent to the
following three equations and that the initial conditions of
\eqref{s-proj} and \eqref{u-proj} are in the corresponding bundles:
\begin{align}
    0 &=  \Pi^c_{\rho}  (1-\vft(\rho))\cdot \vf\circ\xo(\rho)  + \Pi^c_{\rho} \bigl( \B[\vft, \hx](\rho) + \ep \varphi[\vft,\hx](\rho)\bigr) , \label{c-proj}\\
    (\hx ^s)'(\rho)&=
    \dop\vf\circ\xo(\rho) \hx^s(\rho) +  \Pi^s_{\rho} \frac{1}{\vft(\rho)}\bigl( \B[\vft, \hx](\rho) + \ep \varphi[\vft,\hx](\rho)\bigr) , \label{s-proj}\\
    (\hx^u)'(\rho)&= 
    \dop\vf\circ\xo(\rho) \hx^u(\rho) +  \Pi^u_{\rho}  \frac{1}{\vft(\rho)}\bigl( \B[\vft, \hx](\rho) + \ep \varphi[\vft,\hx](\rho)\bigr). \label{u-proj}
\end{align}

\smallskip

Our goal is to transform \eqref{c-proj}--\eqref{u-proj} into a fixed
point equation for the unknowns $(\vft,\hx)$.  We will define
operators $\op_c$, $\op_s$, and $\op_u$ of $\vft$ and $\hx$ whose
fixed point solves \eqref{c-proj}--\eqref{u-proj}.

\smallskip

The operator $\op_c$ comes from isolating $\vft(\rho)$ in \eqref{c-proj}. More explicitly,
\begin{equation}\label{Gammac}
 X(\rho)= \op_c[\vft,\hx](\rho)  \bydef  1 +\frac{\bigl\langle \Pi^c_{\rho} \bigl(\B[\vft, \hx](\rho) + \ep \varphi[\vft,\hx](\rho)\bigr), \vf\circ \xo(\rho)\bigr\rangle}{\langle \vf\circ \xo(\rho),\vf\circ \xo(\rho)\rangle},
\end{equation}
where $\langle \cdot, \cdot \rangle$ denotes the inner product in $\R
^n$.

\smallskip

To solve equations \eqref{s-proj} and \eqref{u-proj}, we apply the
variation of parameters formula on the bundles respectively and take
appropriate limits.  The procedure is very similar to the method of
\cite{Perron29,Cotton11} in the study of invariant manifolds. Note
that, although we are not considering invariant manifolds here, our
ideas can be compared to those used for studying normally hyperbolic
invariant manifolds.

We first obtain that
for $\rho _0 \geq \rho \geq -\rho _0$, 
\begin{equation}\label{Gamma_s_almost}
\begin{split}
 \hx^s(\rho)  &= \int_{-\rho _0}^{\rho} U^s(\rho; v)  \Pi^s_{v} \frac{1}{\vft(v)}\bigl(\B[\vft, \hx](v) + \ep \varphi[\vft,\hx](v)\bigr)   \, dv +   
 U^s(\rho;-\rho_0) \hx^s(-\rho _0), \\
 \hx^u(\rho) &= -\int_{\rho}^{\rho _0} U^u(\rho; v)  \Pi^u_{v} \frac{1}{\vft(v)}\bigl( \B[\vft, \hx](v) + \ep \varphi[\vft,\hx](v)\bigr)   \, dv + U^u(\rho;\rho_0) \hx^u(\rho _0).
 \end{split}
\end{equation}
Using the bounds \eqref{eq.expoU} on the evolution operators $U^s/U^u$
and assuming that $\hx^s $ and $ \hx^u$ are bounded (or that if they
grow, the growth rate is less than $\lambda_s$, $\lambda _u$ respectively, we let
$\rho _0 \to +\infty$ and have
\begin{equation}\label{Gammasu} 
\begin{split}
  \hx^s(\rho) = \op_s[\vft,\hx](\rho) & \bydef \int_{-\infty}^{\rho} U^s(\rho; v)  \Pi^s_{v} \frac{1}{\vft(v)}\bigl(\B[\vft, \hx](v)+ \ep \varphi[\vft,\hx](v)\bigr)   \, dv, \\
  \hx^u(\rho) = \op_u[\vft,\hx](\rho) & \bydef -\int_{\rho}^{+\infty} U^u(\rho; v)  \Pi^u_{v} \frac{1}{\vft(v)}\bigl(\B[\vft, \hx](v)+ \ep \varphi[\vft,\hx](v)\bigr)   \, dv. \\
\end{split}
\end{equation}

Alternatively, one could check that $\op_s$ and $\op_u$ defined in
\eqref{Gammasu} indeed satisfy the equations
\eqref{s-proj}--\eqref{u-proj}. Let us provide the details for the
stable case: Taking derivatives w.r.t. $\rho$ in $\op
_s[\vft,\hx](\rho)$ and using the fundamental theorem of calculus
\begin{equation*}
\begin{split}
    \frac{d}{d\rho}\op _s[\vft,\hx](\rho) &= 
    U ^s(\rho; \rho) 
    \Pi^s_{\rho} \frac{1}{\vft(\rho)}\bigl(\B[\vft, \hx](\rho)+ \ep \varphi[\vft,\hx](\rho)\bigr) \\ &\quad + 
    \int_{-\infty}^{\rho} \dop \vf \circ \xo(\rho) U^s(\rho; v)  \Pi^s_{v} \frac{1}{\vft(v)}\bigl(\B[\vft, \hx](v)+ \ep \varphi[\vft,\hx](v)\bigr) \, dv \\
    &= 
    \Pi^s_{\rho} \frac{1}{\vft(\rho)}\bigl(\B[\vft, \hx](\rho)+ \ep \varphi[\vft,\hx](\rho)\bigr) + \dop \vf \circ \xo(\rho) \op _s[\vft,\hx](\rho)  .  
\end{split}
\end{equation*}
To justify the derivative under the integral sign, we observe that the
integrand decays exponentially. The previous derivation is a standard
argument going back to \cite{Cotton11,Perron29} and it has the
advantage showing that \eqref{Gammasu} is the only solution of the
differential equations \eqref{s-proj}--\eqref{u-proj} with growth rate
smaller than $\lambda _{s,u}$ and, in particular, bounded. In
addition, we observe that given the exponential bounds for $U^s$ and
$U^u$ in \eqref{eq.expoU}, if $\B[\vft, \hx]$ and $\ep
\varphi[\vft,\hx]$ are bounded, the $\hx^s$ and $\hx^u$ produced in
\eqref{Gammasu} are bounded and in the corresponding bundles.

Note that to solve \eqref{s-proj} and \eqref{u-proj}, we are not
specifying any initial condition for $\hx ^s$ or $\hx ^u$ explicitly,
only that the solutions are uniformly bounded by an exponential of
time. Indeed, these boundedness requirement fixes the initial
condition for \eqref{s-proj} and \eqref{u-proj}: if we specified an
initial condition not in the trajectories \eqref{Gammasu}, we would
obtain exponential growth solutions with a rate $\lambda _u$ or
$\lambda _s$ and in particular unbounded.

\section{Precise formulation} \label{sec.formulation}

In this section, we introduce the function spaces and revisit the
construction in Section~\ref{sec.inveqs} with precise formulation. We
also provide our main result, see Theorem \ref{main}.

\subsection{The operator} \label{sec.operator}
We consider an operator $\op$ depending on the perturbative parameter
$\ep$ and whose inputs are:
\begin{enumerate}
\renewcommand{\theenumi}{\arabic{enumi}}
\renewcommand{\labelenumi}{{\tt I}\theenumi)}
 \item A non-zero vector field $\vft$ in $\R$ (whose flow is $\hatt =
   \mathcal{S}[\vft]$);
 \item A stable correction $\hx^s$ to the uniformly hyperbolic orbit
   $\xo$; and
 \item An unstable correction $\hx ^u$ to the uniformly hyperbolic
   orbit $\xo$.
\end{enumerate}
The operator $\op$ has outputs:
\begin{enumerate}
\renewcommand{\theenumi}{\arabic{enumi}}
\renewcommand{\labelenumi}{{\tt O}\theenumi)}
 \item A new vector field $Y=\op_c[\vft, \hx ^s, \hx ^u]$ in $\R$, and
   hence its flow $\psi=\mathcal{S}[Y]$ given by the initial value
   problem
 \[ \frac{d}{dt}  \psi(t) = Y \circ \psi(t), \qquad \psi(0)= 0; \]
 \item A new stable correction $\hy^s=\op_s[\vft, \hx ^s, \hx ^u]$ to
   the orbit $\xo$; and
 \item A new unstable correction $\hy^u=\op_u[\vft, \hx ^s, \hx ^u]$
   to the orbit $\xo$.
\end{enumerate}
Note that in both input and output cases, one could write the unknown
orbit correction as a sum of stable and unstable corrections due to
\eqref{eq.splitting} and \eqref{normalization1}. Indeed, $\hx = \hx ^s
+ \hx ^u$ and since the operators $\op_s$ and $\op_u$ make the
corrections on the stable and unstable bundles respectively, the
output $\hy$ satisfies $\hy = \hy ^s + \hy ^u$ as well.

\smallskip

We define the operator $\op$ with components, i.e.
\[ \op \equiv
\begin{pmatrix}
 \op_c \\ \op_s \\ \op_u
\end{pmatrix}.
\]
The operator $\op$ acts on a function space $\X$, which will be
exhaustively detailed in Section~\ref{sec.space}.

\begin{rmk}[Alternatives to the fixed point operator]
 The operator $\op$ admits alternative versions, for instance, using
 the corrected vector field $Y$ instead of $X$ in $\op _s$ and $\op
 _u$. There are similar variations with the other updated expressions
 as well.
 
 All these alternative operators may have a convergence impact in a
 numerical implementation. Nevertheless, to prove the existence and
 uniqueness of the fixed point, the $\op$ defined here will give
 easier inequalities that otherwise can be bounded by simple triangle
 inequalities.
\end{rmk}

\subsection{The spaces considered} \label{sec.space}

As usual for fixed point problems, one tries to get both existence and
uniqueness of fixed point. We look for fixed points of the operator
$\op$ in a product space $\X$ of finitely differentiable maps. The
existence results become better when considering a space of more
differentiable functions, while the uniqueness results will be better
for a bigger space with lower regularity.

The following definitions, although standard, set the notation for the
statement of our main results, see Section~\ref{sec.main}.

\subsubsection{Spaces of Lipschitz differentiable functions}
Let $\ell \geq 0$ be a fixed integer, let $I\subset\R$ be an open
interval, and let $\dop$ be the differential operator.  The space $C
^\ell(I)$ denotes the space of functions defined on $I$, that are
$\ell$ times differentiable, extend continuously to the closure of
$I$, denoted by $\overline{I}$, and whose derivatives are
bounded. More precisely,
\[
C^\ell(I) = C^\ell( I, \R^n ) \bydef \left\{g\colon  I \rightarrow \R^n \, \left| \,
 \begin{minipage}[c]{.42\textwidth}
  $g$ is $\ell$ times continuously differentiable in $I$, the
   derivatives extend continuously to $\overline{I}$, and $\| g \| _{C^\ell}\bydef \max \limits_{0\le j\le \ell}\bigl\{\sup \limits_{x\in I} |\dop ^j g(x)|\bigr\} < +\infty$
 \end{minipage}\right.\right\},
\]
where $|\cdot|$ denotes a norm in $\R ^n$. By our definition, the
space $C^\ell(I)$ is a Banach space for each $\ell$.  If the domain or
range of the functions are understood, we will suppress it from the
notation. We use the identity $\dop ^0 = \Id$ and define $C ^0$ as the
set of continuous functions with bounded $C ^0$-norm.

\smallskip

We denote the space $C ^{\ell + \Lip}$ as a subspace of $C^\ell$
containing functions whose $\ell$-th derivative is Lipschitz, and we endow
the space with the $\|\cdot\| _{C ^{\ell +\Lip}}$ norm. Explicitly,
\begin{equation*}
 C ^{\ell+\Lip}(I) = C ^{\ell + \Lip}(I, \R ^n) \bydef \bigl\{g \in C ^\ell (I, \R ^n)\colon \| g \| _{C ^{\ell +\Lip}} \bydef \max \bigl\{ \|g \| _{C ^\ell} , \Lip(\dop ^\ell g)\bigr\} < +\infty \bigr\}.
\end{equation*}
The $C ^{\ell + \Lip}$ space is useful for us as it is the $C^0$
closure of the $C^{\ell+1}$ space.

The Lipschitz constant has some properties that we summarize in the
following (known) lemma:
\begin{lem}\label{lem.Lipschitz}
 Let $f, g$ be continuous maps with finite Lipschitz constant and let
 $\lambda$ be a scalar. Then
 \begin{enumerate}
  \item $\Lip(f + g) \leq \Lip (f) + \Lip (g)$;
  \item $\Lip(\lambda f) \leq |\lambda| \Lip (f)$;
  \item $\Lip(f  g) \leq \Lip (f) \| g \| _{C ^0} + \|f \| _{C ^0} \Lip(g) \leq 2 \|f \| _{C ^1} \| g \| _{C ^1}$; and
  \item $\Lip(f\circ g) \leq \Lip( f) \Lip (g)$. 
 \end{enumerate}
\end{lem}

\smallskip

Given $c = (c _0, \dotsc, c _\ell, c _{\ell}^{\Lip}) \in \R _+^{\ell
  +2}$, we denote the ball centered at $l\in C ^{\ell + \Lip}$ with
radius $c $ as
\begin{equation} \label{ballspace}
  \ball ^{\ell+\Lip}_{c}(l) \bydef \bigl\{ g \in C ^{\ell+\Lip} \colon
  | \dop ^j (g-l) |\leq c _j \text{ for } j = 0, \dotsc, \ell\text{
    and } \Lip (\dop ^\ell (g-l)) \leq c _\ell^{\Lip} \bigr\}.
\end{equation}
We will use a similar notation for a ball in $C ^\ell$ space by $
\ball ^{\ell}_{c}(l)$. If $l$ is zero, we will sometimes omit the
center and write $\ball ^\ell_{c}$.

The following straightforward lemma states that the $C ^\ell$ space is
a Banach algebra by our definition. With 
Lemma~\ref{lem.Lipschitz}, we get that $C ^{\ell +\Lip}$ is also a
Banach algebra.

\begin{lem} \label{lem.ClBanach}
The space $C^{\ell}$ is a Banach algebra:
\begin{enumerate}
 \item $f \in \ball^\ell_{c_f}$ and $g \in \ball^\ell_{c_g}$ implies
   $f +g \in \ball^\ell_{c _f + c _g}$;
 \item $\lambda \in \R$ and $g \in \ball^\ell_{c}$ implies $\lambda g
   \in \ball^\ell_{|\lambda| c}$; and
 \item $f \in \ball^\ell_{c_f}$ and $g \in \ball^\ell_{c_g}$ implies
   $fg \in \ball^\ell_{\tilde c}$ with $\tilde c$ only depending on $c
   _f$ and $c _g$.
\end{enumerate}
\end{lem}
\begin{proof}
 \begin{enumerate}
  \item []
    \textit{\stepcounter{enumi}\theenumi}.~and~\textit{\stepcounter{enumi}\theenumi}.~are
    straightforward. We define 
  \[ \tilde c _j \bydef \sum _{k=0}^j \binom j k c _{f,k}c _{g,j-k} \]
  so that \textit{\stepcounter{enumi}\theenumi}. is true by
    Leibnitz rule. \qedhere
 \end{enumerate}
\end{proof}

\paragraph{H\"older space and interpolation inequality}
For the a-posteriori formulation of our main theorem, we briefly
recall the standard notion of H\"older spaces and state the
interpolation inequalities in $C ^r$ spaces \cite{Stein1970}
(A short proof of the interpolation inequalities valid in domains
in Banach spaces can be found in \cite{LO99}). 

\smallskip

Given $\alpha \in (0,1)$, the H\"older space $C ^{0,\alpha}(I)$ is the
set of functions $g \in C^0$ such that the H\"older semi-norm
\begin{equation*}
 [g] _{C ^{0,\alpha}} = [g] _{C ^{0,\alpha}(I)}\bydef \sup
 _{\substack{x,y \in  I \\ x\ne y}} \frac{|g(x) -
   g(y)|}{|x-y| ^{\alpha}}
\end{equation*}
is finite. Note that for $\alpha=1$, this semi-norm is exactly the Lipschitz
constant. Similarly, given $k\geq 0$ an integer, the H\"older space $C
^{k,\alpha}$ is defined by
\begin{equation*}
 C ^{k,\alpha} = C^{k,\alpha}(I) \bydef \bigr\{ g\in C^{k}(I)\colon \|
 g \| _{C ^{k,\alpha}} \bydef \max\{\| g \| _{C ^k(I)}, [\dop ^k g]
 _{C^{0,\alpha}(I)}\} < +\infty \bigr\}.
\end{equation*}
In particular, $C ^{k,1} \equiv C ^{k+\Lip}$. Notice that this yields
the space inclusions
\begin{equation*}
 C ^{\infty} \subset \dotsb \subset C^2 \subset C ^{1+\Lip} \subset C
 ^{1,\alpha} \subset C^1 \subset C^{\Lip} \subset C ^{0,\alpha}
 \subset C ^0.
\end{equation*}
It is then common to define the $C^r$ space when $r >0$ is not an
integer. That is, $C^r \bydef C^{k,\alpha}$ where $r = k + \alpha$,
$k$ being an integer and $\alpha \in (0,1)$.

\smallskip

The space $C ^r$ admits interpolation inequalities. We now quote the
result (omitting domain assumptions) from \cite{Kol49,Had98,LO99}
which states that if $0 \leq r < t$ and $g \in C ^t$, then there is a
constant $M _{r,t}>0$ such that
\begin{equation*} 
    \| g \| _{C ^{ \theta r + (1 - \theta) t}} \leq M _{r,t} \| g \|
    _{C ^r} ^\theta \| g \| _{C ^{t}}^{1-\theta},
\end{equation*}
for any $\theta \in [0,1]$. This is equivalent to consider $s \in
(r,t)$ and $\mu \bydef \frac{t-s}{t-r}$ and rewrite the inequality as
\begin{equation} \label{eq.cellineq}
    \| g \| _{C ^{s}} \leq M _{r,t} \| g \| _{C ^r} ^\mu \| g \| _{C
      ^{t}}^{1-\mu}.
\end{equation}

An interesting remark for the applications is
that the unit ball in $C^{k,\alpha}[a,b]$ where $k \in \mathbb{N}$,
$0 < \alpha \le \Lip$, $-\infty \le a$, $ b \le \infty$
is compact (and therefore closed) in the $C^0$ topology.

Arzela-Ascoli theorem shows that the ball is precompact, and
we also have that $\dop^k u_n$ converges uniformly and are
uniformly $C^\alpha$, then the limit is also $C^\alpha$ with the same constant.

\subsubsection{Contraction space for time-dependent perturbation} \label{sec.contraspace}
When the perturbation in \eqref{model} depends on time in a bounded
manner, we need to bound the difference of time reparametrizations,
that is, to bound flows for two vector fields. Hence, for some fixed $\eta > 0$, we introduce the
Razumikhin norm for continuous functions on an open interval $I
\subset \R$, which is defined as
\begin{equation}\label{razumikhin-norm}
 C _\eta=C _\eta(I, \R ^n) \bydef \left \{g: I\to\R ^n \text{ is
   continuous} \colon \| g \| _{C_\eta} = \| g \| _{\eta} \bydef \sup
 _{\rho \in I} |g(\rho) | e ^{-\eta |\rho|} < +\infty \right\}.
\end{equation}
The parameter $\eta$ will eventually be chosen to ensure that the operator $\op$ is a
contraction. Note that $\| g \| _{C_\eta} \leq \| g \| _{C ^0}$ for
all $g \in C _\eta$. Moreover, under some assumption on $\eta$, the operator ${\cal S}$ in
\eqref{solution} can be bounded as:
\begin{equation*}
 \|{\cal S}[\vft] - {\cal S}[\vfty] \| _{C_\eta} \leq \ct c \| \vft - \vfty \| _{C ^0},
\end{equation*}
for a constant $\ct c$ depending on $\eta$ and the bound for the
Lipschitz constants of the vector fields $\vft$, $\vfty$. This is
formally proved in the following Lemma~\ref{lem.phipsi-fwd-contract}.

\begin{lem} \label{lem.phipsi-fwd-contract}
 Let $\hatt$ and $\hatty$ be flows of the vector fields $\vft$ and $
 \vfty$ in $\R$ in a ball $\ball ^{\Lip}_{(\bt _0, \bt _1)}(1)$
 respectively with zero initial condition at zero. If $\bt _0 \in
 (0,1)$ and $\eta > \bt _1 > 0$, then
 \begin{equation*}
  | \hatt(\rho) - \hatty(\rho) | e ^{-\eta |\rho|} \leq \frac{\|\vft - \vfty \| _{C ^0}}{e(\eta -\bt _1)}.
 \end{equation*}
\end{lem}
\begin{proof}
 Let us define the solution operator of the ODE $\dot \hatt = \vft
 \circ \hatt$ generated by vector field $\vft$ as
 \[
   \Upsilon _{\vft}(\hatt)(\rho) = 
   \begin{cases}
    \displaystyle \int _0 ^\rho \vft \circ \hatt(s)\, ds & \rho \geq 0 \\
    \\
    \displaystyle -\int _0 ^{-\rho} \vft \circ \hatt(-s)\, ds & \rho < 0.
   \end{cases}
 \]
We define operator $\Upsilon _{\vfty}$ similarly. Let $\widetilde
\hatt \colon \R \to \R$ be another function. Then, for $\rho < 0$ (and
similarly for $\rho \geq 0$),
 \begin{align*}
     | \Upsilon _{\vft}(\hatt) (\rho) - \Upsilon _{\vft}(\widetilde \hatt)(\rho)| e ^{-\eta |\rho|} 
     &\leq e ^{\eta \rho}\int _0 ^{-\rho} | \vft \circ \hatt (-s) - \vft \circ \widetilde\hatt(-s) | \, ds  
     \leq  e ^{\eta \rho}\Lip(\vft) \int _0 ^{-\rho} | \hatt (-s) - \widetilde\hatt(-s) | \, ds \\
     & \leq \| \hatt - \widetilde \hatt \| _\eta  \Lip(\vft) e ^{\eta \rho}(e ^{-\eta \rho}-1)  
     \leq \frac{\Lip(\vft)}{\eta} \| \hatt - \widetilde \hatt \| _\eta. 
 \end{align*}
Thus, if $\eta > \bt _1 \geq \max\{\Lip (\vft), \Lip (\vfty)\}$, then
$\Upsilon _{\vft}$ and $\Upsilon _{\vfty}$ are contractions.

\smallskip

 Given $\Upsilon _{\vfty}(\hatty) = \hatty$ a fixed point, the
 a-posteriori estimates for $\Upsilon _{\vfty}$ says:
 \[
   \| \hatt - \hatty \| _\eta \leq \biggl(1 - \frac{\bt _1}{\eta}\biggr)  ^{-1}\| \hatt- \Upsilon _{\vfty} (\hatt) \| _\eta.
 \]
 Now for $\rho < 0$ (and again similarly for $\rho \geq 0$)
 \begin{align*}
     |\hatt(\rho) - \Upsilon _{\vfty} (\hatt)(\rho)| e ^{-\eta |\rho|}
     &\leq e ^{\eta \rho} \int _{0}^{-\rho} |\vft\circ \hatt(-s)-\vfty\circ \hatt(-s) |\, ds \leq e ^{\eta \rho} \frac{-\eta \rho}{\eta} \| \vft - \vfty \| _{C ^0}.
 \end{align*} 
 By taking supremum and using the fact that $\sup \limits_{u < 0} -e
 ^{u} u = \sup \limits_{u \geq 0} e ^{-u} u = e ^{-1} $, we conclude
 that
 \[
   \| \hatt - \hatty \| _\eta \leq \biggl(1 - \frac{\bt _1}{\eta}\biggr)  ^{-1}\frac{1}{\eta e} \| \vft - \vfty \| _{C ^0}. \qedhere
 \]
\end{proof}

\paragraph{Interpolation inequality}

The $C _\eta$ space also admits interpolation inequalities based on
the $C ^r$ interpolation inequalities for functions defined on finite
intervals. Indeed, given $g \colon [a,b] \to \R ^n$, we have
\begin{align*}
    \| g |_{[a,b]} \| _{C_\eta} \leq \| g |_{[a,b]} \| _{C ^0} \leq e ^{\eta \max\{|a|,|b|\}} \| g |_{[a,b]} \| _{C_\eta}.
\end{align*}
Therefore, when $r=0$, we can rewrite the interpolation inequality
\eqref{eq.cellineq} with $\|\cdot\|_\eta$, yielding
\begin{equation*}
 \| g |_{[a,b]} \| _{C^{(1-\theta)t}} \leq M _{0,t} \, \| g |_{[a,b]} \| _{C ^0} ^\theta \| g |_{[a,b]} \| _{C ^{t}}^{1-\theta} \leq M _{0,t} \, e ^{\eta  \theta \max\{|a|,|b|\}} \, \| g |_{[a,b]} \| _{C _\eta} ^\theta \| g |_{[a,b]} \| _{C ^{t}}^{1-\theta}.
\end{equation*}
This interpolation property will be used in the a-posteriori
formulation of the main result, see Theorem~\ref{main} for a formal
formulation and Section~\ref{sec.aposteriori} for a detailed
discussion.

\subsubsection{The operator space for \texorpdfstring{$\op$}{the operator}}
The operator $\op$ takes values from a space $\X$, which is the product space of three spaces, one for each input of the
operator, endowed with the product norm.

\smallskip

The first input of the operator is the vector field $\hvft$ which corrects along the tangent direction
of the orbit. We consider a ball
centered at $1$ of functions $\R \to \R$ that are $\ell$ times
differentiable, with $\ell$-th derivative Lipschitz, that is, $\hvft
\in \ball ^{\ell+\Lip}_{\bt}(1)$ for $\bt \in \R _+^{\ell+2}$.  We
will see in Lemma~\ref{lem.vfregular} that the flow $\hatt={\cal
  S}[\vft]$ gains one regularity.

\smallskip

The other two inputs will be taken in a ball centered at the origin of
functions $\R \to \R ^n$ that are $\ell+1$ times differentiable, with
$(\ell+1)$-th derivative Lipschitz, i.e., let $\hx ^s \in \ball
^{\ell+1+\Lip}_\bs$ and $\hx ^s \in \ball ^{\ell+1+\Lip} _\bu$ for
$\bs, \bu \in \R _+ ^{\ell+3}$. Because of the normalization 
\eqref{normalization1}, $\hx = \hx ^s + \hx ^u$ unequivocally,
therefore $\hx \in \ball ^{\ell+1+\Lip} _{\bs + \bu}$.

\smallskip

Then, the space for $\op$-inputs consists in product of balls:
\begin{equation} \label{eq.ballspaxie4op}
(\hvft, \hx ^s, \hx ^u)\in \X = \X ^{\ell}_{\bt, \bs, \bu} \bydef \ball ^{\ell+\Lip}_\bt (1) \times \ball ^{\ell+1+\Lip}_\bs(0) \times \ball ^{\ell+1+\Lip}_\bu(0).
\end{equation}
By construction, the center projection's range can always be
identified with $\R$, while the stable and unstable projection ranges
are elements in $\R ^n$ belonging to subspaces of $n _s$ and $n _u$
dimension respectively, see Definition~\ref{dfn.unihypersol} and
Remark~\ref{rmk.projections}. Therefore, the ball $\ball
^{\ell+\Lip}_\bt$ lies in the space of functions $\R \to \R$, while
the other two balls in $\X$ are for functions $\R \to \R ^n$. Note
that the space $\X $ is closed under $C^0$-norm.

\smallskip

To fix the notations, we specify the components of the vectors $\bt,
\bs$, and $\bu$ appearing in $\X$:
\begin{equation} \label{eq.tsu-components} 
\begin{split}
 \bt &\bydef (\bt _0, \dotsc, \bt _{\ell}, \bt _{\ell}^{\Lip}) \in \R _+^{\ell+2}, \\
 \bs &\bydef (\bs _0, \dotsc, \bs _{\ell}, \bs _{\ell+1}, \bs
_{\ell+1}^{\Lip}) \in \R _+^{\ell +3}, \\
 \bu &\bydef (\bu _0, \dotsc, \bu _{\ell}, \bu _{\ell+1}, \bu _{\ell+1}^{\Lip}) \in \R _+^{\ell +3}.
\end{split}
\end{equation}
These constants in \eqref{eq.tsu-components}, jointly with the
perturbative parameter $\ep$, are the ones that we will constrain in a
finite set of inequalities to ensure that the operator $\op$, defined
in Section~\ref{sec.operator}, maps $\X$ into itself and is
contractive for some distance.

\medskip

\begin{rmk}[on the special constant $\bt _0$]
  There is a crucial requirement on the value $\bt _0$ to be in the
  interval $(0,1)$ by Lemma~\ref{lem.vfregular}.
\end{rmk}

\smallskip

\begin{rmk}[on the absence of the delay in the space definition]
 We stress that because we are considering a special type of solutions
 (uniformly hyperbolic solutions) under perturbation, we are able to consider a special
 space of functions that is not affected by any delay or functional
 information of the perturbative map $\pp$. In particular, the
 constant $h > 0$ in the ``history domain'' does not even enter in the
 space where we apply the fixed point approach as long as the
 perturbative hypotheses \eqref{Hep1}--\eqref{Hep2}, discussed in
 Section~\ref{sec.main}, hold.
\end{rmk}

\smallskip

\begin{rmk}[on the regularity]
  The fact that $\hvft$ belongs to a space with one degree of regularity less
  than the other functions 
 grants hypothesis  \eqref{Hep1} in
 Theorem~\ref{main}. This allows us to consider
 perturbations $\pp$ that lose one derivative, a property that will be
 exploited in the applications, see Section~\ref{sec.examples}.
 Notably, that property enables the study of neutral equations, equations with small delays
 and the equations of Wheeler-Feynman electrodynamics. 
\end{rmk}

\subsection{The main results} \label{sec.main}
We establish the main results that, under appropriate hypotheses, the
operator $\op$ in Section~\ref{sec.operator} has a fixed point which
is locally unique. Then by the construction of $\op$ in
Section~\ref{sec.formulation}, this fixed point will be a solution of
\eqref{model} under the functional perturbation.

The result has two sets of hypotheses: a first set,
\eqref{H01}--\eqref{H02}, concerning the unperturbed orbit; and a
second set, \eqref{Hep1}--\eqref{Hep2}, on the perturbation $\pp$. The
existence of a solution is ensured by \eqref{Hep1} and its uniqueness
by \eqref{Hep2}.

As indicated in Section~\ref{overview}, the only things to check
are the fact that the operator $\op$ map a smooth ball
into itself and that it is a contraction in low regularity for all functions in
such a smooth ball.  We will show
that this follows from some simple hypothesis on \eqref{model}
and we will verify the hypotheses of Theorem~\ref{main} in concrete
examples of interest. Of course, for each of the models, one could
formulate the operator $\Gamma$ directly and verify the propagated
bounds and the low regularity contraction.

\begin{thm}\label{main}
 We consider the differential equation \eqref{model}. Let $\ell \geq
 0$ be an integer, and $\mu_0\in\R$ a fixed parameter.  Assume that
 the unperturbed system satisfies:
 \begin{enumerate}
\renewcommand{\theenumi}{H\textsubscript{0}\arabic{enumi}}
\renewcommand{\labelenumi}{\sl \theenumi)}  
  \item \label{H01} There is a uniformly hyperbolic
    solution $\{\xo(t)\}_{t\in \R}$, see
    Definition~\ref{dfn.unihypersol}.
  
  \item \label{H02} The function $\vf$ is $C^{\ell+2+\Lip}$ and
    bounded away from zero in a $\delta$-neighborhood of the orbit
    $\{\xo(t)\} _{t\in \R}$.
 \end{enumerate}
 
 \smallskip
 
 Assume that the perturbative map $\pp \colon \R \times
 C^{\ell+1+\Lip}([- h, h], \R ^n) \times (0,1)^2 \to \R ^n$ in
 \eqref{model} defines the operator
 \[
 \Pop \colon C ^{\ell+1+\Lip}(\R, \R ^n) \times (0,1) \to C
 ^{\ell+\Lip}(\R, \R ^n), \qquad \Pop[u,\ep](t) \bydef \pp(t, u
 _t, \ep, \prm_0)
 \]
 such that: 
 \begin{enumerate}
\renewcommand{\theenumi}{H\textsubscript{$\ep$}\arabic{enumi}}
\renewcommand{\labelenumi}{\sl \theenumi)}  
 \item \label{Hep1} For all $\ep\in (0,1)$, $u \in C
   ^{\ell+1+\Lip}(\R, \R ^n)$, $t \in \R$, 
 \begin{equation*}
  \biggl| \frac{d ^j}{dt ^j} \Pop[u, \ep,\prm](t) \biggr| \leq C _j F_j (\| u \| _{C ^{j+1}}) \quad j = 0, \dotsc,
   \ell, \text{ and} \quad 
  \Lip\biggl( \frac{d ^{\ell}}{dt ^{\ell}} \Pop[u,\ep ,\prm] \biggr) \le C _{\ell}^{\Lip} F_\ell^{\Lip} (\| u \| _{C ^{\ell+1+\Lip}}),
 \end{equation*}
 where the constants $C$'s are positive and the functions $F$'s are
 continuous and increasing on $\R$.
 \end{enumerate} 
 Then there exists $\ep _0 \in (0,1)$ such that for all $\ep \in (0,
 \ep _0)$, there are differentiable maps $\hatt $ and $\hx$ such that
 $\dop \hatt$ is in $C ^{\ell+\Lip}$, $\hx$ is in $C ^{\ell+1 +\Lip}$,
 and
 \begin{equation}\label{sol}
  \xx \bydef (\xo + \hx)\circ \hatt
 \end{equation}
 is a $C ^{\ell+1 +\Lip}$ solution of \eqref{model}.
 
 \smallskip
 
 Moreover, if the history value $h < +\infty$ and $\Pop$ also
 satisfies
 \begin{enumerate}
\renewcommand{\theenumi}{H\textsubscript{$\ep$}\arabic{enumi}}
\renewcommand{\labelenumi}{\sl \theenumi)}  
\setcounter{enumi}{1}
 \item \label{Hep2} For all $\ep \in (0,1) $, $u^{1},u ^{2} \in C
   ^{\ell +1+ \Lip}(\R, \R ^n)$, and $t,s \in \R$, there are constants
   $\ct L _1$ and $\ct L _2$ such that
 \begin{equation*}
  |\Pop[u^2,\ep](s) - \Pop[u^1,\ep](t)| \leq \ct L _1 |s - t| + \ct L
  _2 \| u^2_s - u^1_t \| _{C ^1([- h, h])}.
 \end{equation*}
 \end{enumerate}
 Then there exists $\ep
 _0' \leq \ep _0$ such that the maps $\hatt$ and $\hx$ in \eqref{sol} are locally unique for all $\ep \in (0, \ep _0')$.
 
 \medskip 
 
 Furthermore, if \eqref{Hep1}--\eqref{Hep2} hold, we obtain
 a-posteriori result. Given an initial guess $(\vft _{(0)}, \hx
 _{(0)}^s, \hx _{(0)}^u)$ with errors
 \begin{equation*}
 \ee _c \bydef \op _c [\vft _{(0)}, \hx _{(0)}^s, \hx _{(0)}^u] - \vft _{(0)}, \quad
 \ee _s \bydef \op _s [\vft _{(0)}, \hx _{(0)}^s, \hx _{(0)}^u] - \hx _{(0)}^s,\quad
 \ee _u \bydef \op _u [\vft _{(0)}, \hx _{(0)}^s, \hx _{(0)}^u] - \hx _{(0)}^u,
 \end{equation*}
 on any bounded interval $[a,b]\subset\R$, we have
 \begin{align*}
  \| (\vft - \vft _{(0)})|_{[a,b]} \| _{C ^{j}} &\leq \ct c E_\eta^{\frac{\ell+1 -j}{\ell+1}},  & \text{for } 0 \leq j &\leq \ell\\
  \| (\hx^\sigma - \hx ^\sigma_{(0)})|_{[a,b]} \| _{C ^{j}} &\leq \ct c E_\eta^{\frac{\ell +2-j}{\ell+2}},  & \text{for } \sigma =s,u,~\text{and }0 \leq j &\leq \ell+1,
 \end{align*}
 where $E_\eta\bydef\| \ee _c \| _\eta+ \| \ee _s \| _\eta+ \| \ee _u
 \| _\eta+\| \dop\ee _s \| _\eta+\| \dop\ee _u \| _\eta$, and the
 constant $\ct c$ depends on $j$, $a$, $b$, $\ep$, $\xo$, $\vf$, $h$,
 $\pp$.

 Alternatively, when $E_\eta\le 1$, we have that
 \begin{align*}
  \| \dop ^j(\vft - \vft _{(0)})|_{(0,+\infty)} \| _\eta &\leq \ct c E_\eta^{\frac{1}{j+1}},  & \text{for } 0 \leq j &\leq \ell\\
  \| \dop^ j(\hx^\sigma - \hx ^\sigma_{(0)})|_{(0,+\infty)} \| _\eta &\leq \ct c E_\eta^{\frac{1}{j+1}},  & \text{for } \sigma =s,u,~\text{and }0 \leq j &\leq \ell+1,
 \end{align*}
\end{thm}

\medskip

 \smallskip
\begin{rmk}
    Note that from the assumptions \eqref{H01}, \eqref{H02}, and the
    theory of normal hyperbolicity, we have that for $\sigma \in
    \{c,s,u\}$, the maps $\rho \mapsto \Pi ^\sigma _\rho$ are $C
    ^{\ell+1+\Lip}$. We will use this fact in the proof of our main
    result.
\end{rmk}
\begin{rmk}[on $\hatt$ regularity]
 The map $\hatt$ in Theorem~\ref{main} does not belong to $C ^{\ell
   +1+ \Lip}$ because $\hatt$ is not bounded in $C ^0$, see
 Lemma~\ref{lem.vfregular}. However, the composition $\hx \circ \hatt$
 is $C^{\ell+1+\Lip}$ when $\hx$ is $C ^{\ell+1+\Lip}$ and $\dop
 \hatt$ is $C ^{\ell+\Lip}$.
\end{rmk}
 
 \smallskip

\begin{rmk}[on the perturbative parameter]
 We allow the perturbative map $\pp$ to depend on the perturbative
 parameter $\ep$. In some applications treated in
 Section~\ref{sec.examples}, $\pp$ is obtained by power expansion in
 $\ep$ which implies that $\pp$ may have higher order terms in
 $\ep$. In other applications such as the small delay case, the equation will
 be reformulated such that $\pp$ will explicitly appear.
 
 Moreover, notice that $\ep \in (0,1)$ is not necessarily a
 restriction since one can always scale the map $\pp$ or change its
 sign to admit other ranges of $\ep$.
\end{rmk}
 
 \smallskip

\begin{rmk}[on the choice of regularity space]
The smallness of $\ep$ depends on $\ell$, therefore, the method claims
results on finite regularity only.  The results on analytic regularity
are false without extra hypothesis, see \cite{MPNP, MPN}.

Note that if the equation \eqref{model} is smooth enough, once
obtaining that there are $C^1$ solutions in time, one can use the
equation \eqref{model} to bootstrap the regularity of the solution as
the particular problem allows. In some cases, one may obtain
$C^\infty$ solutions.
\end{rmk}
 
 \smallskip

\begin{rmk}[on the a-posteriori formulation]
The proofs we are going to present are constructive, hence they can be
implemented numerically. The operator concatenates several elementary
operations, some of these operations for a 2D model have been
addressed in a numerical toolkit in \cite{GimenoYL21}.

The formulation we adopted in Theorem~\ref{main} admits an
a-posteriori format which states that if there is an approximate
solution, then close to it there is a true solution.
A-posteriori results can be the basis of computer-assisted proofs
(CAP's) because if one is able to estimate rigorously non-degeneracy
conditions and errors, then one concludes existence of the solution.
The error verification in the approximation is a long finite
calculation taking care of round-off and truncation errors. Some cases of
CAP's have already been used in delay equations,
e.g. \cite{GroothedeMJ17,GimenoLMY23,SzczelinaZ18}.
\end{rmk}

\subsubsection{Parameter dependence result}
\label{sec.smoothprm}

Theorem~\ref{main} is on the case where the parameter $\prm = \prm
_0$ is fixed, while we could modify hypotheses \eqref{Hep1}--\eqref{Hep2}
easily to obtain results on smooth dependence on parameters of the
solution in \eqref{sol}.

Indeed, with $\prm \in (0,1)$, we view $\hatt$ and $\hx$ as maps
$\hatt\colon \R \times (0,1) \to \R$ and $\hx \colon \R \times (0,1)
\to \R ^n$. Therefore, the solution is of the form
 \[
  \xx(t,\prm) = \xo \circ \hatt(t,\prm) + \hx\bigl(\hatt(t,\prm),\prm\bigr).
 \] 

For smooth dependence on parameters, in the first hypothesis on the
perturbation \eqref{Hep1}, we need bounds on the partial derivatives
with respect to $t$ and $\prm$ by functions of $\|u\|_{C^{j+1}}$,
where the norms are understood as the norms for $u\colon \R \times
(0,1) \to \R^n$. Meanwhile, the second hypothesis \eqref{Hep2} should
be changed to include the parameter as follows
\begin{equation*}
  |\Pop[u^2,\ep,\prm](s) - \Pop[u^1,\ep,\prm](t)| \leq \ct L _1 |s
  - t| + \ct L _2 \| u^2_s - u^1_t \| _{C ^1([- h, h]\times(0,1))},
 \end{equation*}
for all $\prm\in (0,1)$. With the changes in the
hypotheses, we obtain that for small $\ep$ the solution $x$ is jointly
$C ^{\ell+1 +\Lip}$ in $t$ and $\mu$.

 Notice that $\prm \in (0,1)$ is not necessarily a restriction since
 we can apply an affine transformation to it.  Also, we can generalize
 our result to consider higher dimensional parameter $\mu$ with
 similar argument.

Proving naturally the smooth dependence on parameters is one of the
advantages of our framework. In general, the smooth parameter
dependence is not trivial for solutions of SDDEs, see
\cite{hart,Walther03}, and it involves extra
assumptions. Nevertheless, we admit that we only search for solutions
of a certain form, as in \cite{YangGL21,YangGL22}.

Theorem~\ref{main} also applies when the parameter appears in the
unperturbed system. Indeed, suppose the unperturbed equation takes the
form
\[\dot x(t) = g(x(t), \prm _0 + h).\] If there is $\prm _0$ such that
\eqref{H01}--\eqref{H02} are satisfied for $\vf(x) \bydef g(x,\prm
_0)$, then we define $Q$ as
\begin{equation*}
  Q(x, h) = \int _0 ^1 \dop _\prm g (x, \prm _0 + \sigma h) d \sigma,
\end{equation*} 
and consider $\dot x(t) = g(x(t), \prm _0) + h Q(x(t), h)$. Since $g$
is smooth, we can incorporate $Q$ in the perturbative map $\pp$ of a
model like \eqref{model} to satisfy the assumptions. Note that here $h$ also becomes a perturbative parameter. We could treat the two smallness parameters $\ep$ (for $\pp$) and $h$ (for $Q$) jointly by $\ep$ or separately in the fixed point proof of an operator $\op$.

\section{Main ingredients of the proofs}
\label{sec.proofstrategy}

The proof involves several steps; some of them are standard bounds but
others are strongly related to the type of perturbations we
consider. We will start providing a general overview of the tools and
steps of the proposed proof. In particular, we will provide a sequence
of lemmas, which build up the whole proof.

\subsection{Overview of fixed point arguments}
\label{overview} 
Here we give some ideas on the fixed point theorems used and their
variants. Following the proof strategy of center manifold theorem in
\cite{Lan}, to obtain the conclusions it suffices to show that the
operator $\op$ satisfies two types of bounds:
\begin{enumerate}
\renewcommand{\theenumi}{B\arabic{enumi}}
\renewcommand{\labelenumi}{\sl\theenumi.} 
 \item \label{Bpropageted} Propagated bounds, see
   Section~\ref{sec.propagatedbnds};
 \item \label{Bconstraction} Low regularity contraction, see
   Section~\ref{sec.c0contraction}.
\end{enumerate}

The propagated bounds establish that a ball $B$ in a space of smooth
functions is mapped to itself by $\op$. If moreover the operator is a
contraction in a low regularity norm on $B$, then we conclude that
there exists a unique fixed point in the low regularity closure of the
smooth ball (we will denote this by $\overline B$).

The desired result of existence and uniqueness of fixed points can be
established by two different arguments.

The first argument is to appeal to a version of Schauder (see
\cite[p. 179]{Brezis10}).
\begin{thm}
  Let $X$ be a Banach space, $C \subset X $ nonempty, closed convex, $K
  \subset C$ compact, $\Gamma\colon C \rightarrow C$ continuous, $\Gamma(C)
  \subset K$. Then, $\Gamma$ has a fixed point in $K$.
\end{thm}

In our applications, $C = K$ is a ball in a space of highly differentiable
functions with domain $\R$ and Lipschitz modulus of continuity in the
highest derivative. The space $X$ is a Banach space for functions
with domain $\R$ equipped with a low regularity norm.  The fact that
$K$ is compact is a consequence of an easy version of Arzela-Ascoli
theorem since $\R$ is separable.  Since $C=K$ is a ball, convexity is
obvious.

The propagated bounds in Section~\ref{sec.propagatedbnds} show that
that $\op(K) \subset K$.

A further simplification is that, since $K$ is compact in the low
regularity topology, to prove continuity of $\op: K \rightarrow K$
it suffices to show that the graph of $\op$ is closed in the low
regularity topology. This is very easy to verify.

Application of the Schauder theorem obtains the existence of fixed points
using only the propagated bounds. 

The low regularity contraction shows that the fixed point $x^* $ is
unique and provides -- as we show below -- with a-posteriori bounds
using interpolation inequalities. 

To obtain uniqueness, we could consider using other arguments
(e.g. using that two fixed points satisfy the invariance equation or
other geometric properties).

The contraction in low regularity norm has other consequences besides
the uniqueness of the fixed point.

Given a point $x_0 \in K$ we obtain that $\Gamma^n(x_0)$ converges
exponentially fast to $x^*$ in the low regularity norm.  Furthermore,
because of the propagated bounds, the smooth distance between
$\Gamma^n(x_0)$ and $x^*$ remains bounded.  Using interpolation
inequalities
\cite{Had98,Kol49,LO99}, we also obtain exponential convergence of
$\Gamma^n(x_0)$ to $x^*$ in spaces of regularity in between. This leads to an
a-posteriori result estimating the distance between $x_0$ and $x^*$
based on estimates of $\Gamma(x_0) - x_0$ in spaces of low
regularity. One source of interest is that such estimates for a
numerical approximation $x_0$ can be obtained using a computer assisted
proof.  In our case, there are some extra complications since some of
the norms we use are weighted norms. See
Section~\ref{sec.aposteriori}.

A second method of proof used very often in the theory of center
manifolds is to use the theorem in \cite{Lan}.  The method in
\cite{Lan} uses at the same time the propagated bounds and some other
argument to produce uniqueness of fixed points, It can be applied even
when we are interested in functions whose domain is a non-separable
space (so Arzela-Ascoli requires adaptation).

\medskip

There are many possible variants. For the low regularity contraction
we have several choices. We can use weighted norms (we have used the
Razumikhin norms \eqref{eq.etanorm} in some cases) or contractions in
any bounded interval.  The only role is to get uniqueness.

Note that we need to verify the contraction property in functions
which we already know that are smooth. A notable case which appears a
lot in state dependent delays is the composition operator. Note that
we can use $\| u_1\circ u_2 - u_1 \circ \tilde u_2 \|_{C^0} \le \|\dop
u_1 \|_{C^0} \| u_2 - \tilde u_2\|_{C^0}$ if the functions are defined
in a convex set, or more generally in a balanced domain -- i.e. a
domain in which a multiple of the distance among two points bounds
from above the length of the shortest path joining them.

 \smallskip

\begin{rmk}[Improving the $C^{r+\Lip}$ regularity in the conclusions to $C^{r}$]
In this paper we formulate existence results in $C^{r+\Lip}$ spaces to obtain
solutions in the same regularity space.

  This has the minor inconvenience that one has to present extra
  arguments for the last Lipschitz regularity.  The Lipschitz
  constants, in general do not satisfy formulas such as Fa\`a di
  Bruno. When we have an Euclidean domain, the Lipschitz constant can
  be approximated as limit. Therefore, we could prove the result for
  $C ^{r}$ and add a limit argument for the Lipschitz constant in the
  last regularity level.
\end{rmk}

\begin{rmk}[Weaker alternative uniqueness result]
The approach we adopted for Theorem~\ref{main} also admits a weaker
version. The propagated bounds \eqref{Bpropageted} tells us that there
is $(\vft ^\ast, \hx ^\ast)$ in $\X$ such that $\op [\vft ^\ast,
  \hx^\ast] = (\vft ^\ast, \hx ^\ast)$.

The low regularity contraction step, \eqref{Bconstraction}, provides
the uniqueness in such a ball.  It needs to prove the contraction for
all pair of elements in $\X$. However, known already the existence we
can proceed by contradiction and only check the set of possible
fixed points in $\X$. That is, if $(\vft ^\ast, \hx ^\ast)$ and
$(\vfty ^\ast, \hy ^\ast)$ were two different solutions, if we prove
there is $\kappa \in (0,1)$ such that
\begin{equation} \label{eq.weakcontract}
 d\bigl((\vft ^\ast, \hx ^\ast) , (\vfty ^\ast, \hy ^\ast) \bigr )
 \leq \kappa \phantom{\cdot} d\bigl((\vft ^\ast, \hx ^\ast) , (\vfty
 ^\ast, \hy ^\ast) \bigr ),
\end{equation}
for a suitable distance $d(\cdot,\cdot)$, then the solution is
unique. Notice that $\kappa$ can depend on $\ep$ and also that the
inequality in \eqref{eq.weakcontract} does not need to be
strict. This argument is indeed weaker since
it is does not say anything about other elements in the ball and thus it
does not allow an a-posteriori formulation.
\end{rmk}

\subsection{Estimates on evolution} \label{sec:estivol}

The evolution of a one-dimensional vector field can be estimated in a
completely elementary manner. Even if they are elementary, we collect
the estimates in Lemma~\ref{lem.vfregular} for the ease of reference.
Note that we cannot claim that $\hatt\in C^{\ell}$ because our
definition of $C^{\ell}$ spaces involves uniform boundedness (in
particular even the identity map is not $C^{\ell}$ in our
definition). Even if $\Id \notin C^{\ell} $ and $\hatt\notin
C^{\ell}$, we can ``summarize'' the lemma saying that
\begin{equation*}\label{affinespace}
  \vft-1 \in C^{\ell} \Rightarrow \hatt -\Id \in C^{\ell+1}. 
\end{equation*}
Moreover, we have that $\hx \in C^{\ell+1} \Rightarrow \hx \circ \hatt
\in C ^{\ell+1}$. For higher dimensional vector fields, the estimates
are not so strong and, in fact, even for bounded vector fields the
flows can have exponential growth.
This is a reason why in this work we can only deal with hyperbolic
orbits and not with Normally Hyperbolic Invariant Manifolds.

\begin{lem} \label{lem.vfregular}
 Let $\vft$ be a $C^{\ell+\Lip}$ vector field in $\R$ and $\hatt$ be
 its associated evolution given by $\dot \hatt (t) = \vft \circ
 \hatt(t)$ with initial condition $\hatt(0) = 0$. Define $\hat \vft
 \bydef \vft-1$ and assume that $ \hat \vft \in \ball^{\ell+\Lip}
 _{\bt}(0)$ with $\bt \bydef (\bt _0, \bt _1,\dotsc, \bt _{\ell}, \bt
 _{\ell}^{\Lip}) \in \R _+^{\ell+2}.$ If $\bt _0 < 1$, then
 \begin{enumerate}
  \item $\hatt$  and $\hatt^{-1}$ are strictly 
    increasing functions.
  
  \item For all $t$ and $s$ in $\R$, 
  \begin{align}\label{distortion}
  \begin{split}
     (1 - \bt _0) |t-s| \leq |\hatt(t) - \hatt(s) | \leq (1 + \bt _0) |t - s|,\\
      \frac{1}{1+\bt _0} |t-s| \leq |\hatt^{-1}(t) - \hatt^{-1}(s) | \leq
   \frac{1}{1-\bt _0} |t - s|. 
  \end{split}
  \end{align} 
   
  In particular, $(1 - \bt _0) |t| \leq |\hatt(t)| \leq (1 + \bt _0)
  |t|$ and $\frac{1}{1+\bt _0} |t| \leq |\hatt^{-1}(t)| \leq
  \frac{1}{1-\bt _0} |t|$.
  
  \item \label{lem.vfregularIII} $|\dop ^{j+1} \hatt(t)| \leq \tilde
    \bt _{j}$ and $|\dop ^{j+1} (\hatt^{-1})(t)| \leq \hat \bt _{j}$
    for all $j = 0, \dotsc, \ell$, where $\tilde \bt _j$ and $\hat \bt
    _j$ only depend on $\bt _0, \bt _1, \dotsc, \bt _{j}$.
  
  \item $\Lip (\dop ^{\ell+1} \hatt )\leq \tilde \bt _{\ell} ^{\Lip}$
    and $\Lip (\dop ^{\ell+1} (\hatt^{-1})) \leq \hat \bt _{\ell}
    ^{\Lip}$, where $\tilde \bt _{\ell} ^{\Lip}$ and $\hat \bt _{\ell}
    ^{\Lip}$ only depend on $\bt$.
 \end{enumerate}
 In particular, $\dop \hatt \in \ball ^{\ell+\Lip}_{\tilde \bt}$ and
 $\dop (\hatt^{-1}) \in \ball ^{\ell+\Lip}_{\hat \bt}$.
\end{lem}
\begin{proof}
 Let us first observe that $\dop (\hatt^{-1})(t) =
 \frac{1}{X(t)}$, hence $\hatt^{-1}(t) = \int _0 ^t \frac{d
   \sigma}{\vft (\sigma)}$.
 \begin{enumerate}
  \item Since $\bt _0 < 1$, we have $\vft > 0$, which implies the monotonicity.
  \item Note that
  \[
  \hatt(t) - \hatt(s) = \int _s^t \vft \circ \hatt(u)\, du .
  \]
  By the assumption on $\hat \vft$, we obtain $1-\bt _0\leq X \leq
  1+\bt _0$, then the first inequality in \eqref{distortion} is
  proved. We can prove the second inequality for $\hatt^{-1}$
  similarly. The last argument is true since $\hatt(0)=0$.
  \item Clearly, $|\dop\hatt(t)| \leq 1 + \bt _0$. We now prove that
    $|\dop ^j(\dop \hatt(t))|\leq \tilde \bt_j$ for $j = 1, \dotsc,
    \ell$.
  \begin{enumerate}
  \renewcommand{\theenumii}{\roman{enumii}}
  \renewcommand{\labelenumii}{\sl\theenumii)}
   \item $\dop ^2\hatt(t) = \dop \hat \vft\circ \hatt(t) (1 + \hat
     \vft \circ \hatt(t))$ which is bounded by $\bt _ 1(1 + \bt _0)$.
   \item By the Fa\`a di Bruno Formula, we have an expression of the form
   \[
   \dop ^{r+1} \hatt = \dop ^r(\hat \vft \circ \hatt) = \sum
   _{\substack{(m _1, \dotsc, m _r)\in \N ^r \\ \sum\limits_{j=1}^r j
       m _j = r}} C _{m _1, \dotsc, m _r} \dop ^{m _1 + \dotsb + m _r}
   \hat \vft \circ \hatt \prod _{j=1}^r (\dop ^j \hatt)^{m _j},
   \]
   where $C _{m _1, \dotsc, m _r}$'s are  combinatorial numbers. By the 
   induction hypothesis $|\dop ^{j+1}\hatt(t)| \leq \tilde \bt _j$ for
   $j = 0, \dotsc, r-1$ and triangle inequality, we prove that $\dop ^{r+1}\hatt$
   is bounded by $\tilde \bt _r$  depending on $\bt _0, \bt _1, \dotsc, \bt
   _{r}$.
  \end{enumerate}
Since
  \begin{equation}
        \dop ^{r+1} (\hatt^{-1}) = \dop ^r\biggl( \frac{1}{\vft} \biggr),
  \end{equation}
   we derive that $\dop ^{r+1}
   (\hatt^{-1})(t) $ is bounded by a $ \tilde \bt _{r} $ only depending
   on $\bt _0, \bt _1, \dotsc, \bt _{r}$.
   \item Straightforward by using Lemma~\ref{lem.Lipschitz}. \qedhere
 \end{enumerate}
\end{proof}

\subsection{Propagated bounds} \label{sec.propagatedbnds}

In this section, we show that for small enough $\ep$, we can find the
parameters $\bt$, $\bs$, and $\bu$ of the space $\X$, see
\eqref{eq.ballspaxie4op} and \eqref{eq.tsu-components}, so that if the
inputs $(\vft, \hx ^s, \hx ^u)$ are in the space $\X$, then its image
under $\op$ also lies in $\X$.

To reach this goal, we bound $\op[\vft, \hx ^s, \hx ^u]$ and its
derivatives by algebraic expressions of $\bt$, $\bs$, and $\bu$ in the
following Lemmas. We use basic tools like triangle inequalities and
rules of differentiation, including the Leibnitz product formula,
Fa\`a di Bruno formula, etc.

At the end of the section, we discuss the choices for the parameters
$\bt$, $\bs$, and $\bu$ for small enough $\ep$.

\begin{lem} \label{lem.bndTB}
 Let $\ell \geq 0$ be a fixed integer. There are constants $a _0,
 \dotsc, a _\ell, a_\ell ^{\Lip}$, $b _0, \dotsc, b _\ell, b_\ell
 ^{\Lip}$ such that $T[\xo,\hx]$ and $\B[\hvft,\hx]$ defined in
 \eqref{eq.T} and \eqref{Bdefined} have the following bounds for all
 $j = 0, \dotsc, \ell$.
 \begin{enumerate}  
   \item $|\dop ^j T[\xo, \hx](\rho)| \leq a _j$ for $a _j$'s
     depending on $\bs _0, \dotsc, \bs _j$, $\bu _0, \dotsc, \bu _j$,
     $\| \xo \| _{C ^j}$, and $\| \vf \| _{C ^{j+2}}$;
      
  \item $|\dop ^j \B[\hvft,\hx](\rho)| \leq b _j$ for $b _j$'s
    depending on $\bt _0, \dotsc, \bt _j$, $\bs _0, \dotsc, \bs _j$,
    $\bu _0, \dotsc, \bu _j$, $\| \xo \| _{C ^j}$, and $\| \vf \| _{C
      ^{j+2}}$;
  
  \item 
   $\Lip (\dop ^{\ell} T[\xo, \hx])\leq a _\ell ^{\Lip} $ and
    $\Lip(\dop ^\ell \B[\hvft,\hx]) \leq b _\ell ^{\Lip}$ for $a _\ell
    ^{\Lip}$ and $b _\ell ^{\Lip}$ depending on $\bt$, $\bs$, $\bu$,
    $\| \xo \| _{C ^{\ell+\Lip}}$, and $\| \vf \| _{C
    ^{\ell+2+\Lip}}$.
 \end{enumerate}
 Where the norms of $\vf$ are evaluated on a neighborhood of
 $\{\xo(t)\}_{t\in \R}$.
\end{lem}
\begin{proof}
This lemma is proved by using Leibnitz product formula and Fa\`a di
Bruno's formula. By the classical Taylor error bound,
  \begin{equation}\label{Ttaylor}
   T[\xo, \hx](\rho) = \int _0^1 \int^\sigma_0 \dop ^2 \vf\circ (\xo + s \hx)(\rho) \hx(\rho)^2 \, d s\, d\sigma,
  \end{equation}
  where we use $\dop ^2 \vf\circ (\xo + s \hx)(\rho) \hx(\rho)^2 $ to
  denote the result of the bilinear operator $\dop ^2 \vf\circ (\xo +
  s \hx)(\rho) $ acting on $\hx(\rho)$ and $\hx(\rho)$. Then, we can
  bound $\dop ^j T [\xo, \hx]$ in terms of $\bs _0 + \bu _0, \dotsc,
  \bs _j + \bu _j$, $\| \xo \| _{C ^j}$, and $\|\vf\|_{C
    ^{j+2}}$. Similarly, $\Lip (\dop ^{\ell} T[\xo, \hx])$ can be
  bounded by $\bs$, $\bu$, $\| \xo \| _{C ^{\ell+\Lip}}$, and $\| \vf
  \| _{C ^{\ell+2+\Lip}}$.
 
   Since 
   \[\dop ^j \B[\hvft,\hx]=\dop ^j[(1 - \vft) \dop \vf \circ \xo \hx]+\dop ^j T [\xo, \hx],\]
   $\dop ^j \B[\hvft,\hx]$ is bounded by an algebraic expression of
   $\bt _0, \dotsc, \bt _j$, $\bs _0, \dotsc, \bs _j$, $\bu _0,
   \dotsc, \bu _j$, $\| \xo \| _{C ^j}$, and $\| \vf \| _{C ^{j+2}}$.
   \qedhere

\end{proof}

\smallskip

The operator $\op$ defined in Section~\ref{sec.operator} involves
quotients. To bound the center component of $\op$, we need to assume
that the vector field $f$ is bounded away from zero on the unperturbed
hyperbolic solution $\{x_0(t)\}$, i.e. there is $b > 0$ such that $b
\leq \inf \{|\vf\circ x_0 | \}$. For the stable and unstable
components, we use the fact that $1 - \bt _0<\vft$ for
$\vft\in\ball_{\bt}^{\ell+\Lip}(1)$.

We use rules of differentiation, Cauchy-Schwartz inequality, and Lemma
\ref{lem.bndTB} to prove:
\begin{pro}[center correction] \label{prop.b1center}
There are constants $\ct b _{c, j}$ and $\ct d _{c, j}$ such that for
all $\rho\in\R$ and $(\vft, \hx ^s, \hx ^u)\in \X$, the operator $\op
_c$ defined in \eqref{Gammac} satisfies
 \begin{enumerate}
  \item $| \op _c [\vft, \hx ^s, \hx ^u](\rho) - 1 | \leq \ct b _{c,0}
    + \ep \ct d _{c,0}$;
 \item $ |\dop ^j \op _c [\vft, \hx ^s, \hx ^u](\rho) | \leq \ct b
   _{c,j} + \ep\ct d _{c,j} $ and $j = 1, \dotsc, \ell$;
 \item $ \Lip\left(\dop ^\ell \op _c [\vft, \hx ^s, \hx ^u]\right)
   \leq \ct b^{\Lip} _{c,\ell} + \ep\ct d ^{\Lip}_{c,\ell} $,
 \end{enumerate}
 where $\ct b _{c, j}$'s and $\ct d _{c, j}$'s depend on
 $\|\xo\|_{C^j}$, $\|\Pi ^c_t\|_{C^j}$, $\|\vf\|_{C^{j+2}}$, $\bt _0,
 \dotsc, \bt _{j}$, $\bs _0, \dotsc, \bs _{j}$, and $\bu _0, \dotsc,
 \bu _{j}$. The constants $b^{\Lip} _{c,\ell}$ and $\ct d
 ^{\Lip}_{c,\ell} $ depend on $\|\xo\|_{C^{\ell+\Lip}}$, $\|\Pi
 ^c_t\|_{C^{\ell+\Lip}}$, $\|\vf\|_{C^{\ell+2+\Lip}}$, $\bt$, $\bs$,
 and $\bu$.
\end{pro}
\begin{proof}
   For the first bound, we use \eqref{Hep1} and Lemma \ref{lem.bndTB},
   and note that the terms involving $\B[\hvft,\hx]$ and
   $\varphi[\vft,\hx]$ are bounded by $\ct b _{c,0}$ and $\ep \ct d
   _{c, 0}$ respectively, where
\[
 \ct b _{c,0} = \frac{C _\Pi\|\vf\|_{C^0}}{b^2} \biggl[\bt _0 \| \vf
   \| _{C ^1}(\bs _0 + \bu _0) + \frac{1}{2} \|\vf \|_{C ^2} (\bs _0 +
   \bu _0 )^2\biggr] \qquad \ct d _{c, 0} = \frac{C
   _\Pi\|\vf\|_{C^0}}{b^2} \|\varphi[\vft,\hx]\| _{C ^0}.
\]
   Then we use \eqref{Hep1}, Fa\`a di Bruno formula, Leibnitz product
   formula, the quotient rule, Lemma~\ref{lem.vfregular}, and
   Lemma~\ref{lem.bndTB} to obtain the bounds for the derivatives of
   $\op _c$.
\end{proof}
The operators $\op _s$ and $\op _u$ in \eqref{Gammasu} gain one
derivative thanks to the integration, which is the reason why we
define the space $\X$ in \eqref{eq.ballspaxie4op} with different
regularities in the components. Here we bound the derivatives of $\op
_s$ and $\op _u$ up to order $\ell+1$.

\begin{pro}[stable and unstable corrections] \label{prop.b1stabunstab}
There are constants $\ct b _{\sigma, k}$ and $\ct d _{\sigma, k}$,
$\sigma \in \{s,u\}$, such that for all $\rho\in\R$ and $(\vft, \hx
^s, \hx ^u)\in \X$, the operators $\op _s$ and $\op _u$ defined in
\eqref{Gammasu} satisfy
 \begin{enumerate}
  \item  
  $
   | \op _\sigma [\vft, \hx ^s, \hx ^u](\rho)  | \leq \ct b _{\sigma,0} + \ep \ct d _{\sigma,0}
  $.
 \item $ |\dop ^{j} \op _\sigma [\vft, \hx ^s, \hx ^u](\rho) | \leq
   \ct b _{\sigma,j} + \ep \ct d _{\sigma,j} $ for $j = 1, \dotsc,
   \ell+1$;
\item $ \Lip\left(\dop ^{\ell+1} \op _\sigma [\vft, \hx ^s, \hx ^u]
  \right) \leq \ct b _{\sigma,\ell+1}^{\Lip} + \ep \ct d
  _{\sigma,\ell+1}^{\Lip} $,
 \end{enumerate}
 where for $j>0$, $\ct b _{\sigma, j}$'s and $\ct d _{\sigma, j}$'s
 depend on $\|\xo\|_{C^{j-1}}$, $\lambda _\sigma$, $\|\Pi
 ^\sigma_t\|_{C^{j-1}}$, $\|\vf\|_{C^{j+1}}$, $\bt _0, \dotsc, \bt
 _{j-1}$, $\bs _0, \dotsc, \bs _{j-1}$, and $\bu _0, \dotsc, \bu
 _{j-1}$. Moreover, $\ct d _{\sigma, j}$'s also depend on $ \bs _j$
 and $\bu _j$. Similarly, the constants $b _{\sigma,\ell+1}^{\Lip}$ and
 $\ct d _{\sigma,\ell+1}^{\Lip}$ depend on $\|\xo\|_{C^{\ell+\Lip}}$,
 $\lambda _\sigma$, $\|\Pi ^\sigma_t\|_{C^{\ell+\Lip}}$,
 $\|\vf\|_{C^{\ell+2+\Lip}}$, $\bt$, $\bs$, and $\bu$, where the dependence
 on $ \bs ^{\Lip} _{\ell+1}$ and $\bu^{\Lip} _{\ell+1}$ is only for $\ct d
 _{\sigma,\ell+1}^{\Lip}$.
\end{pro}
\begin{proof}
Note that
 \[\ct b _{\sigma,0} = \frac{C _\Pi C _U}{\lambda _\sigma (1 - \bt _0)} \biggl[\bt _0 \| \vf \| _{C ^1}(\bs _0 + \bu _0) + \frac{1}{2} \|\vf \|_{C ^2} (\bs _0 + \bu _0 )^2\biggr] \qquad \ct 
 d _{\sigma, 0} = \frac{C _\Pi C _U}{\lambda _\sigma (1 - \bt _0)} \|
 \varphi[\vft,\hx] \| _{C ^0}.\] In order to bound the derivatives of
 $\op _\sigma$, we use the fact that $\op _\sigma$ solves the
 differential equation
\begin{equation}
\label{eq.opderiv}
    \dop \op _\sigma [\vft, \hx ^s, \hx ^u](\rho) = \dop \vf \circ
    \xo(\rho) \op _\sigma [\vft, \hx ^s, \hx ^u](\rho) + \Pi ^\sigma
    _\rho \frac{1}{\vft(\rho)} \Bigl[\B[\vft,\hx](\rho) + \ep
      \varphi[\vft,\hx](\rho) \Bigr], \quad \sigma \in \{s,u\}.
\end{equation} 
 Then we can obtain the bounds with Fa\`a di Bruno formula, Leibnitz
 product formula, the quotient rule, Lemma~\ref{lem.vfregular}, and
 Lemma~\ref{lem.bndTB}.
\end{proof}

It remains to show that it is possible to choose the components of
$\bt$, $\bs$, and $\bu$ so that $\op$ maps $\X$ into itself as long as
$\ep$ is small enough. As we will see, we need to choose a small
$\bt_0$. Without loss of generality, we assume that $\bt_0\leq
\frac{1}{2}$ so that $\frac{1}{1-\bt_0}\le 2$ and we do not need to
worry about $1-\bt_0$ in the denominator.
 
 Indeed, the zero order constants should satisfy

\begin{align} \label{eq.bnd0t0}
\ct b _{c,0} + \ep  \ct d _{c,0} &\leq \bt _0,\nonumber\\
\ct b _{s,0} + \ep \ct d _{s,0} &\leq \bs _0,\\
\ct b _{u,0} + \ep \ct d _{u,0} &\leq \bu _0.\nonumber
\end{align}                                                               
where the left sides of the inequalities come from
Proposition~\ref{prop.b1center} and
Proposition~\ref{prop.b1stabunstab}. As $\ct b _{c,0}$, $\ct b _{s,0}$,
and $\ct b _{u,0}$ being quadratic in $\bt_0$, $\bs_0$, and $\bu_0$,
we can choose small enough $\bt_0$, $\bs_0$, and $\bu_0$ so that when
$\ep$ is small enough, the set of inequalities \eqref{eq.bnd0t0} are
satisfied.

For the $i$-th order, the following inequalities should be satisfied.
\begin{align}\label{eq.bndi}
    G_c^i(\ep, \bt_0,\dotsc,\bt_{i}, \bs_0,\dotsc,\bs_{i},\bu_0,\dotsc,\bu_{i})&\leq \bt_{i},\nonumber\\
    G_s^i(\ep, \bt_0,\dotsc,\bt_{i-1}, \bs_0,\dotsc,\bs_{i-1},\bu_0,\dotsc,\bu_{i-1})&\leq \bs_{i},\\
    G_u^i(\ep, \bt_0,\dotsc,\bt_{i-1}, \bs_0,\dotsc,\bs_{i-1},\bu_0,\dotsc,\bu_{i-1})&\leq \bu_{i},\nonumber
\end{align}
where $G_\sigma^i$, $\sigma\in \{c,s,u\}$ are polynomials of $\ep$ and
the components of $\bt$, $\bs$, and $\bu$. Moreover, one factor in the
coefficient of $\bt_{i}$ in $G_c^i$ is $(\bs_0+\bu_0)$. In order to
guarantee \eqref{eq.bndi}, we first fix $\bs_i$ and $\bu_i$, and then
choose $\bt_i$. Similar arguments hold for $\bt^{\Lip}_\ell$, $\bs^{\Lip}_{\ell+1}$, and $\bu^{\Lip}_{\ell+1}$. Indeed, in this process, we may have to ask for
smaller $\ep$, $\bs_0$, and $\bu_0$ at each step, so we will not be
able to obtain $C^\infty$ result with our method in general.

\subsection{Low regularity contraction} \label{sec.c0contraction}
The operator $\op$ defined in Section~\ref{sec.operator} is a contraction on $
\X ^\ell_{\bt, \bs, \bu}$ in \eqref{eq.ballspaxie4op} if there is $\kappa \in (0,1)$ such that
\begin{multline} \label{eq.contraction}
d  \Bigl( \bigl(\op_c[\vft, \hx ^s, \hx ^u],\op _s[\vft, \hx ^s, \hx ^u],\op _u[\vft, \hx ^s, \hx ^u]\bigr), \bigl(\op_c [\vfty, \hy ^s, \hy ^u],\op_s [\vfty, \hy ^s, \hy ^u], \op _u [\vfty, \hy ^s, \hy ^u] \bigr) \Bigr)
 \\
< \kappa \phantom{\cdot} d \Bigl(\bigl( \vft, \hx ^s, \hx ^u \bigr) ,\bigl(\vfty, \hy
^s, \hy ^u \bigr) \Bigr),
\end{multline}
for all $(\vft, \hx ^s, \hx ^u) $ and $(\vfty, \hy ^s, \hy ^u) $ in $
\X ^\ell_{\bt, \bs, \bu}$ and $d(\cdot,\cdot)$ a distance.

\smallskip

We consider the distance in a low regularity space, $ \X ^\ell_{\bt,
  \bs, \bu}$ where $\ell =0$. The space $\X ^0_{\bt, \bs, \bu}$ has
information of the vector field $\vft$, $\hx$, $\dop \hx$. Because we
accept time-dependence in the perturbative map $\pp$, our construction
requires to bound the difference of backward flows associated to the
center correction. More precisely, if $\hatt$ and $\hatty$ are flows
of $\vft$ and $\vfty$ respectively, we need to bound
$ \hatt^{-1}(\rho) - \hatty^{-1}(\rho)$. Neverthelss, this
may fail to be bounded in $C^0$  (e.g. if the vector fields
differ by a constant) but, 
however, it can be  bounded in $C _\eta$ space for $\eta > 0$, see
\eqref{razumikhin-norm}. Therefore, we consider the distance:
\begin{equation*} \label{eq.distance-contraction}
 d \Bigl(\bigl( \vft, \hx ^s, \hx ^u \bigr) ,\bigl(\vfty, \hy
^s, \hy ^u \bigr) \Bigr)
 \bydef \|\vft - \vfty\| _\eta + \|\hx ^s - \hy ^s\| _\eta + \|\hx ^u
 - \hy ^u\| _\eta + \|\dop \hx ^s - \dop \hy ^s\| _\eta + \|\dop \hx ^u
 - \dop \hy ^u\| _\eta,
\end{equation*}
where the norm $\| \cdot \| _{\eta}$ on $\R$ is defined by
\begin{equation*} \label{eq.etanorm}
\| x - y \| _\eta \bydef \sup _{\rho \in \R} |x(\rho) - y(\rho) | e
^{- \eta|\rho|}.
\end{equation*}
In what follows and for typographical reasons, we may just write $\hx$
for $(\hx ^s, \hx ^u)$ to express, for instance, $\| \hx - \hy \|$
instead of $\| \hx ^s - \hy ^s \| + \| \hx ^u - \hy ^u \|$. Similarly
for $\dop \hx$ and $\dop \hy$.

\medskip

Assuming that for $(\vft, \hx) $ and $(\vfty, \hy) $ in $\X$, we have
bounds in the differences involving $\B$ defined in \eqref{Bdefined}
and $\varphi$ defined in \eqref{eq.Px} as in inequalities
\eqref{eq.Bcontract} and \eqref{eq.varphicontract}.  Then if $0< b
\leq \inf \{|\vf\circ x_0 | \}$, we can obtain
\begin{equation}\label{eq.contrac}
 \|\op _c[\vft,\hx] - \op _c[\vfty,\hy]\|_\eta \leq \frac{C _\Pi \|\vf \| _{C ^0}}{b^2}\Bigl[  (\ct d _\B + \ep \ct d _\varphi) \| \vft - \vfty \| _\eta +
 (\ct c _\B + \ep \ct c _\varphi) \| \hx - \hy \| _\eta + \ep \ct e _\varphi \| \dop \hx - \dop \hy \| _\eta \Bigr], 
\end{equation}
and for $\sigma \in \{s,u\}$,
\begin{equation}\label{eq.contras-u}
 \|\op _\sigma[\vft,\hx] - \op _\sigma[\vfty,\hy]\|_\eta \leq \frac{2C _\Pi C _U}{(\lambda _\sigma -\eta ) (1 - \bt _0)}\Bigl[ (\ct d _\B + \ep \ct d _\varphi) \| \vft - \vfty \| _\eta +
 (\ct c _\B + \ep \ct c _\varphi) \| \hx - \hy \| _\eta + \ep \ct e _\varphi \| \dop \hx - \dop \hy \| _\eta \Bigr],
\end{equation}
where $\ct c _\B$, $\ct d _\B$, $\ct c _\varphi$, $\ct d _\varphi$,
and $\ct e _\varphi$ are constants specified in
Propositions~\ref{pro.Bcontract} and
\ref{pro.varphicontract}. Notably, we can make the constants $\ct c
_\B$ and $\ct d _\B$ small. For the stable and unstable components of
$\op$, we have used the bounds in \eqref{eq.expoU}. We discuss the
idea for the stable one. With the Razumikhin norm, we have to bound an
integral of the form:
\begin{equation*}
    I_1=\int^\rho_{-\infty} e^{-\lambda_s(\rho-v)}e^{-\eta |\rho|}e^{\eta|v|} dv.
\end{equation*}
We consider the cases when $\rho\leq 0$ and $\rho> 0$, and obtain that
$I_1\le \frac{2}{\lambda_s-\eta}$ for both cases, provided $\eta <
\lambda _s$. The unstable direction could be estimated similarly when
$\eta < \lambda _u$. Therefore, we derive the estimates in
\eqref{eq.contras-u} when $\eta < \min \{ \lambda _s, \lambda _u\}$.

Using the expression of $\dop \op _\sigma$ for $\sigma \in \{s,u\}$ in
\eqref{eq.opderiv}, the differences of $\dop \op _\sigma$ can be
estimated by
\begin{multline}
\label{eq.contrad}
 \|\dop \op _\sigma[\vft,\hx] - \dop \op _\sigma[\vfty,\hy]\|_\eta \leq \| \vf\| _{C ^1} \| \op _\sigma [\vft, \hx] - \op _\sigma [\vfty ,\hy] \| _\eta \\
 +\frac{C _\Pi}{1 - \bt _0}\Bigl[ (\ct d _\B + \ep \ct d _\varphi) \| \vft - \vfty \| _\eta + (\ct c _\B + \ep \ct c _\varphi) \| \hx - \hy \| _\eta + \ep \ct e _\varphi \| \dop \hx - \dop \hy \| _\eta \Bigr].
\end{multline}
With the bounds in \eqref{eq.contrac}--\eqref{eq.contrad}, we prove
that the operator is a contraction if $\ep$ is small enough.

\medskip

We now bound the differences in $\B$ and $\varphi$ for $(\vft, \hx) $
and $(\vfty, \hy)$ in Propositions~\ref{pro.Bcontract}
and~\ref{pro.varphicontract}.

\begin{pro} \label{pro.Bcontract}
 Let $\B$ be the map defined in \eqref{Bdefined} and let $\vft$,
 $\vfty$ be in $\ball _{\bt _0}^0(1)$ and let $\hx$, $\hy$ be in
 $\ball_{\bs _0 + \bu _0}^0(0)$. Then there are constants $\ct c _\B$
 and $\ct d _\B$ only depending on $\|\vf \| _{C ^{2+\Lip}}$, $\bt
 _0$, $\bs _0$, and $\bu _0$ such that
 \begin{equation}\label{eq.Bcontract}
      |\B[\vft, \hx](\rho) - \B[\vfty,\hy](\rho)| e ^{-\eta|\rho|} \leq \ct c _\B \| \hx - \hy \| _\eta + \ct d _\B \|\vft - \vfty \| _\eta.
 \end{equation}
\end{pro}
\begin{proof}
 $\B$ consists of two terms. The one coming from the Taylor error is
  bounded using the integral formulation, 
 \[ T[\xo, \hx](\rho) = \int _0^1 \int^\sigma_0 \dop ^2 \vf\circ (\xo + s \hx)(\rho) \hx(\rho)^2 \, d s\, d\sigma. \]
 Hence by adding and subtracting
 \[
  |T[\xo,\hx](\rho) - T[\xo,\hy](\rho)| e ^{-\eta|\rho|}\leq   (\bs _0 + \bu _0)\left(\Lip(\dop ^2 \vf) (\bs _0 + \bu _0)+ \|\dop ^2 \vf\| \right)\| \hx - \hy \| _\eta.
 \]
 The other term in $\B$ is also bounded similarly, by adding and
 subtracting, which ends up to the final bound
\begin{align*}
    |\B[\vft, \hx](\rho) - \B[\vfty,\hy](\rho)| e ^{-\eta|\rho|}&\leq   \biggl(\| \vf \| _{C ^1}\bt _0 +(\bs _0 + \bu _0)\left(\Lip(\dop ^2 \vf) (\bs _0 + \bu _0)+ \|\dop ^2 \vf\| \right)\biggr) \| \hx - \hy \| _\eta\\
    &\qquad +\| \vf \| _{C ^1}(\bs _0 + \bu _0) \| \vft - \vfty \| _\eta.
\end{align*}
Defining
\begin{align*}
    \ct c _\B & \bydef \| \vf \| _{C ^1}\bt _0 +(\bs _0 + \bu _0)\left(\Lip(\dop ^2 \vf) (\bs _0 + \bu _0)+ \|\dop ^2 \vf\| \right),\\
    \ct d _\B & \bydef \| \vf \| _{C ^1}(\bs _0 + \bu _0),
\end{align*}
we have the desired inequality. Moreover, the constants $\ct c _\B$
and $\ct d _\B$ are small if $\bt _0$, $\bs _0$, and $\bu _0$ are
small.  \qedhere
\end{proof}

\smallskip
\begin{rmk}
    Notice that the smallness of $ \ct c _\B$ and $ \ct d _\B$ is
    ensured by choosing small enough $\bt_0$, $\bs _0$, and $\bu _0$.
\end{rmk}
To bound the difference in $\varphi$, we prove two preliminary lemmas.
Lemma~\ref{lem.phipsi-contract} shows how to bound the difference of
two backward flows, which motivates our choice of Razumikhin norm.
Lemma~\ref{lem.alphbeta-contract} bounds difference of two forward
flows composed with backward ones. The last result is essential for
the type of functional perturbations we are interested in.

\begin{lem} \label{lem.phipsi-contract}
 Let $\vft$ and $ \vfty$ be vector fields in $\R$ in a ball $\ball
 ^0_{\bt _0}(1)$ with $\bt _0 \in (0,1)$ and let $\eta > 0$. If $\dot
 \hatt = \vft \circ \hatt$ and $\dot \hatty = \vfty \circ \hatty$ with
 zero initial conditions at zero, then
 \begin{equation*}
  | \hatt^{-1}(\rho) - \hatty^{-1}(\rho) | e ^{-\eta |\rho|} \leq \frac{\|\vft - \vfty \| _\eta}{\eta (1 - \bt _0)^2}.
 \end{equation*}
 In particular, $\| \hatt^{-1} - \hatty^{-1} \| _\eta \leq
 \tfrac{1}{\eta (1 - \bt _0)^2} \|\vft - \vfty \| _{C^0}$.
\end{lem}
\begin{proof}
 Since $\hatt(0) = \hatty(0)=0$,
 \begin{equation*}
  \hatt ^{-1} (\rho) = \int _0^\rho \frac{d\sigma }{\vft (\sigma)} \quad \text{and} \quad
  \hatty ^{-1} (\rho) = \int _0^\rho \frac{d\sigma }{\vfty (\sigma)}.
 \end{equation*}
 Therefore
 \begin{equation*}
   | \hatt^{-1}(\rho) - \hatty^{-1}(\rho) | e ^{-\eta |\rho|} \leq \frac{\|\vft - \vfty \| _\eta}{(1 - \bt _0)^2} \int _0 ^1 |\rho | e ^{(\sigma - 1)\eta|\rho|}\, d\sigma =  \frac{\|\vft - \vfty \| _\eta}{\eta(1 - \bt _0)^2}(1 - e ^{-\eta |\rho|}) \leq \frac{\|\vft - \vfty \| _\eta}{\eta(1 - \bt _0)^2}. \qedhere
 \end{equation*}
\end{proof}

\begin{lem} \label{lem.alphbeta-contract}
 Let $\hatt$, $\hatty$ be flows of vector
 fields $\vft$, $\vfty \in \ball ^{\Lip}_{(\bt _0,\bt _1)}(1)$ respectively with zero initial conditions at zero. For all
 $s \in [-h,h]$, define
 \begin{equation*}
  \alpha(\rho,s) \bydef \hatt(\hatt^{-1}(\rho) + s) \quad \text{and} \quad
  \beta(\rho,s) \bydef \hatty(\hatty^{-1}(\rho) + s).
 \end{equation*}
 Then there is a constant $\ct z$ depending on $\bt _0, \bt _1$, $\eta$, $h$ such that
 \begin{equation*}
  \sup _{s \in [-h,h]}|\alpha(\rho,s) - \beta(\rho,s)| e ^{-\eta
    |\rho|} \leq \ct z \|\vft - \vfty \| _\eta .
 \end{equation*}
\end{lem}

\begin{proof}
In order to consider different signs of $s \in [-h,h]$, we define
$\alpha _\pm(\rho, s) \bydef \alpha (\rho, \pm s)$ and $\beta
_\pm(\rho, s) \bydef \beta(\rho, \pm s)$ for $s \in [0,h]$. Then
 \begin{equation*}
 \begin{split}
 \sup _{s \in [-h,h]}|\alpha(\rho,s) - \beta(\rho,s)| e ^{-\eta
    |\rho|} = \max \Biggl\{ &\sup _{s \in [0,h]}|\alpha_+(\rho,s) - \beta_+(\rho,s)| e ^{-\eta
    |\rho|}, \\
    &\sup _{s \in [0,h]}|\alpha_-(\rho,s) - \beta_-(\rho,s)| e ^{-\eta
    |\rho|}\Biggr\}.
\end{split}
\end{equation*}
 By expanding in $s$,
 \begin{equation*}
  \alpha_\pm(\rho, s) = \rho \pm \int _0^s \vft\circ \alpha _\pm(\rho, \sigma)  \, d\sigma \qquad \text{and} \qquad 
  \beta _\pm(\rho, s) = \rho \pm \int _0^s \vfty\circ \beta_\pm (\rho, \sigma)\, d\sigma.
 \end{equation*}
 Adding and subtracting,
 \begin{equation*}
  |\alpha_\pm( \rho,s) - \beta _\pm(\rho,s)| e ^{-\eta |\rho|} \leq \int _0^s
  |\vft\circ \alpha _\pm(\rho,\sigma ) - \vfty\circ \alpha _\pm(\rho,\sigma )|e
  ^{-\eta |\rho|} +\Lip(\vfty) |\alpha _\pm(\rho,\sigma ) -
  \beta _\pm(\rho,\sigma )|e ^{-\eta |\rho|} \, d\sigma.
 \end{equation*}
 Notice that
 \begin{align*}
  \int _0^s |\vft\circ \alpha _\pm(\rho,\sigma ) - \vfty\circ
  \alpha _\pm(\rho,\sigma )|e ^{-\eta |\rho|} d \sigma &\leq \| \vft - \vfty \| _\eta
  \int _0^s e ^{\eta \bigl[ |\alpha _\pm (\rho, \sigma )| - |\rho| \bigr]}\,
  d \sigma \\
  &\leq \| \vft - \vfty \| _\eta \int _0^s e ^{\eta (1 + \bt
    _0) \sigma } \, d \sigma =  
    \| \vft - \vfty \| _\eta \frac{e ^{\eta (1 + \bt _0) h} - 1 }{\eta (1 + \bt _0)}.
 \end{align*}
 By Gr\"onwall's inequality, 
 \begin{equation*}
  |\alpha _\pm(\rho,s) - \beta _\pm(\rho,s)| e ^{-\eta |\rho|} \leq e ^{\bt _1
    h} \frac{e ^{\eta (1 + \bt _0) h} - 1 }{\eta (1 + \bt _0)}\| \vft - \vfty \| _\eta  .\qedhere
 \end{equation*}
\end{proof}

\begin{pro} \label{pro.varphicontract}
There are constants $\ct c _\varphi$, $\ct d _\varphi$, and $\ct e
_\varphi$ such that for all $\vft$, $\vfty \in \ball _{(\bt _0,\bt
  _1)}^{\Lip}(1)$ and $\hx$, $\hy \in \ball _{(\bs_0 + \bu _0,\bs _1 +
  \bu _1,\bs _2 + \bu _2)}^{1+\Lip}(0)$, the following inequality
holds for the map $\varphi$ defined in \eqref{eq.Px}.
\begin{equation}\label{eq.varphicontract}
    |\varphi[\vft,\hx](\rho) - \varphi[\vfty,\hy](\rho)| e
  ^{-\eta|\rho|} \leq \ct c _\varphi \| \hx - \hy \| _\eta + \ct d
  _\varphi \|\vft - \vfty \| _\eta  + \ct e
  _\varphi \|\dop \hx - \dop \hy \| _\eta.
\end{equation}
  
\end{pro}
\begin{proof}
 By the definition of $\varphi$ and the assumption \eqref{Hep2}, we
 have that
 \[
  |\varphi[\vft,\hx](\rho) - \varphi[\vfty,\hy](\rho)| \leq \ct L
  _1 | \hatt^{-1}(\rho) - \hatty^{-1}(\rho)| + \ct L _2 \Bigl\| \bigl((\xo
  + \hx)\circ \hatt \bigr)_{\hatt^{-1}(\rho)} - \bigl((\xo + \hy)\circ
  \hatty \bigr)_{\hatty^{-1}(\rho)} \Bigl\| _ {C ^1([-h,h])}.
 \]
Using Lemma~\ref{lem.phipsi-contract}, $| \hatt^{-1}(\rho) -
\hatty^{-1}(\rho)|e ^{-\eta |\rho|}$ is bounded by a constant multiple
of $\|\vft - \vfty \| _\eta$.  In order to bound the second part of
the above inequality, we first consider
 \begin{equation} \label{C0-4-Hep2}
 \sup _{s \in [-h,h]} |(\xo + \hx)\circ \hatt (\hatt^{-1}(\rho)+s) -
 (\xo + \hy)\circ \hatty (\hatty^{-1}(\rho)+s)| e ^{-\eta |\rho|}.
 \end{equation}
 

 Using the $\alpha$, $\beta$ notation from
 Lemma~\ref{lem.alphbeta-contract} and adding/subtracting,  \eqref{C0-4-Hep2} is equivalent to
 \begin{equation*}
  \sup _{s \in [-h,h]}\bigl|(\xo \circ \alpha - \xo \circ \beta) + 
  (\hx \circ \alpha - \hy \circ \alpha) + 
  (\hy \circ \alpha - \hy \circ \beta)\big|_{(\rho,s)}\bigr| e^{-\eta |\rho|}.
 \end{equation*}
 The first and third terms are bounded using
 Lemma~\ref{lem.alphbeta-contract} and Lipschitz property of $\xo$ and
 $\hy$. The second term is controlled as follows
 \[
 |\hx \circ \alpha(\rho,s) - \hy \circ \alpha(\rho, s)| e ^{-\eta
   |\rho|} \leq e ^{\eta (|\alpha(\rho, s)| - |\rho|)}\| \hx - \hy \|
 _\eta \leq e ^{\eta (1 + \bt _0) h} \| \hx - \hy \| _\eta.
 \]
 Now we consider the derivative 
 \begin{multline*}
  \frac{d}{d s} \Bigl[ \bigl((\xo
  + \hx)\circ \hatt \bigr)_{\hatt^{-1}(\rho)}(s) - \bigl((\xo + \hy)\circ
  \hatty \bigr)_{\hatty^{-1}(\rho)}(s) \Bigr] 
  \\
  =(\xo + \hx)'\circ \alpha(\rho, s) \vft \circ \alpha (\rho, s) - (\xo + \hy)'\circ \beta(\rho, s) \vfty \circ \beta (\rho, s),
 \end{multline*}
which equals to the following sum evaluated at $(\rho,s)$ by
adding/subtracting
 \begin{align}
  &(\xo' \circ \alpha - \xo' \circ \beta) \vft \circ \alpha  \tag{\texttt{L1}} \label{pro.varphicontract.L1} \\
  &+\xo'\circ \beta (\vft \circ \alpha - \vfty \circ \alpha)  \tag{\texttt{L2}} \label{pro.varphicontract.L2} \\
  &+\xo'\circ \beta( \vfty \circ \alpha - \vfty \circ \beta)  \tag{\texttt{L3}} \label{pro.varphicontract.L3} \\
  &+(\hx' \circ \alpha - \hx'\circ \beta) \vft \circ \alpha   \tag{\texttt{L4}} \label{pro.varphicontract.L4} \\
  &+(\hx'\circ \beta - \hy'\circ \beta )\vft\circ \alpha   \tag{\texttt{L5}} \label{pro.varphicontract.L5} \\
  &+\hy'\circ \beta (\vft \circ \alpha - \vft \circ \beta)   \tag{\texttt{L6}}  \label{pro.varphicontract.L6} \\
  &+\hy'\circ \beta (\vft \circ \beta - \vfty \circ \beta).  \tag{\texttt{L7}}  \label{pro.varphicontract.L7}
 \end{align}
Each line can be bounded directly or by using
Lemma~\ref{lem.alphbeta-contract}. Indeed,
\begin{align*}
 |\eqref{pro.varphicontract.L1} (\rho,s)| e ^{-\eta |\rho|} &\leq (1+\bt _0) \Lip(\xo ' ) \ct z \|\vft - \vfty \| _\eta, \\
 |\eqref{pro.varphicontract.L2} (\rho,s)| e ^{-\eta |\rho|} &\leq \| \xo \| _{C ^1} e ^{\eta(1+\bt _0)h} \| \vft - \vfty \| _\eta, \\
 |\eqref{pro.varphicontract.L3} (\rho,s)| e ^{-\eta |\rho|} &\leq \| \xo \| _{C ^1} \bt _1  \ct z  \| \vft - \vfty \| _\eta, \\
 |\eqref{pro.varphicontract.L4} (\rho,s)| e ^{-\eta |\rho|} &\leq (1 + \bt _0) (\bs _2 + \bu _2) \ct z  \| \vft - \vfty \| _\eta, \\
 |\eqref{pro.varphicontract.L5} (\rho,s)| e ^{-\eta |\rho|} &\leq (1+\bt _0) e ^{\eta (1 + \bt _0) h} \| \hx' - \hy' \| _\eta, \\
 |\eqref{pro.varphicontract.L6} (\rho,s)| e ^{-\eta |\rho|} &\leq (\bs _1 + \bu _1) \bt _1  \ct z  \| \vft - \vfty \| _\eta, \\
 |\eqref{pro.varphicontract.L7} (\rho,s)| e ^{-\eta |\rho|} &\leq (\bs _1 + \bu _1)  e^{\eta (1 + \bt _0)h} \| \vft - \vfty \| _\eta.
\end{align*} 
 
 Collecting all the intermediate bounds we have explicit $\ct c
 _\varphi$, $\ct d _\varphi$, $\ct e _\varphi$ depending on
 $\|\xo\|_{C^{1+\Lip}}$, $\ct L _1$, $\ct L _2$, $h$, $\eta$, $\bt
 _0$, $\bt _1$, $\bs _1$, $\bs _2$, $\bs _1$, and $\bu _2$. 
\end{proof}

\subsection{A-posteriori results} \label{sec.aposteriori}
By the propagated bounds \eqref{Bpropageted}, there is a fixed point $
v^\ast$ of the operator $\op$. For an initial guess $v\bydef (\vft
_{(0)}, \hx ^s _{(0)}, \hx ^u_{(0)})$ of the fixed point method, we
have
\begin{equation*}
 \begin{split}
    \| \vft _{(0)} - \Pi ^c [v ^\ast] \| _{C ^{\ell+\Lip}} &\leq M _c <+\infty ,\\
    \| \hx ^s _{(0)} - \Pi ^s [v ^\ast] \| _{C ^{\ell+1+\Lip}} &\leq M _s <+ \infty ,\\
     \| \hx ^u _{(0)} - \Pi ^u [v ^\ast] \| _{C ^{\ell+1+\Lip}}  &\leq M _u < +\infty .
 \end{split}
\end{equation*}

On the other hand, from the low regularity contraction
\eqref{eq.contraction} and the Banach fixed point Theorem,
\begin{equation} \label{eq.apost-banach}
   d ( v , v^\ast ) \leq (1-\kappa)^{-1} d(v , \op [v]),
\end{equation}
where $\kappa$ is the contraction rate.

The a-posteriori formulation consists in controlling derivatives of $v
- v ^\ast $ by the low regularity norm of the initial error. If the
initial error is small, this formulation assures that there is a true
solution close to such initial guess in the sense of $C^j$. The
estimation is done using interpolation inequalities.

\subsubsection{A-posteriori argument on a bounded interval}

On an interval $[a,b]$, the inequality \eqref{eq.apost-banach} implies
\begin{equation*}
    \| (v - v^\ast)|_{[a,b]} \| _{C ^0} \leq e ^{\delta \eta} (1-\kappa)^{-1}d(v , \op [v]),
\end{equation*}
where $\delta = \max\{|a| , |b|\}$. Thus, by using the interpolation
inequalities in \eqref{eq.cellineq}, we deduce that there are
constants $\ct c _c$, $\ct c _s$, and $\ct c _u$ such that
\begin{equation*}
\begin{split}
    \| (\vft _{(0)} - \Pi ^c [v ^\ast]) |_{[a,b]} \| _{C ^j} &\leq
    \ct c _c e ^{\delta \eta \frac{\ell +1-j}{\ell+1}} (1 - \kappa) ^{- \frac{\ell +1-j}{\ell+1}}  d(v , \op [v]) ^{ \frac{\ell +1-j}{\ell+1}} \qquad 0 \leq j \leq \ell, \\
    \| (\hx ^s _{(0)} - \Pi ^s [v ^\ast]) |_{[a,b]} \| _{C ^j} &\leq
    \ct c _s e ^{\delta \eta \frac{\ell +2-j}{\ell+2}} (1 - \kappa) ^{- \frac{\ell +2-j}{\ell+2}}  d(v , \op [v]) ^{ \frac{\ell +2-j}{\ell+2}} \qquad 0 \leq j \leq \ell+1, \\
    \| (\hx ^u _{(0)} - \Pi ^u [v ^\ast]) |_{[a,b]} \| _{C ^j} &\leq
    \ct c _u e ^{\delta \eta \frac{\ell +2-j}{\ell+2}} (1 - \kappa) ^{- \frac{\ell +2-j}{\ell+2}}  d(v , \op [v]) ^{ \frac{\ell +2-j}{\ell+2}} \qquad 0 \leq j \leq \ell+1,
\end{split}
\end{equation*}
where $\Pi ^\sigma [w](t) \bydef \Pi ^\sigma _t w$ for $\sigma \in
\{c,s,u\}$. Note that for the stable and unstable directions, we could
use the interpolation with $C^1$ and $C^{\ell+1+\Lip}$ spaces as well.

\subsubsection{A-posteriori argument on semi lines}
Let $g \colon \R \to \R ^n$ be a smooth function and define $g _\eta \colon \R \to \R ^n$ as 
\[
 g _\eta(t) \bydef e ^{-\eta |t|} g(t) = 
 \begin{cases}
     e ^{-\eta t} g(t) & t > 0 \\
     e ^{\eta t} g(t) & t \leq 0.
 \end{cases}
\]
In general, the function $g _\eta$ is not differentiable at $t = 0$
for $\eta > 0$. Therefore, we provide interpolation inequalities for
$t > 0$ in $C _\eta$ space (recall Section~\ref{sec.contraspace}) in
the following Lemma~\ref{lem.cetaineq}. Similar results hold for $t <
0$.

\begin{lem} \label{lem.cetaineq}
 Let $g \colon (0,+\infty) \to \R ^n$ be a $C ^{\ell+\Lip}$ function ($\ell \geq 0$). Then 
 \[
  \| g \| _{C _\eta} \leq 1 \quad \text{implies} \quad \| \dop ^j g\| _{C _\eta} \leq \ct c _j \| g \| _{C _\eta} ^{\frac{1}{j+1}},
 \]
 for all $j = 0, \dotsc, \ell$ and some constants $\ct c _j$'s
 depending on $\eta$, $j$, and $\|g\|_{C^{j+1}}$ (
 $\|g\|_{C^{\ell+\Lip}}$ when $j=\ell$).
 
\end{lem}
\begin{proof}
  Let us prove the result by induction: 
  \begin{enumerate}
  \renewcommand{\theenumi}{\roman{enumi}}
  \renewcommand{\labelenumi}{\theenumi.)}
      \item For $j = 1$, we use the interpolation inequality
        \eqref{eq.cellineq}. Noticing that $\| g _\eta \| _{C^0} = \|
        g \| _{C _\eta}$ and $\|g \| _{C _\eta} \leq \|g \| _{C _\eta}
        ^{1/2}$ as $\|g \| _{C _\eta} \leq 1$, we have
        \begin{equation*}
          \begin{split} 
          | e^{-\eta t} \dop g(t) | &= | \dop g _\eta (t) + \eta e ^{-\eta t} g(t) | \leq \| g _\eta \| _{C ^1} + \eta \|g \| _{C _\eta} 
          \leq  M _{0,2} \|g _\eta \| _{C ^0} ^{1/2} \| g _\eta \| _{C ^2} ^{1/2} + \eta \| g _\eta \| _{C^0}  \\
          &=\bigl( M _{0,2} \| g _\eta \| _{C ^2} ^{1/2} + \eta \bigr) \| g _\eta \| _{C^0}^{1/2}.
      \end{split} 
      \end{equation*}
      We take $\ct c _1 = M _{0,2} \| g _\eta \| _{C ^2} ^{1/2} +
      \eta$, so that the proof of the case $j = 1$ is done.
      \item Assume that the result is true up to $j-1$. By Leibnitz product formula, we have 
      \[
       \dop ^j g _\eta(t) = \sum _{k = 0}^j \binom{j}{k} \dop ^k(e ^{-\eta t}) \dop ^{j-k}g (t) = e ^{-\eta t} \biggl(\dop ^{j}g (t) + \sum _{k = 1}^j \binom{j}{k} (-\eta) ^k \dop ^{j-k}g (t)\biggr).
      \]
      Then by induction hypotheses and interpolation inequality,
      \begin{multline*}
          | e^{-\eta t} \dop ^j g(t) | = \biggl| \dop ^j g _\eta + \sum _{k = 1}^j \binom{j}{k} (-\eta) ^{k+1} \dop ^{j-k}g (t
          )\biggr| \leq \| g _\eta \| _{C ^j} + \sum _{k = 1} ^j \ct a _k \|g _\eta \| _{C ^0}^{\frac{1}{k}}  
            \\
          \leq 
          \biggl( M _{0,j+1} \| g _\eta \| _{C ^{j+1}} ^{\frac{j}{j+1}} + \sum _{k = 1} ^j \ct a _k \|g _\eta\| _{C ^0}^{\frac{1}{k}- \frac{1}{j+1}} \biggr) \| g _\eta \| _{C^0} ^{\frac{1}{j+1}},
      \end{multline*}
      for some $\ct a _k$ involving combinatorial numbers. We let $ c
      _j = M _{0,j+1} \| g _\eta \| _{C ^{j+1}} ^{\frac{j}{j+1}} +
      \sum _{k = 1} ^j \ct a _k$ so that the result for $j$ is
      proved. \qedhere
  \end{enumerate}
\end{proof}

Then from \eqref{eq.apost-banach} and Lemma~\ref{lem.cetaineq}, we conclude that for a good enough initial guess, 
\begin{equation*}
\begin{split}
    \| \dop ^j (\vft _{(0)} - \Pi ^c [v ^\ast]) |_{(0,+\infty)}  \| _{C _\eta} &\leq
    \ct c _c (1 - \kappa) ^{-\frac{1}{j+1}}  d(v , \op [v]) ^{\frac{1}{j+1}} \qquad 0 \leq j \le \ell, \\
    \| \dop ^j (\hx ^s _{(0)} - \Pi ^s [v ^\ast]) |_{(0,+\infty)} \| _{C _\eta} &\leq
    \ct c _s (1 - \kappa) ^{-\frac{1}{j+1}}  d(v , \op [v]) ^{\frac{1}{j+1}} \qquad 0 \leq j \le \ell+1, \\
    \| \dop ^j(\hx ^u _{(0)} - \Pi ^u [v ^\ast])  |_{(0,+\infty)} \| _{C _\eta} &\leq
    \ct c _u (1 - \kappa) ^{-\frac{1}{j+1}}  d(v , \op [v]) ^{\frac{1}{j+1}} \qquad 0 \leq j \le \ell+1,
\end{split}
\end{equation*}
where $\Pi ^\sigma [w](t) \bydef \Pi ^\sigma _t w$ for $\sigma \in
\{c,s,u\}$. Similar results hold on $(-\infty,0)$.

\section{Further results}\label{sec.further}
In this section, we discuss bootstrap of the regularity and non-autonomous unperturbed  systems. 

\subsection{Estimates on the growth of higher derivatives} \label{sec.exposol}

For ODEs, one could bootstraps the regularity of the solution: An
initial value problem of an ODE, say $\dot y(t) = g\circ y(t)$ with $y
(0) = y _0$, has the property that if one is able to find a $C ^1$
solution and $g$ is $C ^{\ell+1}$, then automatically such a solution
will be $C ^{\ell +2}$ for $\ell \geq 0$.  If we considered $\Pop
\colon C ^{\ell+1} \to C ^{\ell+1}$ in Theorem~\ref{main}, we would
have the same bootstrap property as in ODE's and we could first find
solution in $C ^1$ space. However, this setting would not cover
applications with neutral or small delays.

Instead, we consider $\Pop\colon C ^{\ell+1} \to C
^{\ell}$. Therefore, we are not able to bootstrap regularity directly
in this case.  Hence, the fixed point method should be performed on a
suitable space up to the right regularity level (beyond $C ^1$), see
Section~\ref{sec.space}. Nevertheless, once the Theorem~\ref{main} is
proved, we can bootstrap other type of solution behaviors; solutions
with exponential derivative growth.

\smallskip

We stress that we are looking for $C^1$ solutions such that higher
derivatives can arbitrarily grow. A simple example is the function
$z(t) = \int _0 ^t \sin e^s \, ds$, which is $C ^1$ but from the
second derivative on grows exponentially. Because we will have $C ^1$
solutions, then the perturbation $\pp$ will at least be $C ^0$
(otherwise it would not be possible to be controlled by the
  perturbative parameter $\ep$).

Let us now deduce how adding some slightly different assumptions
  to $\pp$ we can include new type of solutions. Indeed, given a $C ^1$
solution $\xx(t)$ of \eqref{model}, if we consider the second and
third derivatives, then we have
\begin{equation*}
 \begin{split}
 \ddot \xx(t) &= \dop \vf \circ \xx(t) \dot \xx(t) + \ep \frac{d}{dt}\pp(t,\xx _t, \ep, \prm), \\
 \dddot \xx(t) &= \dop ^2 \vf \circ \xx(t) \dot \xx(t) ^{\otimes 2} + \dop \vf \circ \xx(t) \ddot \xx(t) + \ep \frac{d^2}{dt^2}\pp(t,\xx _t, \ep, \prm).  
 \end{split}
\end{equation*}
If $\ddot \xx$ has an exponential growth, then it necessarily comes
from the perturbation $\pp$ since $\dop \vf$ and $\dot \xx$ are
bounded. More precisely, if there are $\gamma \geq 0$ and  $C
_j >0$ such that
\begin{equation*}
 \biggl|\frac{\partial ^j}{\partial t ^j}\pp(t,\xx _t, \ep, \prm)
 \biggr| \leq C _j e ^{\gamma j|t|},
\end{equation*}
then for $j=1, \dotsc, \ell$, the solution $\ddot \xx(t)$ will also be
bounded exponentially and, in general, what we have is that $|\dop
^{j+1} \xx (t) | e ^{-\gamma j|t|} < +\infty $ for $j \geq 0$.

\smallskip

To provide a formal statement, let us define the exponential
derivative growth space:
\begin{dfn}[Finitely differentiable space with exponential growth] 
 Let $\gamma \geq 0$ and let $C ^\ell _\gamma(I, \R ^n)$ be the space
 of $\ell$ times differentiable functions on the interior of interval
 $I \subset \R$ and with finite norm:
 \begin{equation*}
  \| g \| _ {C ^\ell_\gamma} \bydef \max _{j = 0, \dotsc, \ell} \sup
  _{t \in I} | \dop ^j g (t) | e ^{-\gamma j |t|} \qquad \text{ for
    all } g \in C ^\ell _\gamma.
 \end{equation*}
\end{dfn}

Notice that when $\gamma > 0$, an element in $C^\ell_\gamma$ does not
have Lipschitz boundedness in all its
derivatives. Corollary~\ref{thm.expmain} is a bootstrap result in the
space with exponential growth defined above.

\begin{cor}\label{thm.expmain}
Let $\gamma \geq 0$, let $\xx(t)$ be a solution from
Theorem~\ref{main} for $\ell =1$, and let $\ell' \geq \ell$.

Assume that the unperturbed system satisfies:
\begin{enumerate}
\renewcommand{\theenumi}{H\textsubscript{$0,\gamma$}\arabic{enumi}}
\renewcommand{\labelenumi}{\sl \theenumi)} 
 \item The hypothesis \eqref{H02} in Theorem~\ref{main} holds for
   $\ell'$,
\end{enumerate} 

and that the perturbation $\Pop$ in Theorem~\ref{main} satisfies:
\begin{enumerate}
\renewcommand{\theenumi}{H\textsubscript{$\ep,\gamma$}\arabic{enumi}}
\renewcommand{\labelenumi}{\sl \theenumi)}
\setcounter{enumi}{0}
 \item \label{Hep1exp} For all $\ep \in (0,\ep _0)$, $t \in \R$ and
   for $j = 1, \dotsc, \ell '$, $u \in C ^1(\R, \R ^n)$, $\dop u \in
   C ^{j-1}_\gamma(\R, \R ^n)$,
 \begin{equation*}
  \biggl| \frac{d ^j}{dt ^j} \Pop[u,\ep,\prm](t) \biggr| \leq C _j
  e^{\gamma j|t|} F (\| u \| _{C ^1}, \| \dop  u \| _{C ^{j-1}_\gamma}),
 \end{equation*}
 where $C _j > 0$ and $F:\R\to\R$ is  increasing and continuous.
\end{enumerate} 
 Then the $C ^1$-solution $\xx(t)$ of \eqref{model} is such that $\dop
 \xx$ is in $C ^{\ell'}_\gamma$ .
\end{cor}

Notice that Corollary~\ref{thm.expmain} does not change the range of
$\ep\leq \ep _0$ and it can incorporate smooth parameter dependence
with the approach introduced in Section \ref{sec.smoothprm} using
$C^{\ell}_\gamma$ spaces.

\subsection{Non-autonomous unperturbed system}

Our set up can incorporate non-autonomous systems using the standard
method of adding an extra variable. Consider a non-autonomous system
$\dot x(t) = g(x(t), t)$, where $g$ is $\ell$-times differentiable and
Lipschitz. We introduce a new variable $s$, and let $y \bydef (x,t)$,
then
\begin{equation*}\label{nonaut}
 y'(s) = \frac{d}{ds} 
 \begin{pmatrix}
  x(s) \\ t(s)
 \end{pmatrix}= G\circ y(s) =
 \begin{pmatrix}
  g(x(s), t(s)) \\
  1
 \end{pmatrix}.
\end{equation*}

Let us define an affine differentiable space $\widetilde C^\ell \bydef
\Id + C ^\ell$. This space has a well-defined Lipschitz constant. 
Thus, a solution $y$ belongs to the product space $C^{\ell+1} \times
\widetilde C ^{\ell+1}$. 

\begin{rmk} 
Note that our setting for hyperbolic orbits does not involve the orbit to lie
on a bounded set, it only requires that the vector field are 
bounded in a neighborhood of uniform size of the orbit. In the
non-autonomous case, this amounts to uniform for all the derivatives
of small enough order of $g$ -- including derivatives with respect to
time -- in a neighborhood of uniform size of the orbit.  Hence, we can
remake all the unperturbed hypothesis admitting these affine
differentiable spaces and derive a similar result like in
Theorem~\ref{main} that explicitly includes non-autonomous unperturbed
systems.

In the applications to delay equations, we will include, for
technical reasons that the delays are bounded. 
\end{rmk}

\smallskip

Even if this very direct approach gives results for many applications,
it can be improved. Indeed, It is well known \cite{Obaya} that one can
obtain a theory of evolutions of the equation $\dot x(t) = g(x(t), t)$
by assuming only that $g$ is measurable with respect to $t$ (several
mild integrability assumptions are needed). This is usually called
\emph{Caratheodory theory}. Under rather mild assumptions, the
Caratheodory theory allows to write variational equations and the
remainder. The operator $\op$ in this paper can then be formulated
just as well. At the moment, we are not aware of any significant
applications.

\section{Some models covered by the general results}\label{sec.examples}

This section is devoted to providing examples of perturbations $\pp$
which satisfy the assumptions of our main theorem. We show how to
verify the hypotheses in Theorem~\ref{main} and we add some important
remarks.

\subsection{ODE Perturbation}

A very particular case of Theorem~\ref{main} is when $P(t, x_t) = g(t,
x(t) )$, where the history segment $x_t$ is evaluated at zero to
obtain $x(t)$. This case corresponds to ODE perturbations.

When there is a hyperbolic orbit in the unperturbed system satisfying
Definition \ref{dfn.unihypersol}, we obtain that there is a solution
close to the unperturbed hyperbolic orbit applying
Theorem~\ref{main}. The hyperbolicity of the perturbed solution can be
seen from \cite{Moser69}. This is a version of Anosov shadowing
theorem \cite{Anosov69}. The precise version is close to the version
in \cite{Moser69} as modified in \cite{LlaveMM86}.

We show that hyperbolic orbits have a counterpart in the perturbed
system. As a corollary of our formalism (as in \cite{LlaveMM86}), we
obtain smooth dependence on parameters, see
Section~\ref{sec.smoothprm}.

Note that the range of perturbation parameters for which the orbit
persists depends on the hyperbolicity parameters of the orbit.  Also
the size affected by perturbations on an orbit depends on the
hyperbolicity parameters. For Anosov systems for which all the orbits
have uniform hyperbolicity constants, the validity range of
perturbations is uniform and the size of the perturbation effects is
uniform.  In non-uniformly hyperbolic sets, the allowed values of the
perturbation and the size of the responses will depend a lot on the
orbits.

\subsection{State and time dependent delay equations} \label{sec.sdde}
Let us consider the model
\begin{equation} \label{sddemodel}
 \dot x(t) = \vf \circ x(t) + \ep Q \bigl(t, x(t + \rr (t, x(t))) \bigr).
\end{equation}
where the perturbative map $\pp$ is 
\begin{equation*} \label{sddepp}
 \pp (t, \vartheta, \ep) \bydef Q \bigl(t, \vartheta \circ \rr (t, \vartheta (0)) \bigr).
\end{equation*}
The perturbative hypotheses in Theorem~\ref{main} are satisfied by
considering $\rr \colon \R \times \R^n \to \R$ a $C ^{\ell+\Lip}$ map
and $Q \colon \R \times \R ^n \to \R ^n$ a $C ^{\ell+\Lip}$ map. In this
case, the history segment can be taken as $h \bydef \| \rr \| _{C
  ^0}$.

The hypothesis \eqref{Hep1} can be verified using the chain rule and
Fa\'a di Bruno formula. To check hypothesis \eqref{Hep2}, we analyze
$\Pop$ in Theorem~\ref{main}:
\begin{equation*} \label{Popsdde}
    \Pop[u,\ep](t) =  \Pop[u](t) = Q\bigl(t, u(t + \rr (t, u(t)))\bigl) ,
\end{equation*}
where we have omitted the $\ep$ and $\prm$ in $\Pop$ since in this
example $Q$ does not depend on them. Hence,
\begin{equation} \label{Popsdde.Hep2}
    |\Pop[u ^2](s) - \Pop[u ^1](t)| \leq \Lip (Q) |s - t| + \Lip (Q) | u ^2(s + \rr(s, u ^2(s)) - u ^1(t + \rr(t, u ^1(t))|.
\end{equation}
The second term is bounded by adding/subtracting and triangle inequality, that is,
\begin{equation*}
\begin{split}
    | u ^2(s + \rr(s, u ^2(s)) - u ^1(t + \rr(t, u ^1(t)) | &\leq 
    | u ^2(s + \rr(s, u ^2(s)) - u ^2(s + \rr(t, u ^1(t)) | \\
    &\qquad +    | u ^2(s + \rr(t, u ^1(t)) - u ^1(t + \rr(t, u ^1(t)) |\\ &\leq
    \| u ^2 \| _{C ^1} \| \rr \| _{C ^1} \bigl[ |s - t| + \| u ^2 _s - u ^1 _t \| _{C ^{0}[-h,h]}\bigr]  \\
    & \qquad + \| u ^2 _s - u ^1 _t \| _{C ^{0}[-h,h]}.
\end{split}
\end{equation*}
Then we can take constants $\ct L _1$ and $\ct L _2$ so that \eqref{Hep2} is true for all $u ^1, u ^2$ in a ball of $C ^{\ell + 1 +\Lip}(\R, \R ^n)$.

\subsection{Nested delay equations} \label{sec.nested}

Let us consider a differential equation with nested delay/advance terms
\begin{equation*} \label{nestedsddemodel}
 \dot x(t) = \vf \circ x(t) + \ep Q \bigl(t, x(t + \rr (t, x(t + \rr _1 \circ x(t)))) \bigr).
\end{equation*}
In this case, the perturbative map for \eqref{model} is
\begin{equation*} 
 \pp (t, \vartheta, \ep) \bydef Q(t, \vartheta \circ \rr (t, \rr_1 \circ \vartheta (0))),
\end{equation*}
and the ``history segment'' is $h \bydef \max\{\| \rr \| _{C ^0}, \|
\rr _1\| _{C ^0} \}$. If $\rr \colon \R \times \R^n \to \R$ is a $C
^{\ell +\Lip}$ map, $\rr _1 \colon \R ^n \to \R ^n$ a $C ^{\ell+\Lip}$
map, and $Q \colon \R \times \R ^n \to \R ^n$ a $C ^{\ell+\Lip}$ map, then
the perturbative hypotheses in Theorem~\ref{main} are satisfied.

The idea is similar to Section~\ref{sec.sdde}. Now we need to bound 
\begin{equation} \label{Popnestedsddee.keyterm}
    | u ^2(s + \rr(s, u ^2(s + \rr _1 \circ u ^2(s))) - u ^1(t + \rr(t, u ^1(t + \rr _1 \circ u ^1(t)))|.
\end{equation}
We obtain that
\begin{equation*}
    \begin{split}
        \eqref{Popnestedsddee.keyterm} &\leq 
        | u ^2(s + \rr(s, u ^2(s + \rr _1 \circ u ^2(s))) - u ^2(s + \rr(t, u ^1(t + \rr _1 \circ u ^1(t)))| \\
        & \qquad + 
        | u ^2(s + \rr(t, u ^1(t + \rr _1 \circ u ^1(t)) - u ^1(t + \rr(t, u ^1(t + \rr _1 \circ u ^1(t))| \\ &\leq 
        \| u ^2 \| _{C ^1} \| \rr \| _{C^1} \bigl[|s - t| + \| u ^2 \| _{C ^1} \| \rr _1 \| _{C^1} \| u ^2 _s - u ^1 _t\| _{C ^0[-h,h]} + \| u ^2 _s - u ^1 _t\| _{C ^0[-h,h]} \bigr] \\
        & \qquad +  \| u ^2 _s - u ^1 _t\| _{C ^0[-h,h]}.
    \end{split}
\end{equation*}
Therefore, for all $u ^1, u ^2$ in a ball of $C ^{\ell + 1 +\Lip}(\R,
\R ^n)$, there are constants $\ct L _1$ and $\ct L _2$ such that
\eqref{Hep2} is satisfied.

\subsection{Neutral delay equations} \label{sec.neutral}

As an example of neutral delay/advance equation, we consider 
\begin{equation} \label{neutralddemodel}
  \dot x(t) = \vf \circ x(t) + \ep Q \bigl(t, x(t + \rr (t, \tfrac{d }{dt}x(t))) \bigr),
\end{equation}
Where $Q: \R \times \R^n$ is a smooth function. 

This can be made into the form \eqref{model} taking. 
\begin{equation*} 
 \pp (t, \vartheta, \ep, \prm) \bydef Q\bigl(t, \vartheta \circ \rr (t, \tfrac{d}{ds}\vartheta(0))\bigr),
\end{equation*}
that depends on time and on the derivative of the state. Note that we
used the fact that $\frac{d x _t}{ds}(0) = \frac{d x}{dt}(t) $. The
history segment in this case is $h \bydef \| \rr \| _{C^0}$.
Note that we are not assuming that that the sign of $r$ is
negative, so that we can just as well have advanced equations. 

Using the standard adding of extra variables, the unperturbed equation
could be an equation of order $n+1$, but the R.H.S cannot introduce
derivatives of order higher than $n+1$, 

To apply Theorem~\ref{main}, we assume regularities on $Q$ and $\rr$
such that the perturbation in \eqref{neutralddemodel} satisfies
\eqref{Hep1}. In particular, if $\rr \colon \R \times \R^n \to \R$ is
a $C ^{\ell + \Lip}$ map and $Q \colon \R \times \R ^n \to \R ^n$ is a
$C ^{\ell+\Lip}$ map, \eqref{Hep1} is verified.

To check \eqref{Hep2}, we bound the term
\begin{equation} \label{Popneutral.Hep2.keyterm}
    \Bigl| u ^2\bigl(s + \rr(s,  \tfrac{d}{ds }u ^2(s)\bigr) - u ^1\bigl(t + \rr(t, \tfrac{d}{dt} u ^1(t))\bigr) \Bigr|,
\end{equation}
and obtain
\begin{equation*}
    \begin{split}
    \eqref{Popneutral.Hep2.keyterm} &\leq 
    \Bigl| u ^2\bigl(s + \rr(s,  \tfrac{d}{ds }u ^2(s)\bigr) - u ^2\bigl(s + \rr(t, \tfrac{d}{dt} u ^1(t))\bigr) \Bigr| \\
    & \qquad + 
    \Bigl| u ^2\bigl(s + \rr(t, \tfrac{d}{dt} u ^1(t))\bigr) - u ^1\bigl(t + \rr(t, \tfrac{d}{dt} u ^1(t))\bigr) \Bigr| \\ &\leq
    \| u ^2 \| _{C ^1} \| \rr \| _{C ^1} \bigl[ |s - t| + \| u ^2 _s - u ^1 _t \| _{C ^{1}[-h,h]}\bigr]  \\
    & \qquad + \| u ^2 _s - u ^1 _t \| _{C ^{0}[-h,h]}.
    \end{split}
\end{equation*}
Therefore, \eqref{Hep2} is satisfied.

The modification of the verification for several delays/advances is
left to the reader. Note that we can let some of the $r$'s be delays and others
be advances. 

\smallskip

Similar to Section~\ref{sec.nested}, one can also consider nested
delays involving first derivative in the state, or more generally
$\frac{d}{ds} \vartheta(s)$ for $s \in [-h,h]$. Of course, particular
cases such as constant delays satisfy the assumptions of our result.

\subsection{Small delays and small advances}
\label{sec:small} 
There are problems in the literature in which the time at which the
solution needs to be evaluated contains very small time changes. An
important case, which involves special challenges, is the motion of
charged particles, studied in more detail in
Section~\ref{sec:electrodynamics}.

In this section we will show that the terms with small delay or small
advances can be included in the formalism of Theorem~\ref{main}.  The
delays allowed are very general and could be functionals on the
history segment. Only some mild regularity assumptions will be
imposed. In particular, we do not need to assume that the delays
are positive, so we can also consider advanced perturbations (or
perturbations that include both advanced and retarded terms).  This
generality becomes useful in the treatment of motion of point charges
where the delay can depend on the whole trajectory.  Some of the
physical theories proposed in \cite{WheelerF49,WheelerF45} involve both advanced
and retarded term and, hence, they could be included in our framework.

A small delay is a singular perturbation because the nature of the
problem changes completely.  Heuristically, the expansions on the
perturbative parameter, involve derivatives of the function. If the
delay is not zero -- even if small -- the phase space may be an
infinite dimensional solution manifold or something more complicated.

\medskip

The simplest non-trivial case is 
\begin{equation}\label{smalldelay} 
  \dot x(t) =  f \circ x(t - \ep\tau( x_t) ),
\end{equation}
where $\tau$ is a functional of the trajectory segment.
This case \eqref{smalldelay}  fits into our framework 
by rewriting it as
\begin{equation*}\label{smalldelay2} 
  \dot x(t) =  f \circ x(t)  +
  \ep \biggl[ \frac{1}{\ep} f\circ x(t - \ep\tau( x_t) )- \frac{1}{\ep}f\circ x(t) \biggr].
\end{equation*}
To estimate the perturbation and verify \eqref{Hep1}--\eqref{Hep2} we use 
the heuristic idea that 
\[
  \frac{1}{\ep} f\circ x(t -\ep\tau( x_t) )- \frac{1}{\ep}f\circ x(t)
\approx -\dop f \circ x(t) x'(t) \tau(x_t).
\]
Therefore, the effect of the small delay is similar to
including a functional losing one derivative, which is incorporated in Theorem~\ref{main} fortunately. 

More precisely: 
\[  
 \frac{1}{\ep} f\circ x(t - \ep\tau( x_t) )- \frac{1}{\ep}f\circ x(t) = -\int_0^1 \dop f\circ x(t - \sigma \ep \tau(x _t)) x'(t-\sigma \ep\tau(x _t)) \tau (x _t) \, d\sigma.
\]
Hence, equation \eqref{smalldelay} with small delays or advances fits our setting where
the functional perturbative map is 
\begin{equation} \label{Qop}
    \mathcal{Q}[\vartheta,\ep] \bydef -\int_0^1 \dop f\circ \vartheta(-\sigma \ep \tau\circ \vartheta) \vartheta'(-\sigma \ep\tau\circ \vartheta) \tau \circ \vartheta \, d\sigma,
\end{equation}
where $\tau \colon C ^{\ell+1+\Lip}([-h,h],\R^n) \to \R$. Note that,
we can also consider several delays $\tau _i$ in \eqref{smalldelay}.

Applying Theorem~\ref{main}, we obtain
that
\begin{thm}\label{thm.smalldelay}
  Consider the equation
  \begin{equation} \label{modelsmall}
    \dot x(t) = \vf \circ \Bigl( x \bigl(t-\ep \tau_1(t, x_t) \bigr), \dotsc, x \bigl(t -\ep\tau_L(t, x_t)  \bigr) \Bigr)
    + \ep \pp(t, x_t, \ep, \prm).
  \end{equation}
  Assume \eqref{H01}--\eqref{H02} in Theorem~\ref{main} and that for all $i=1,\dotsc, L$, if $x _t$
  is in a ball in $C ^{\ell + 1+\Lip}$ space, then $\tau _i(t,x _t)$ ranges
  in a ball in $C^{\ell+\Lip}$ space, and that $\tau_i$'s have
  Lipschitz properties.

  Then, hypotheses \eqref{Hep1} and \eqref{Hep2} hold for the $\mathcal{Q}_i$ given as 
\[
\mathcal{Q}_i(t,\vartheta, \ep)\bydef\int _{0}^1 \dop _i \vf \Bigl( \vartheta \bigl( -\sigma \ep \tau _1 (t, \vartheta) \bigr), \dotsc, \vartheta \bigl(-\sigma \ep \tau _L (t,\vartheta) \bigr) \Bigr) \vartheta' \bigl(-\sigma \ep \tau _i(t,\vartheta) \bigr)\tau _i (t, \vartheta)
\, d\sigma,\]

Assume in addition that, $\Pop$ defined in \eqref{eq.Px}
  corresponding to $\pp$ in equation \eqref{modelsmall} satisfies  \eqref{Hep1} and \eqref{Hep2} in
  Theorem~\ref{main}, then we have the same conclusions as Theorem~\ref{main}.
\end{thm}

Maybe the most unclear part of the proof of
Theorem~\ref{thm.smalldelay} is that $\mathcal{Q}_i$ satisfies \eqref{Hep2}
(since \eqref{Hep1} comes from the fact that $\tau _i$ maps a ball in
$C^{\ell+1+\Lip}$ to another ball in $C^{\ell+\Lip}$). We illustrate
\eqref{Hep2} with $\mathcal{Q}$ in \eqref{Qop}. Let $\vartheta$ and
$\varrho$ in a $C ^{\ell +1 + \Lip}$ ball and let us bound
$\mathcal{Q}[\vartheta,\ep] - \mathcal{Q}[\varrho,\ep]$ whose
integrand is (after adding and subtracting)
\begin{align} 
&\dop f\circ \vartheta(-\sigma \ep \tau\circ \vartheta) \vartheta'(-\sigma \ep\tau\circ \vartheta) \tau \circ \vartheta - \dop f\circ \varrho(-\sigma \ep \tau\circ \varrho) \varrho'(-\sigma \ep\tau\circ \varrho) \tau \circ \varrho \nonumber \\ 
=&  \bigl[\dop f\circ \vartheta(-\sigma \ep \tau\circ \vartheta) - \dop f\circ \varrho(-\sigma \ep \tau\circ \vartheta)\bigr] \vartheta'(-\sigma \ep\tau\circ \vartheta) \tau \circ \vartheta \tag{\texttt{S1}} \label{S1} \\
&+\bigr[\dop f\circ \varrho(-\sigma \ep \tau\circ \vartheta) - \dop f\circ \varrho(-\sigma \ep \tau\circ \varrho)\bigr] \vartheta'(-\sigma \ep\tau\circ \vartheta) \tau \circ \vartheta \tag{\texttt{S2}}\label{S2} \\
&+\dop f\circ \varrho(-\sigma \ep \tau\circ \varrho)\bigl[ \vartheta'(-\sigma \ep\tau\circ \vartheta) - \varrho'(-\sigma \ep\tau\circ \vartheta) \bigr]\tau \circ \vartheta \tag{\texttt{S3}} \label{S3}\\
&+\dop f\circ \varrho(-\sigma \ep \tau\circ \varrho)\bigl[ \varrho'(-\sigma \ep\tau\circ \vartheta) - \varrho'(-\sigma \ep\tau\circ \varrho) \bigr]\tau \circ \vartheta \tag{\texttt{S4}}\label{S4}  \\ 
&+ \dop f\circ \varrho(-\sigma \ep \tau\circ \varrho) \varrho'(-\sigma \ep\tau\circ \varrho) \bigl[\tau \circ \vartheta - \tau \circ \varrho\bigr].  \tag{\texttt{S5}}\label{S5}
\end{align}
The individual intermediate lines admit straightforward bounds in
terms of the inputs, i.e. for $s \in [-h,h]$,
\begin{align*}
 |\eqref{S1}(s)| &\leq  \| \vf \| _{C ^2} \| \vartheta \| _{C ^1} \| \tau \| _{C ^0} \| \vartheta - \varrho \| _{C ^0}\\
 |\eqref{S2}(s)| &\leq \sigma \ep \Lip (\dop \vf \circ \varrho) \Lip (\tau) \| \vartheta - \varrho \| _{C ^0} \| \vartheta \| _{C ^1} \| \tau \| _{C ^0}\\
 |\eqref{S3}(s)| &\leq \| \vf \| _{C ^1} \| \vartheta - \varrho \| _{C ^1} \| \tau \| _{C ^0} \\
 |\eqref{S4}(s)| &\leq \sigma \ep \| \vf \| _{C ^1} \| \varrho \| _{C ^2} \Lip (\tau) \| \vartheta - \varrho \| _{C ^0} \| \tau \| _ {C^0} \\
 |\eqref{S5}(s)| &\leq  \| \vf \| _{C^1} \| \varrho \| _{C ^1} \Lip (\tau) \| \vartheta - \varrho \| _{C^0}.
\end{align*}
Thus, adding up all these bounds, there is a constant, say $\ct L_2$, such that $|\mathcal{Q}[\vartheta,\ep] - \mathcal{Q}[\varrho,\ep]| \leq \ct L _2 \| \vartheta - \varrho \| _{C^1}$.

\subsection{Delays implicitly defined by the solution. Applications 
  to electrodynamics}
\label{sec:electrodynamics}
 
In this section, we show that our framework applies to the problem of
electrodynamics of point charges and formulate
Theorem~\ref{thm:electrodynamics} which is obtained by applying
Theorem~\ref{main} to \eqref{electrodynamics}, the model of particles
moving on the electromagnetic fields generated by others.

From the mathematical point of view, the main new problem is that the
delays that appear in the equation -- the time the signals emitted by
one particle take to reach another -- depend on the trajectory. The
delays can only be found by solving an implicit equation that involves
the whole trajectory, see \eqref{delaydefined}. We refer to this
situation as implicitly defined delays.

\subsubsection{Formulation of the mathematical problem}
The motion of point charges under the electromagnetic field generated
by others has several physical problems due to \emph{self-energy} (see
\cite{Spohn04} and also \cite[Chapter 16]{Jackson}) which we will not
discuss.

We will follow the formulation of on \cite{WheelerF49} which avoids
the self-energy problems.  The basic idea of
\cite{WheelerF45,WheelerF49} is that each charge moves in the field
generated by the others and by external sources (not in the field
generated by themselves!).  The expression of the electromagnetic
fields generated by charges in motion is obtained by solving Maxwell
equations. The explicit expression of the solution of Maxwell
equations generated by charges (knowing their positions and
velocities) is well known since the turn of the XX century and is
called Li\'enard-Weichert potentials (\cite{LandauL,Rohrlich, Jackson,
  Zangwill}).  The motion of the particle in this potential is given
by Newton's laws taking the relativistic expression of the mass.  One
can think of the Li\'enard-Wiechert potentials as the standard
Coulomb/Ampere expressions taking into account delays and Fitzgerald
contractions.  As in the wave equation, the solutions of the Maxwell
equations can be advanced or retarded, or convex combinations of both.
In Physics literature, it is customary to take only retarded
solutions, but this restriction does not follow from Maxwell equations
or the boundary conditions and we do not need it for the results in
this paper.

Hence, the equation of the $i$th particle are of the form: 
\begin{equation}\label{electrodynamics} 
  \ddot q_i(t) = A_{\text{ext}}(t, q_i(t), \dot q_i(t)) + \sum_{j \ne i}
  A_{i,j}\bigl(
  q_i(t), \dot q_i(t), 
  q_j(t - \tau_{ij}), \dot q_j(t - \tau_{ij}), 
  q_j(t + \sigma_{ij}), \dot q_j(t + \sigma_{ij})\bigr),
\end{equation}
The $q_i$ represents the position of the $i$ point charge, The term
$A_{\text{ext}}$ denotes the external force, and $A _{i,j}$ is an
explicit expression given by the Li\'enard-Weichert potentials that
depend on the time delay/advance defined by solving the implicit
equations.
\begin{equation}\label{delaydefined}
  \begin{split}
    & \tau_{ij}(t)  = \frac{1}{c} \big| q_i(t) -   q_j(t -\tau_{ij}(t) ) \big| ,\\
    & \sigma_{ij}(t)  = \frac{1}{c} \big| q_i(t) -   q_j(t + \sigma_{ij}(t) ) \big| ,\\
    \end{split}
    \end{equation} 
with $c$ being the speed of light.  We think of $c$ as large so that
$\ep = \frac{1}{c}$ is a small parameter.

The most salient mathematical feature of \eqref{electrodynamics} is
that it involves delays (or advances) which correspond to the time
that the light takes to travel from the source particle to the dynamic
particle.  This delay depends on the whole trajectory of both
particles (one needs to solve implicitly for the trajectory of a light
ray to intersect the trajectory of the source particle).  Since the
Li\'enard-Wiechert potentials can be retarded or advanced (or convex
combinations of both)\footnote{ In \cite{WheelerF49}, it is suggested
that combining the delay and advance with coefficients $1/2$ is
physically important.}  we get that the resulting equations can be
retarded or advanced also.

Some minor complications are that \eqref{electrodynamics} presents some
singularities when $q_i(t) = q_j(t); i \ne j$ or when $|\dot q_i(t)| =
c$.

The explicit form of the expressions of $A _{i,j}$ and $A
_{\text{ext}}$ in \eqref{electrodynamics} can be found in any advanced
textbook and the detailed expression is not relevant for our
treatment.

Our treatment is rather general and applies to other models of the
same structure. Due to relativity, all models of particles interacting
by pairwise interactions are of the form \eqref{electrodynamics}.  The
structure \eqref{electrodynamics} includes models incorporating
gravity or more manageable approximations of Li\'enard-Wiechert
potentials (it is common to keep the first
order in $\ep$ and ignore higher orders in $\ep$ such as Fitzgerald
contractions).

The treatment presented here extends to modifications of
\eqref{electrodynamics} that involve interactions among 3 or more
bodies.
\begin{multline}\label{electrodynamics2} 
  \ddot q_i(t) = A_{\text{ext}}(t, q_i(t), \dot q_i(t)) +
  F_i \bigl( \{ q_j(t -\tau_{jk}(t))\}_{j,k=1}^N, 
    \{ \dot q_j(t -\tau_{jk}(t))\}_{j,k=1}^N, \\
    \{ q_j(t + \sigma_{jk}(t))\}_{j,k=1}^N, 
        \{ \dot q_j(t + \sigma_{jk}(t))\}_{j,k=1}^N \bigr) ,
\end{multline}
where $F$ is an expression with the same type of singularities and the
delays/advances $\tau_{jk}$, $\sigma_{jk}$ are defined in
\eqref{delaydefined} or even by more general procedures that 
involve all the trajectories.  The only requirement is that the 
delays range in a  $C ^{\ell + \Lip}$ ball if the trajectories 
range in a $C^{\ell+ \Lip } $
 ball.

If there are no external electromagnetic forces (or if the external
electromagnetic forces are time independent), the classical model
conserves energy, so that it does not have any hyperbolic orbits. On
the other hand, in time dependent (e.g. periodic) external fields,
many interesting models (e.g. ion traps, mirror magnet machines) have
many hyperbolic orbits \cite{Ghosh96,Kajita22}.

In this paper, we show that near hyperbolic orbits of the unperturbed model
\eqref{electrodynamics}, i.e. $\ep = 0$ in \eqref{delaydefined}
\footnote{For $\ep = 0$, the delays and advances vanish so that 
\eqref{electrodynamics} is an ODE.}, and if the perturbed model avoids
the singularities and it has bounded delays, Theorem~\ref{main}
applies. Moreover, constructed solutions are similar to those
hyperbolic orbits for the unperturbed model \eqref{electrodynamics},
as long as $0 < \ep \ll 1$ in \eqref{delaydefined}. See
Theorem~\ref{thm:electrodynamics} for a precise formulation.

\medskip

Before proceeding to the detailed analysis, let us make some comments
on the equations and their physical properties.

\begin{enumerate}
\item
  The expressions defining the forces are algebraic expressions
  (arithmetic operations and square roots). They have singularities
  when there are collisions ($q_i(t) = q_j(t)$ for some $i\ne j$) or
  when a particle reaches the speed of light ($|\dot q_i(t) | = c$ for
  some $i$).

  We will assume that the unperturbed solutions we consider are a
  finite distance away from these singularities so that the
  expressions given the second derivatives and given the positions and
  velocities are smooth functions in a neighborhood of the unperturbed
  solution.  See Definition~\ref{nonsingular}.
  
\item 
  The delays $\tau_{ij}$ or advances $\sigma_{ij}$ as in
  \eqref{delaydefined} involve solving implicit equations that involve
  the pairwise trajectories $q _i(t)$ and $q _j(t)$.

  The physical meaning of the delays/advances are the time it takes a
  light signal to travel from the particle $i$ to the particle
  $j$. Since the speed of the particles is bounded away from the speed
  of light, this time exists and is unique. Nevertheless, finding the
  delay requires to solve an equation that depends on the pairwise
  trajectories.  See \eqref{delaydefined}.

  Fortunately, for the method used in this paper, the main property
  needed is that $\tau_{ij}$ and $\sigma_{ij}$ are uniformly smooth
  assuming that the trajectories are smooth (and that they are away
  from the singularities).
  
  Note that in general $\tau_{ij} \ne \tau_{ji} $ even if $\tau_{ij} -
  \tau_{ji} = O(\ep ^2)$. Indeed, after expanding the solution of
  \eqref{delaydefined} up to first order, i.e.
  \begin{equation}\label{appoxdelay} 
   \tau _{ij}(t) = \ep | q _i (t) - q _j(t) | + \ep ^2 \bigl(q _{i}(t) - q _j(t)\bigr) \cdot \dot q _j(t) + O (\ep ^3).
  \end{equation}
 (hint: It is easy to consider the expansions for $\tau_{ij}^2 $ obtained by squaring both sides of
\eqref{delaydefined} and express the square of length as an inner product). 
The above derivation \eqref{appoxdelay}  of an approximate form of the delay is  purely formal. 
Note that it is only valid for differentiable  $q$ 
and that the error incurred in the approximation depends on higher  derivatives of $q$. 
Of course in a  delay equation, modifying the form of the delay is a very singular 
perturbation so, substituting the approximation above may lead to equations with different solutions.   

In Theorem~\ref{main} the main objects are sets of uniformly differentiable $q$. 
In these sets of uniformly differentiable functions, the approximations in 
\eqref{appoxdelay} are uniform (in a norm involving one less derivative 
than the uniform).

\item
Theorem~\ref{main} also applies to many modifications of 
the model and produces solutions. For example, ignoring The Figtgerald contractions \cite{Ver2} or changing 
the delay by its state dependent approximation \eqref{appoxdelay} \cite{Driver} (they 
are formally second order in $\ep$). The exact solutions of these models obtained 
applying Theorem~\ref{main} 
will be approximate solutions of \eqref{electrodynamics}. 

In theoretical Physics, there are also other methods to produce
approximate solutions using formal power series
\cite{CasalCL20,GimenoLY25,MS} or using numerical approximations.

The a-posteriori formulation of Theorem~\ref{main} shows that these
approximate solutions are close  to solutions of \eqref{electrodynamics}.

\item
  The delays/advances $\tau_{ij}, \sigma_{ij} $ may be unbounded if
  the particles $i,j$ get far apart. This is not covered by our
  Theorem~\ref{main}. We assume in this paper that the
  unperturbed orbits remain in a bounded region.

  One can hope that this assumption can be weakened because when the
  delays are large, the interactions are weak.  The problem of charges
  scattering is, of course, of great physical importance and has been
  considered many times.  Starting with the pioneering work of
  \cite{DriverElec}.  We refer to the recent \cite{BauerDDH17}
  for an account of progress in this line of research and
  several other approximate models of
  electrodynamics. 
\item
  Other recent advances in the theory of constructing solutions
  for electrodynamics exploit
that periodic solutions 
  of these equations have a variational structure. This allows
  to use deep variational tools such as Floer theory.
  We refer to \cite{AlbersFS20,Fraunfelder21}.

  The approach for periodic orbits in \cite{YangGL22}
    does not require a variational structure, but on the other hand,
    requires proximity to an ODE. 
\end{enumerate}

Denoting $y(t) \bydef (q_1(t),\dotsc, q_N(t), \dot q_1(t),\dotsc, \dot
q_N(t) )$, we can write the equation \eqref{electrodynamics} in the
form of \eqref{smalldelay} with the delays being implicitly
defined. Note that there are $N(N-1)$ delays and $N(N-1)$ advances.
Of course, we can consider cases where the expression of the forces
does not depend on the advances.

Note that formally, the effect of delays is $O(1/c)$ whereas the
effect of Lorentz-Fitzgerald contractions is $O(1/c^2)$. Therefore,
it is common to consider models in which only the delays are
considered and the Lorentz-Fitzgerald contraction effects are ignored
\cite{Ver2, ChiconeM07}.

These models are of the form considered in \eqref{electrodynamics}. So, the results of this paper on persistence of hyperbolic solutions
of the non-relativistic models apply.  Furthermore, thanks to the
a-posteriori formulation of Theorem~\ref{main} we obtain that the
exact solutions in these models are at a distance $O(1/c^2)$ from the
solutions to the full model.

Similarly, the hyperbolic solutions produced by other models that
solve the relativistic equations to order $n$ (e.g. \cite{MS} to order
$4$) will be $O(\ep^{n+1})$ close to the solutions produced here.

However, it is important to note that the justification of 
the results here is only for hyperbolic solutions and that the
quantitative aspects of the corrections needed
depend on the hyperbolicity constants.  In
particular, in non-uniformly hyperbolic sets, where the hyperbolic
constants deteriorate, the range of validity of
the approximations
will become smaller and the errors will be affected by larger
constants.  This is consistent with the impossibility of formulating
effective equations valid everywhere \cite{CurrieS}, which still
allows formulating approximations in some sets.

\subsubsection{The result}
\label{sec:mathematicalresult}

Since the equation~\eqref{electrodynamics} has singularities, we have
to assume that the unperturbed solution is away from the singularities.

\begin{dfn}
  \label{nonsingular}
  We say that a solution of the classical equations of motion is
  non-singular when there exist $0<\xi_1 <1$ and $\xi_2>0$, such that
  for all $t$:
\begin{equation*} 
\label{trajectory_condition}
\begin{split} 
& |\dot q_j(t)| \le \xi_1 c, \qquad \text{for all } j,  \\ 
& |q_i(t) - q_j(t) | \ge \xi_2, \qquad \text{for all } i \ne j.
\end{split} 
\end{equation*}
\end{dfn}

The internal forces and the
masses are analytic around non-singular solutions, by Definition~\ref{nonsingular}.  Therefore, the regularity
assumptions of Theorem~\ref{main} on the unperturbed equation concern only the
external fields.

Assuming that the solutions remain in a bounded set and that the
trajectories are uniformly away from collisions, we obtain that the vector
fields giving the evolution are uniformly differentiable.

\begin{thm}
\label{thm:electrodynamics} 
Consider the model \eqref{electrodynamics} with the delays or advances defined in
\eqref{delaydefined}. Denote $1/c$ as $\ep$, and treat it as a
parameter.

Assume that the external fields are 
$A _{\mathrm{ext}}$ are $C^{\ell + 2+\Lip}$.

Assume that for $\ep = 0$, the ODE \eqref{electrodynamics} has a
solution such that:
\begin{enumerate}
\item
  It is hyperbolic in the sense of Definition~\ref{dfn.unihypersol},
\item
  It is non singular in the sense of Definition~\ref{nonsingular}.
\item
  It lies in a bounded set.
\end{enumerate}

Then, Theorem~\ref{main} applies to the problem given by
\eqref{electrodynamics} and for small enough values of $\ep$, we can
find solutions of \eqref{electrodynamics} of the form in \eqref{sol}.
\end{thm} 

The proof follows from Theorem~\ref{thm.smalldelay} once we have the
estimates on regularity bounds and low regularity contraction of the
delays and advances \eqref{delaydefined}. These estimates will be
obtained in the following section.

Given the formulation of the fixed point argument, we only need to
prove that the delays are differentiable when the history segments are
assumed to be differentiable (propagated bounds) and to show that the
operator $\Gamma$ is a contraction in a low regularity space when the
history segments are differentiable.

\medskip

\paragraph{Results on the regularity properties
  of the delay}

In this section, we study \eqref{delaydefined} as an equation for
$\tau_{ij}(t)$ and $\sigma_{ij}(t)$ when we prescribe the trajectories $q_i$ and $q_j$.
This makes precise the notion that the delay/advance is a functional of the
whole trajectory and it shows that 
Theorem~\ref{thm:electrodynamics} follows from
Theorem~\ref{thm.smalldelay}.

\begin{pro} \label{prop:tauestimates}
Let $q_i$ and $q_j$ be continuously differentiable trajectories that
satisfy Definition~\ref{trajectory_condition} and have 
derivatives bounded uniformly away from the speed of light $c$.

Then, for all $t\in \R$, we can find unique $\tau_{ij}(t),
\sigma_{ij}(t) > 0 $ solving \eqref{delaydefined}.

Furthermore, let $\tilde \tau_{ij}, \tilde \sigma_{ij}$ be the delay and advance of trajectories $\tilde q_i$ and $\tilde q_j$, there exists a
constant $C$ such that
\[
\| \tau_{ij} - \tilde \tau_{ij}\|_{C^0}, 
\| \sigma_{ij} - \tilde \sigma_{ij}\|_{C^0} \le C \bigl( \| q _i - \tilde q _i\| _{C ^0} + \| q _j - \tilde q _j\| _{C ^0} \bigr).
\]

Moreover, if the trajectories $q_i$ and $q_j$ are $C^{\ell+\Lip}$ and
satisfy Definition~\ref{nonsingular}. Then the $\tau_{ij}, \sigma_{ij}$ are $C^{\ell
  + \Lip}$, and there is an explicit algebraic expression $g$ such
that
\begin{equation*} 
\label{taurbounds} 
\| \tau_{ij} \|_{C^{\ell + \Lip}},
\| \sigma_{ij} \|_{C^{\ell + \Lip}} 
\le g( \| q_i\|_{C^{\ell + \Lip}},  \| q_j\|_{C^{\ell + \Lip}}, \xi_1, \xi_2).
\end{equation*}

\end{pro}

\begin{proof}  
The first part of Proposition~\ref{prop:tauestimates} follows from the
standard contraction mapping theorem applied to \eqref{delaydefined}.
The second part also follows from the contraction principle,
remembering that we are assuming uniform differentiability of the
$q_i$.

\smallskip

Let us define the operator $\mathcal{N}[\tau; q_i, q _j](t) \bydef \ep
|q _i(t) - q _j(t - \tau(t))|$ and let us first prove that
$\mathcal{N}$ is a contraction for small enough $\ep$. That is,
\begin{multline*}
  |\mathcal{N}[\tau](t) - \mathcal{N}[\widetilde \tau](t)| = 
  \ep \bigl| |q _i(t) - q _j(t - \tau(t))| - |q _i(t) - q _j(t - \widetilde \tau(t))|\bigr| \\
  \leq \ep |q _j(t-\tau(t)) - q _j(t - \widetilde \tau(t))| 
  \leq \ep \Lip(q _j) \|\tau - \widetilde \tau\| _{C ^0}.
\end{multline*}
Let $\ep < 1 / \| q _j\|_{C^1}$, and define $\kappa
\bydef \ep \| q _j\|_{C^1}$, then $\mathcal{N}$ is a contraction with rate $\kappa$.

\smallskip

We first bound $\|\mathcal{N}[\tilde \tau; q_i,
  q _j] - \tilde \tau \| _{C^0} \leq B$ for $\tilde \tau$ a
fixed point of $\mathcal{N}[\cdot; \tilde q _i, \tilde q _j]$
of some particles $\tilde q _i$ and $\tilde q _j$, and $B$ related to the difference between $q$ and $\tilde q$. This
implies $\| \tau - \tilde \tau \| _{C ^0} \leq
\frac{B}{1-\kappa}$.

\begin{align*}
 |\mathcal{N}[\tilde\tau; q_i, q _j](t) - \tilde \tau(t) | 
 = \bigl| \ep |q _i(t) &- q _j(t - \tilde \tau(t))|  - \ep |\tilde q _i(t) - \tilde q _j(t - \tilde \tau(t))| \bigr | \\
 &\leq \ep | q _i(t) - \tilde q _i(t) | 
 + \ep | q _j(t - \tilde \tau(t)) - \tilde q _j(t - \tilde \tau(t))| \\
 &\leq \ep \bigl( \| q _i - \tilde q _i\| _{C ^0} + \| q _j - \tilde q _j\| _{C ^0} \bigr).
\end{align*}

The estimates on the derivatives of $\tau_{ij}, \sigma_{ij}$ follow from applying the
implicit function theorem for the solutions of \eqref{delaydefined}.
\end{proof}

By Proposition \ref{prop:tauestimates}, for $\tau_{ij}, \sigma_{ij}$, we have  Lipschitz estimates and that they are  $C^{\ell+\Lip}$ if  $q_i, q_j$ are
$C^{\ell+\Lip}$. Therefore,
 we can use the procedure in
Section~\ref{sec:small} and apply Theorem \ref{thm.smalldelay}. Note that the $\tau $'s in
Theorem~\ref{thm.smalldelay} correspond to $\tau _{ij}/\ep$ or $\sigma
_{ij}/\ep$ since both $\tau _{ij},\sigma _{ij}$ are of order $O(\ep)$.

\appendix

\section*{Acknowledgments}

The project has been supported with the Spanish grant
PID2021-125535NB-I00 (MICINN/AEI/FEDER, UE)
and the Catalan grant 2021 SGR 01072.  The project that gave
rise to these results also received the support of the fellowship from
``la Caixa'' Foundation (ID 100010434), the fellowship code is
LCF/BQ/PR23/11980047. J.G. also thanks the School of Mathematics of GT
for its hospitality in 2022 year. This project has been supported by the Fundamental Research Funds for the Central Universities.

\section*{Statements and Declarations}
The authors declare that they have no conflict of interest.

\bibliographystyle{alpha}
\addcontentsline{toc}{section}{References}
{\small \bibliography{ref2,hyperbolic} }

\end{document}